\documentclass[final,3p,times, num,sort&compress,]{elsarticle}

\usepackage{amssymb}
\usepackage{amsmath}
\usepackage{multibib}
\usepackage{mathtools}
\usepackage{caption}
\usepackage{subcaption}
\usepackage{float}
\usepackage{color,soul}
\usepackage[document]{ragged2e}
\usepackage{microtype} 
\usepackage[T1]{fontenc}
\numberwithin{equation}{section}
\graphicspath{{Results/}}
\usepackage{soul}
\usepackage[usenames,dvipsnames,svgnames,table]{xcolor}

 \usepackage{lineno}

\journal{Computational Physics}
\newcommand{\dx}{{\Delta{x}}}
\newcommand{\dt}{{\Delta{t}}}

\newcommand{\RN}[1]{%
  \textup{\uppercase\expandafter{\romannumeral#1}}%
}

\begin{document}
\begin{frontmatter}



\title{Numerical Artifacts in the Generalized Porous Medium Equation: Why Harmonic Averaging Itself Is Not to Blame}
\author[label1]{Danielle C. Maddix}
\author[label2]{Luiz Sampaio}
\author[label1,label2]{Margot Gerritsen}
\address[label1]{Institute of Computational and Mathematical Engineering}
\address[label2]{Energy Resources Engineering \\Stanford University}
\justify

\begin{abstract} 
The degenerate parabolic Generalized Porous Medium Equation (GPME) poses numerical challenges due to self-sharpening and its sharp corner solutions.  For these problems, we show results for two subclasses of the GPME with differentiable $k(p)$ with respect to $p$, namely the Porous Medium Equation (PME) and the superslow diffusion equation.  Spurious temporal oscillations, and nonphysical locking and lagging have been reported in the literature.  These issues have been attributed to harmonic averaging of the coefficient $k(p)$ for small $p$, and arithmetic averaging has been suggested as an alternative.  We show that harmonic averaging is not solely responsible and that an improved discretization can mitigate these issues.  Here, we investigate the causes of these numerical artifacts using modified equation analysis.  The modified equation framework can be used for any type of discretization.  We show results for the second order finite volume method.  The observed problems with harmonic averaging can be traced to two leading error terms in its modified equation.  This is also illustrated numerically through a Modified Harmonic Method (MHM) that can locally modify the critical terms to remove the aforementioned numerical artifacts.  

\end{abstract}


\end{frontmatter}


\section{Introduction}
\justify
\label{intro}
We discuss how to discretize the self-sharpening Generalized Porous Medium Equation (GPME) without the temporal oscillations, the locking, and the lagging reported in previous works \cite{lipnikov2016, inl2008, nie2013, vandermeer2016}.  In these works, the numerical artifacts have been attributed to harmonic averaging, and arithmetic averaging has been proposed to resolve them.    The exact causes of the artifacts were not published in the literature.  In this paper, we aim to understand the exact causes.  We demonstrate that harmonic averaging is not solely to blame and that the artifacts also depend on both the spatial and temporal discretizations.  The critical terms are identified through modified equation analysis and numerical evidence is provided to show that counteracting these terms removes the artifacts.  We refer to this demonstration method as the Modified Harmonic Method (MHM). 

The GPME is given in its multi-dimensional, conservative form as, 
\begin{equation} 
			\begin{aligned}
				&p_t = \nabla \cdot (k(p) \nabla p)   \hspace{.1 cm} \text{ in } \hspace{.1 cm} \Omega, \\
				&p(x,0) = h(x), \\
				&p(x,t) = g(x,t), \text{ } \forall x \in \partial \Omega,
			\end{aligned}
			\label{eq:GPME}
		\end{equation}
where $\partial \Omega$ is the boundary of the domain $\Omega$,  $h(x)$ is the initial condition and $g(x,t)$ is the Dirichlet boundary condition. 
The GPME can be seen as a subset of the fully general equation class that we are interested in, where $k(p)$ in Eqn. \eqref{eq:GPME} is replaced by $k(p,x)$.  The fully general equation class includes the variable coefficient problem, which will be discussed in a subsequent paper \citet{maddix_applied}.

\subsection{Applications of the GPME}
The Porous Medium Equation (PME), a subclass of the GPME, is used in a number of applications.
 In the PME, 
 \begin{equation}
 	k(p) = p^m,
	\label{eq:PME}
\end{equation}
for $m \ge 1$.  
The name PME originates from modeling gas flow in a porous medium.  The PME can be derived from the continuity equation, Darcy's Law and the equation of state for perfect gases  \cite{vazquez2007, ngo2016}.  The polytropic exponent $m$ relates the pressure $k(p)$ to the density $p$.  The case $m = 1$ models isothermal processes, while $m > 1$ for adiabatic processes.  From experimental data, \citet{vazquez2007} found that $m = 1.405$ corresponds to airflow at normal temperatures.

Nonlinear heat transfer of plasmas (ionized gases), mainly by radiation, drove much of the initial theoretical research on the PME \cite{vazquez2007}.  \citet{zeldovich66} proposed a model of this heat transfer at very high temperatures given by Eqn. \eqref{eq:GPME}.  Here, the coefficient $k(p)$ is also a monomial of $p$ where in this case, $k(p)$ represents the radiation thermal conductivity and $p$ the temperature.  In the optically thick limit approximation, $m = 3$.  In more complex models with multiple ionized gases, this exponent can increase to $4.5-5.5$.  

The PME with $m = 1$ has several applications, including groundwater flow with Boussinesq's equation \cite{boussinesq1903} and population dynamics in biology  \cite{vazquez2007}.


An example of the GPME that we will also consider is superslow diffusion \cite{galaktionov2004, vazquez2007} with
			\begin{equation}
				k(p) = \exp(-1/p).
				\label{eq:superslow}
			\end{equation}
The name superslow diffusion was introduced because the diffusivity $k(p) \rightarrow 0$ as $p  \rightarrow 0$ faster than any power of $p$ \cite{galaktionov91}.   This equation is used to model the diffusion of solids at different absolute temperatures $p$.  The coefficient $k(p)$ now represents the mass diffusivity and is connected with the Arrhenius law in thermodynamics \cite{galaktionov2006}.  

In the above examples of the GPME, $k(p)$ is a continuous function with respect to $p$.  An example of the GPME with discontinuous $k(p)$ is a Stefan problem \cite{vazquez2007}, where $k(p)$ is a step function, e.g.
\begin{equation}
					k(p) = \begin{dcases}
  									1,           & \text{if} \hspace{.1cm} p \ge 0.5,  \\
    									0,              & \text{otherwise}.
						\end{dcases}
					\label{eq:foam}
\end{equation}
The GPME with discontinuous coefficients is the topic of a subsequent paper in \citet{maddix_foam}.
 
In this paper, we focus on the GPME with continuous $k(p)$.  In particular, we will use the PME and the superslow diffusion equation.  For the PME, we consider the physical values for $m$ up to 3, which is large enough to exhibit the numerical artifacts of interest.  

\subsection{Degeneracy of the GPME and its Numerical Challenges} 
\label{degen}
Theoretical properties and behaviors of the GPME have been studied in many works, dating back to the 1950s with \citet{barenblatt56, barenblatt58, oleinik58, kalasnikov67} and \citet{aronson69} to more recently with \citet{shmarev03, shmarev05}.  In \cite{barenblatt56, barenblatt58}, the Barenblatt-Prattle self-similar solution of the PME and its finite propagation property are derived.  Self-similar solutions of the GPME have also been found \cite{galaktionov83}.  \citet{vazquez2007, ngo2016} provide review papers and detailed references.  

\begin{figure}[H]
		\center
		\begin{subfigure}[H]{.42\textwidth}  
			\includegraphics[width =\textwidth]{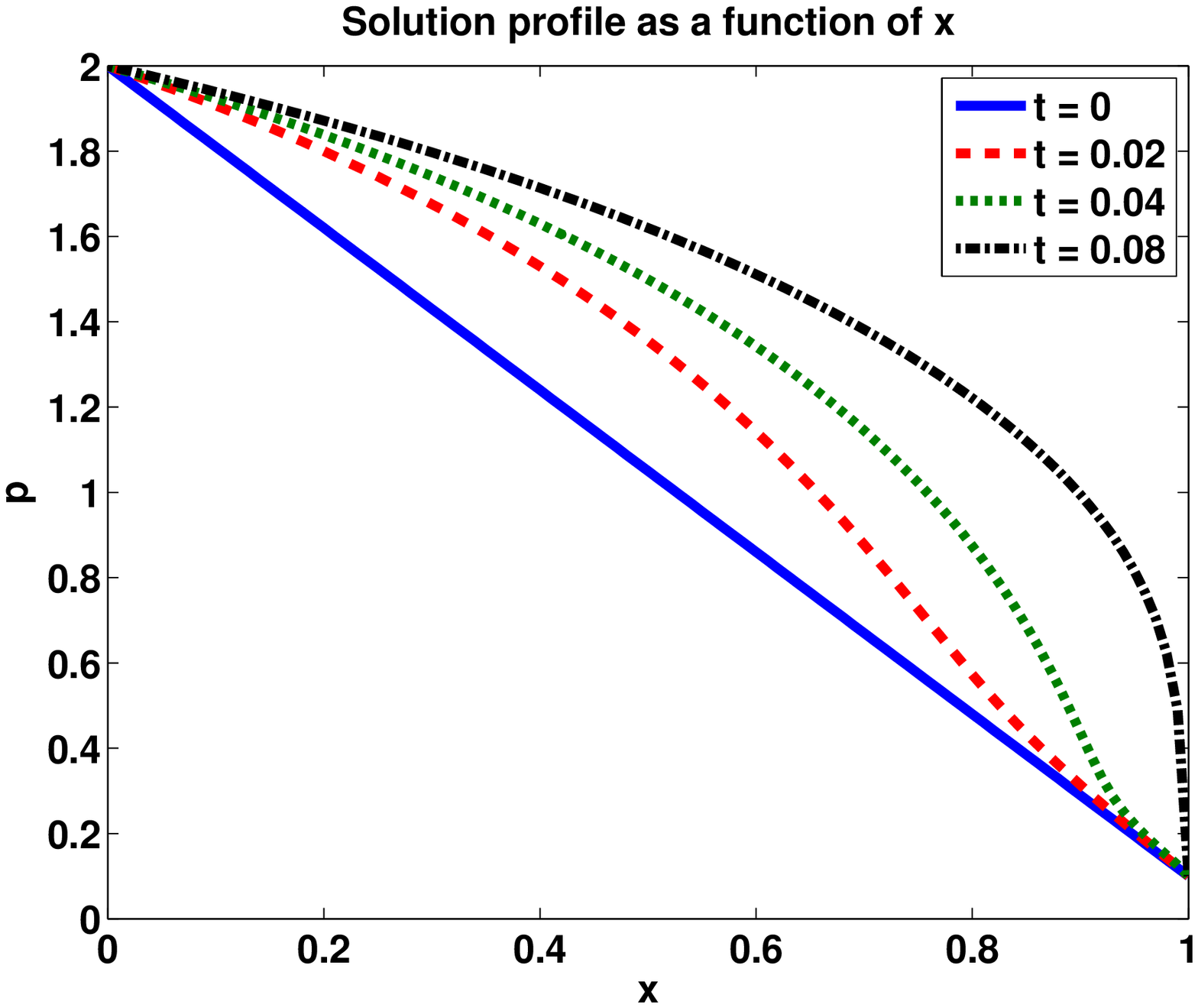}
			\caption {Self-sharpening example}
			\label{fig:selfsharp}
		\end{subfigure}
		\begin{subfigure}[H]{.42\textwidth}  
			\includegraphics[width =\textwidth]{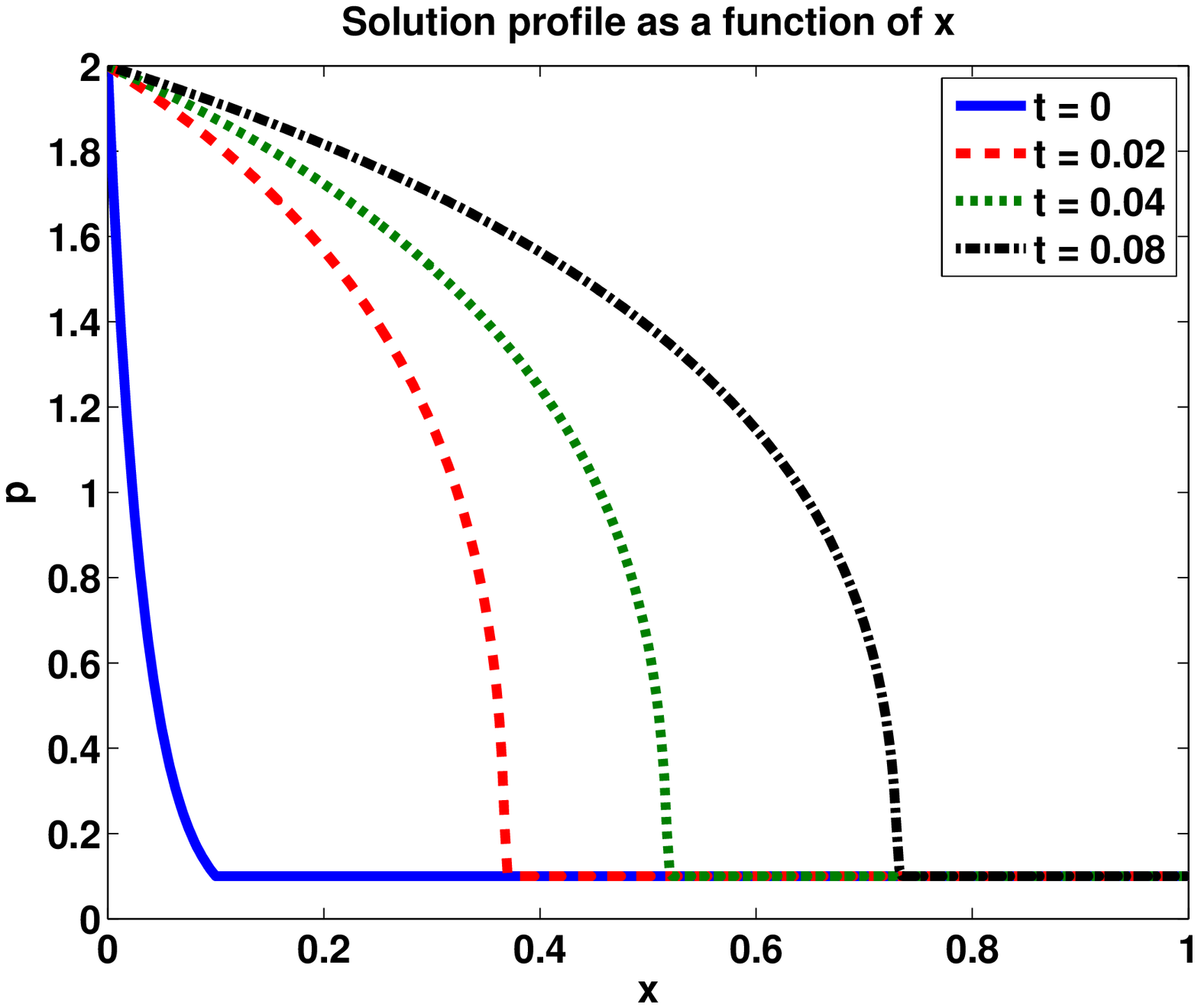}
			\caption {Moving interface example}
			\label{fig:movinginterface}
		\end{subfigure}
		\caption{Plots of reference solutions to the PME for $m\!\!=$ 3.  In Figure \ref{fig:selfsharp}, the initial condition leads to self-sharpening.  In Figure \ref{fig:movinginterface}, the initial condition leads to a moving interface.  The curve for $t = 0$ in Figure \ref{fig:movinginterface} will be used in Section \ref{MH_Results} for an initial condition in the numerical simulations.}
\end{figure}

Some of the behaviors of the GPME may be surprising since at first glance it appears to be a straightforward variation of the heat equation.  For example, self-sharpening can occur even with smooth initial data.  This self-sharpening effect is illustrated in Figure \ref{fig:selfsharp} for $k(p) = p^3$ $(m = 3)$, in which a linear initial condition develops a sharp gradient over time.  Figure \ref{fig:movinginterface} depicts a moving interface that can be developed for compactly supported initial conditions. The speed of the interface can be calculated exactly with Darcy's Law \cite{vazquez2007, ngo2016}.
In this solution, we see a sharp corner develop at the front.  Because of these behaviors, the GPME is referred to as degenerate parabolic.  
In \ref{degenerate}, we extend the PME results in \cite{vazquez2007} to the GPME and show that the governing equation \eqref{eq:GPME} can be expressed in terms of $k(p)$ as
\begin{equation}
							k_t = (1-C(p)) |\nabla k|^2 + k\Delta k,
	\label{eq:k_eqtn}
\end{equation}
where $C(p) = kk_{pp}/k_p^2$.
From Eqn. \eqref{eq:k_eqtn}, we see that parabolic degeneracy occurs when the following conditions are satisfied:
\begin{enumerate}
	\item $\lim_{p \to 0} k(p) = 0$ 
	\item $C(p) < 1$, $\forall p > 0$ and $\lim_{p \to 0} C(p) \le 1.$
\end{enumerate}
Condition 1 reveals that for small $p$, the coefficient $k(p)$ of the parabolic term $\Delta k$ tends to zero, resulting in an Eikonal equation for Eqn. \eqref{eq:k_eqtn} \cite{vazquez2007}.  By Condition 2, the coefficient of $|\nabla k|^2$ in the Eikonal equation is always positive for positive $p$. 
The faster the decay of $k(p)$ for small $p$, the more the parabolic character of the equation is masked.  We will see that this leads to challenges in the numerical results (Section \ref{MH_Results}).

The PME and superslow diffusion clearly satisfy both of the above conditions.  For the PME, $C(p) = 1 - \frac{1}{m} < 1$ for all $p$.   For superslow diffusion, it can be verified that 
\[
	\begin{aligned}
	&k_p^2 = k^2p^{-4}, \\
	&kk_{pp} = k^2(p^{-4} - 2p^{-3}) = k_p^2(1 - 2p),
	\end{aligned}
\]
and so $C(p) = 1 - 2p <  1$ for all $p > 0$.

\begin{figure}[H]
		\center
			\includegraphics[width =0.42\textwidth]{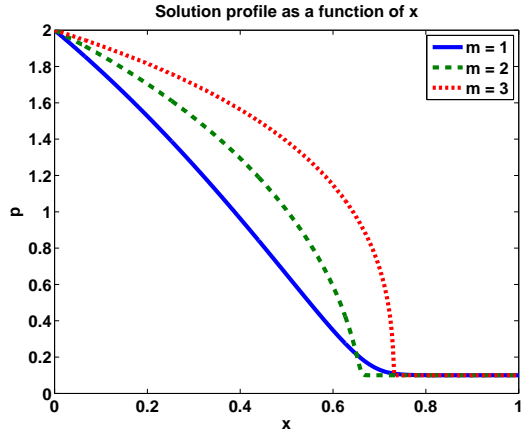}
	\caption{Comparison of reference solutions at time $t = 0.08$ to the PME for $m = 1,2,3$.}
	\label{fig:m}
\end{figure}

Apart from the degeneracy, there are distinct behaviors for the PME depending on the value of $m$, as illustrated in Figure \ref{fig:m}.  The sharp corners shown in Figure \ref{fig:movinginterface} do not necessarily form for lower values of $m$.  The gradient is finite for $m = 1$, and approaches infinity near the front for $m > 1$ \cite{vazquez2007}.


\subsection{Strategies to Overcome the Numerical Artifacts in Discretizations of the GPME} 
The degeneracy of the GPME poses interesting numerical challenges and, if care is not taken, locking, lagging and temporal oscillations can result \cite{lipnikov2016,inl2008, nie2013, vandermeer2016}.  
These numerical artifacts are distinct from the typical sharp gradient challenges caused by dispersive errors in hyperbolic conservation laws (see Section \ref{arith_harm}).
Locking, lagging and temporal oscillations occur for many types of discretizations, including the traditional Forward in Time, Central in Space (FTCS) discretization.  In this paper, we focus on understanding why these artifacts are occurring with the commonly used FTCS discretization and harmonic averaging.  We also apply a few implicit and higher-order temporal schemes.  We first discuss the main approaches suggested in the literature to alleviate the numerical artifacts.  

The first is to replace harmonic averaging of the coefficient with arithmetic averaging.  This works well for the continuous coefficient GPME, as illustrated in \cite{lipnikov2016,inl2008, nie2013}.  For applications in reservoir simulations, where the governing equation is parabolic with variable coefficients, arithmetic averaging can result in overly diffusive profiles and nonphysical solutions with flow allowed in and out of impermeable layers \cite{inl2008}.  Also, arithmetic averaging does not always remove the temporal oscillations, as seen in \cite{vandermeer2016, maddix_foam} for the GPME with discontinuous $k(p)$ (Eqn. \eqref{eq:foam}).  This discontinuous GPME is challenging and has led to several other approaches, e.g. the integral average in \cite{vandermeer2016}.


The second is Adaptive Mesh Refinement (AMR), where the mesh is refined in areas of steep gradients.  In \cite{maddix_foam}, it is shown that AMR can reduce these temporal artifacts but not remove them.

In the Finite Element community, two additional approaches have been proposed.  The first utilizes the Galerkin Finite Element Method with moving meshes, leading to grid refinement near the shock \cite{ngo2016}.  The numerical solution still has artifacts, including spatial oscillations near the shock resulting in the nonphysical loss of positivity of the solution.  \citet{zhang09} propose a high-order local discontinuous Galerkin (LDG) method on a uniform Cartesian mesh that removes any nonphysical spatial oscillations, but temporal oscillations are not discussed in this work.

An alternate approach not discussed in the literature that we examine is to change the temporal discretization.  An implicitly stable scheme and a total-variation diminishing (TVD) scheme are implemented to test their effect on the temporal oscillations.  The results in the Section \ref{TVD} show that the numerical artifacts remain with these choices of temporal schemes.  The Modified Equation Analysis in Section \ref{BE} can be used to explain the cause of the temporal oscillations of Backward Euler with a smaller time step and their absence with a larger time step.   
Our goal is to understand why the numerical artifacts occur with a traditional scheme (FTCS) and harmonic averaging.  
The differences, described in Section \ref{arith_harm}, between the schemes with arithmetic and harmonic averaging can be explained by their modified equations.  We present the Modified Equation Analysis in Section \ref{MEA}.  The Modified Harmonic Method (MHM) was developed in Section \ref{MH} to illustrate the impact of the modified equations.  
  The MHM can be seen as a Lax-Wendroff \cite{lax60} like approach that counteracts leading terms in the truncation error from the scheme.  We then show the results in Section \ref{MH_Results} that the MHM removes the numerical artifacts associated with harmonic averaging, confirming the modified equation analysis.  

\justify
\section{Numerical Discretization and Challenges} 

\label{arith_harm} 


For the reasons outlined above, we discretize with FTCS.  The domain $\Omega = [0,1]$ is partitioned into $N\!+\!1$ equally spaced cell-centered grid points $x_i$.  This defines a uniform Cartesian mesh with spatial step size $\dx = (x_{N+1} - x_{1})/N$, where 
\[ 
	0 = x_1 < x_1 + \dx <  \cdots < x_{1} + N\dx = x_{N+1} = 1.
\] 
There are $N\!-\!1$ control volumes (cells), referred to as $CV_i$ with cell faces that are at a distance of $\dx/2$ from the cell centers.  
The unknowns defined at each $x_i$ are $p_i$ for $i = 2, \dots, N$.  There are $N\!-\!1$ degrees of freedom, since $p_1$ and $p_{N+1}$ are fixed by the Dirichlet boundary conditions.  The unknowns $p_i$ and the corresponding coefficients $k_i$ are defined at the cell centers.  The fluxes 
\begin{equation}
	F_{i+1/2} \equiv k_{i+1/2}u_{i+1/2},
	\label{eq:fluxesdef}
\end{equation}
velocities\footnote{The convention for the velocity is not the standard Darcy velocity.} 
\begin{equation}
	u_{i+1/2} \equiv -\frac{p_{i+1}-p_i}{\dx},
	\label{eq:vel}
\end{equation} and numerical averages $k_{i+1/2}$ of $k_i$ and $k_{i+1}$ are defined at the $N$ cell faces.  We restrict $k_{i+1/2}$ to be a local average of its neighboring coefficients, as common in applications of interest.  Analogous definitions are used for $i\!-\!1/2$.  In subsequent subsections, we analyze and compare the arithmetic and harmonic averages. 


\subsection{Arithmetic and Harmonic Averages}
\label{arithharm_trunc}
The arithmetic average ($k^A_{i+1/2}$) of $k_i$ and $k_{i+1}$ on a uniform Cartesian grid is defined as
\begin{equation}
	k^A_{i+1/2} \equiv \frac{k_i + k_{i+1}}{2}.
		\label{eq:arith_avg}
	\end{equation}
Using Taylor Series expansions in space about $x_i + \dx/2$ for $k_i \equiv k(x_i, t)$ and $k_{i+1} \equiv k(x_i + \dx, t)$ at any time $t$, it can be shown that arithmetic averaging is second order accurate.  We define 
\[
	k^+ \equiv  k(x_i + \dx/2, t) \equiv k_{i+1/2}.
\]  
Let 
\[
	v =  k^+  +  \frac{\dx^2}{8}k_{xx}^++ \mathcal{O}(\dx^4)
\] be the even terms and 
\[
	w =   \frac{\dx}{2} k_x^+ + \mathcal{O}(\dx^3)
\]
be the odd terms in the Taylor expansion of $k_{i+1}$.  Then 
\[
	k^A_{i+1/2} = \frac{(v-w) + (v+w)}{2} = v = k^+ + \tau^A,
\] where
\begin{equation}
	\tau^A =  \frac{\dx^2}{8}k_{xx}^++ \mathcal{O}(\dx^4).
	\label{eq:arith_trun}
\end{equation}



The harmonic average ($k^H_{i+1/2}$) of $k_i$ and $k_{i+1}$ on a uniform Cartesian grid is defined as
					\begin{equation}
						k^H_{i+1/2} \equiv \frac{2k_i k_{i+1}}{k_i + k_{i+1}}.
						\label{eq:harm_avg}
					\end{equation}
The truncation error for the harmonic average can be written in terms of the truncation error for the arithmetic average in Eqn. \eqref{eq:arith_trun} as
\[
	k^H_{i+1/2} = \frac{v^2-w^2}{v} = v - \frac{w^2}{v} =  k^+  + \tau^H,
\] 
where
			\begin{equation}
				\tau^H =  \tau^A - \frac{\dx^2}{4} \frac{(k_x^+)^2}{k^+} + \mathcal{O}(\dx^4).
				\label{eq:harm_trun}
			\end{equation}	
The above equation is derived through substituting the expressions for $k_i = v-w$ and $k_{i+1} = v+w$ into Eqn. \eqref{eq:harm_avg}.

When there are smooth variations in $k$,  the two averages in Eqns. \eqref{eq:arith_avg} and \eqref{eq:harm_avg} are similar.  With sharp gradients in $k$, it is clear from Eqn. \eqref{eq:harm_trun} that the averages can significantly differ.  For example, if $k_{i+1} = \epsilon$, for some small $\epsilon > 0$, and $k_i > 1$, then $k_i + \epsilon \approx k_i,$ and 
	\begin{equation}
		k^H_{i+1/2} \approx \frac{2k_i\epsilon}{k_i} = 2\epsilon <<  \frac{k_i}{2} \approx k^A_{i+1/2}.
		\label{eq:compare}
	  \end{equation}
This shows that with a small coefficient at any one point, there can be orders of magnitude differences between the two averages.  
\subsection{Numerical Artifacts in the 1D PME}
\label{TLP}
The one-dimensional PME is an example where the harmonic average results in artificial and undesirable blocking.  The numerical solution is artificially damped, by not allowing enough flux across the moving interface.  The numerical artifacts encountered are lagging, temporal oscillations and locking.  

\begin{figure}[H]
		\center
		\begin{subfigure}[H]{.42\textwidth}  
			\includegraphics[width =\textwidth]{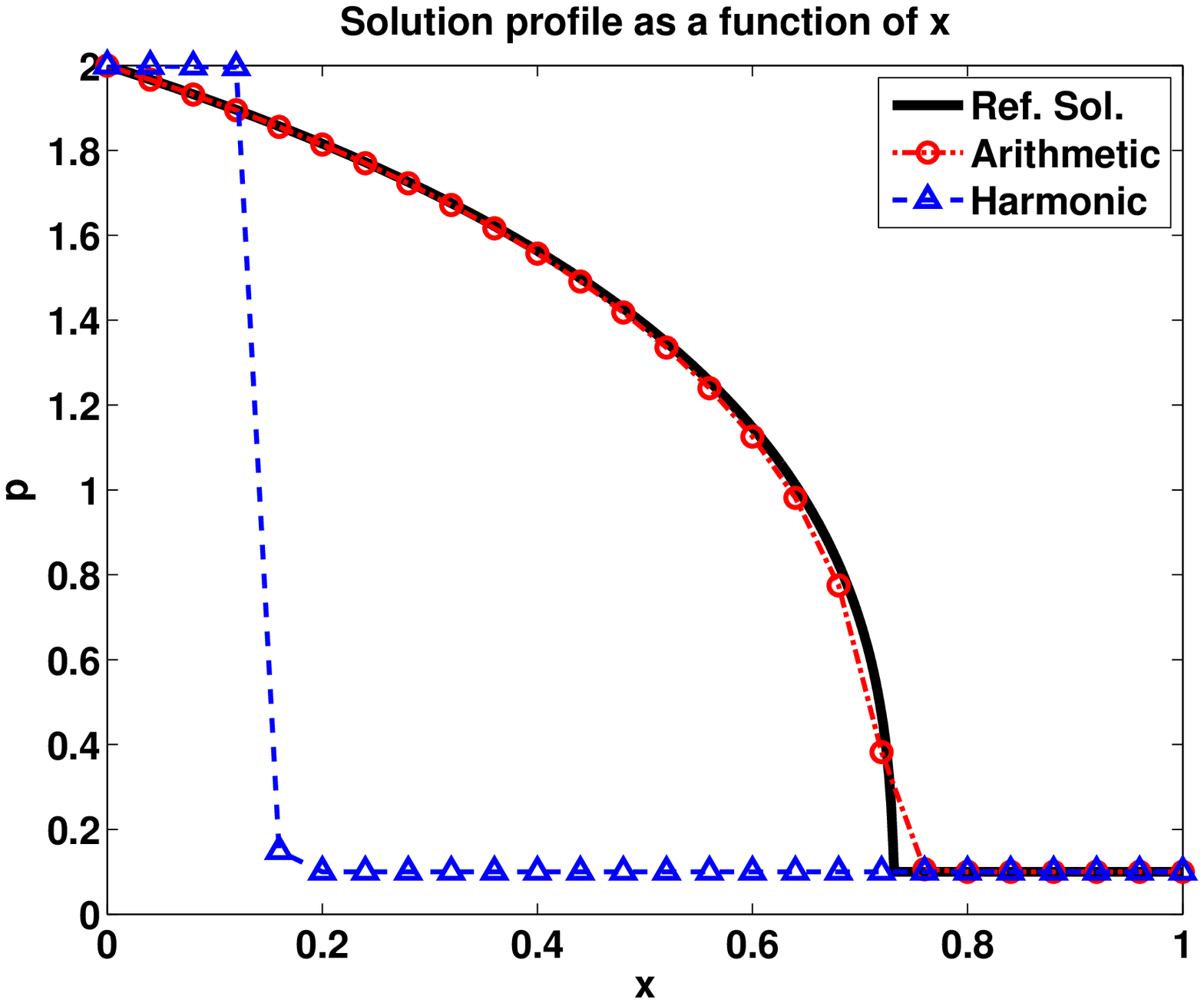}
			\caption {$t = 0.08$}
			\label{fig:pme_arithvsharm_space}
		\end{subfigure}
		\begin{subfigure}[H]{.42\textwidth}  
			\includegraphics[width =\textwidth]{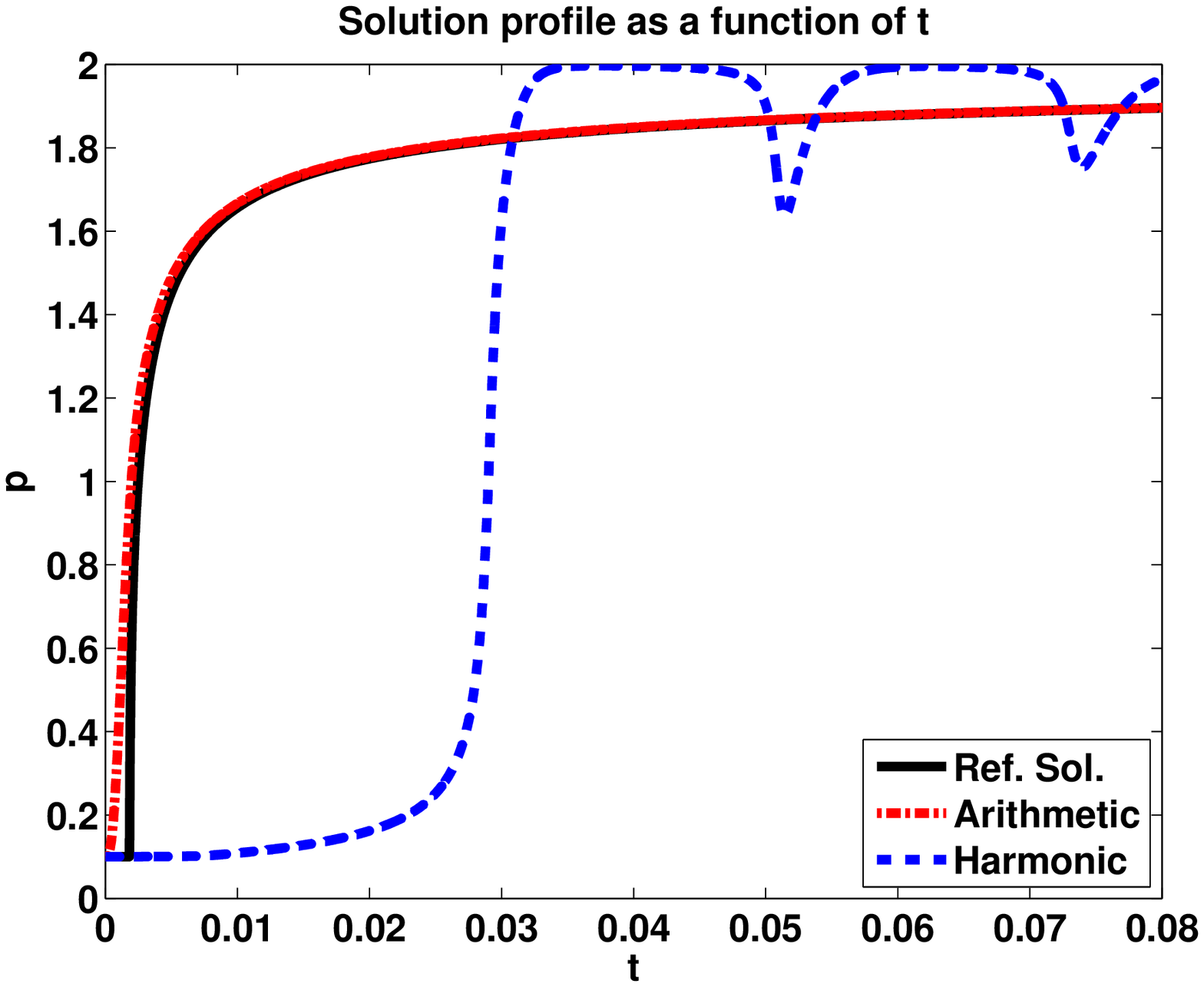}
			\caption {$x = 0.12$}
			\label{fig:pme_arithvsharm_time}
		\end{subfigure}
		\caption{Comparison of arithmetic and harmonic averaging for $k(p) = p^3$.  The simulations are initialized with the same initial condition shown in Figure \ref{fig:movinginterface}.  The spatial results match Figure 3 in \cite{inl2008} for $N = 25$ grid points.}
		\label{fig:pme_arithvsharm}
\end{figure}

Lagging is defined as the numerical interface moving slower than the analytical interface.  
Figure \ref{fig:pme_arithvsharm_space} illustrates the lagging present in the solution with harmonic averaging.  The solution with arithmetic averaging more accurately captures the front location on coarse grids.  A fine grid of $N = 800$ grid points is required for the numerical solution with harmonic averaging to overcome this lagging phenomenon \cite{inl2008}.  
Figure \ref{fig:pme_arithvsharm_space} also shows that the numerical speed of the moving interface with arithmetic averaging is a little too fast and that this solution is a bit diffusive near the discontinuity.     

Temporal oscillations are measured by tracking the unknown $p_i$ at a specific grid point $x_i$ over time.  By self-similarity, the temporal plots have the same profile for any $x$-coordinate that the front has passed through.  Here, we track $p(x,t)$ at $x = 0.12$.  Figure \ref{fig:pme_arithvsharm_time} shows the spurious temporal oscillations in the solution with harmonic averaging.  Figures \ref{fig:pme_arithvsharm_space} and \ref{fig:pme_arithvsharm_time} reveal that the temporal oscillations tend to occur even when there are no visible spatial oscillations.  In Figure \ref{fig:pme_arithvsharm_time}, we see that the numerical solution with arithmetic averaging is smooth.  

The temporal plots provide additional information about the behavior of the numerical solution than provided in the spatial plots.  For example, on refined grids, the spatial behavior of the solution with harmonic averaging appears correct, whereas in time high frequency oscillations are present.  The presence or absence of these temporal oscillations can then be used as a metric of accuracy for a numerical method. 

  \begin{figure}[H]
		\center
		\begin{subfigure}[H]{.42\textwidth}  
			\includegraphics[width =\textwidth]{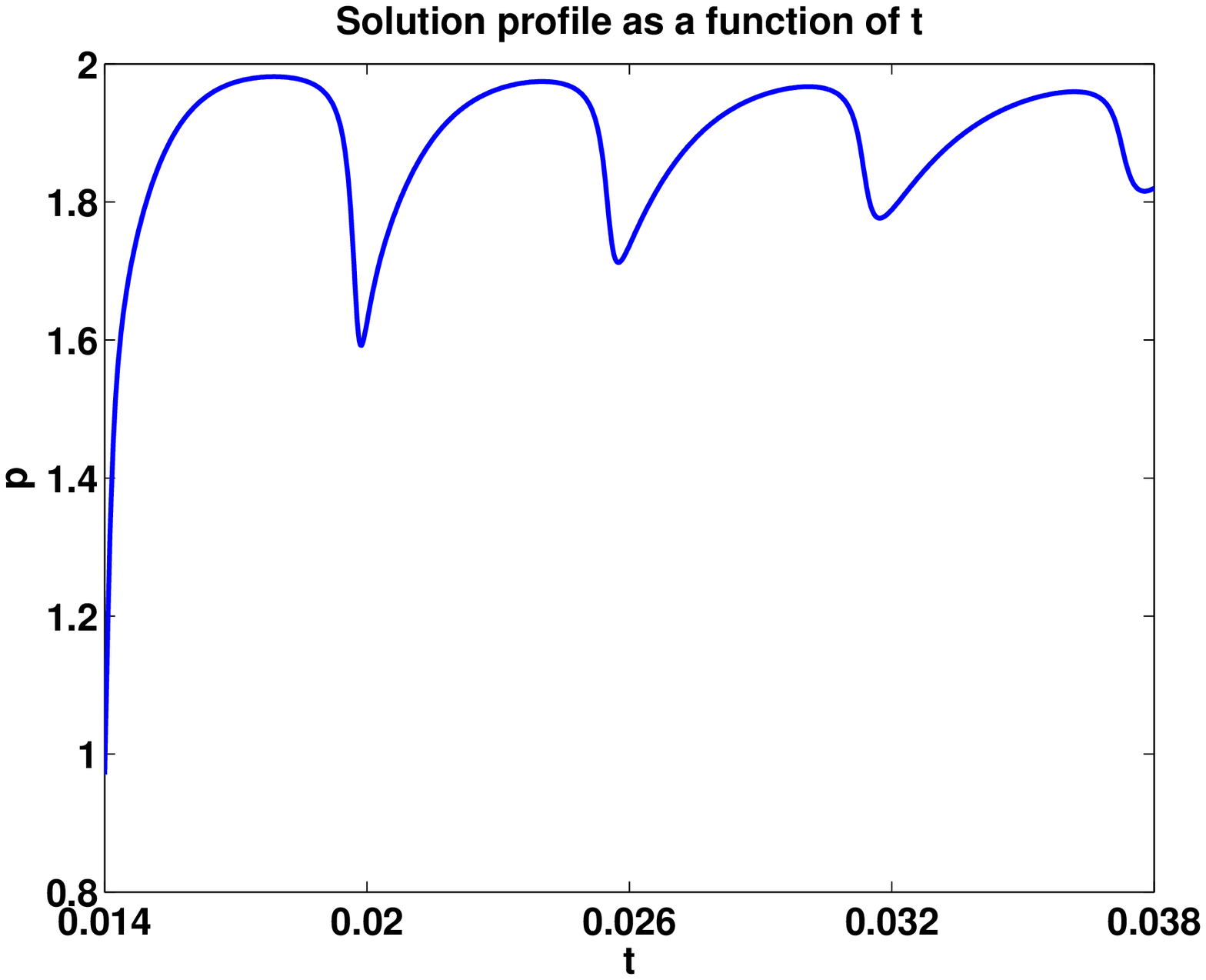}
			\caption{Period of temporal oscillations is 0.06.}
			\label{fig:time_nx50}
		\end{subfigure}
		\begin{subfigure}[H]{.42\textwidth}  
			\includegraphics[width =\textwidth]{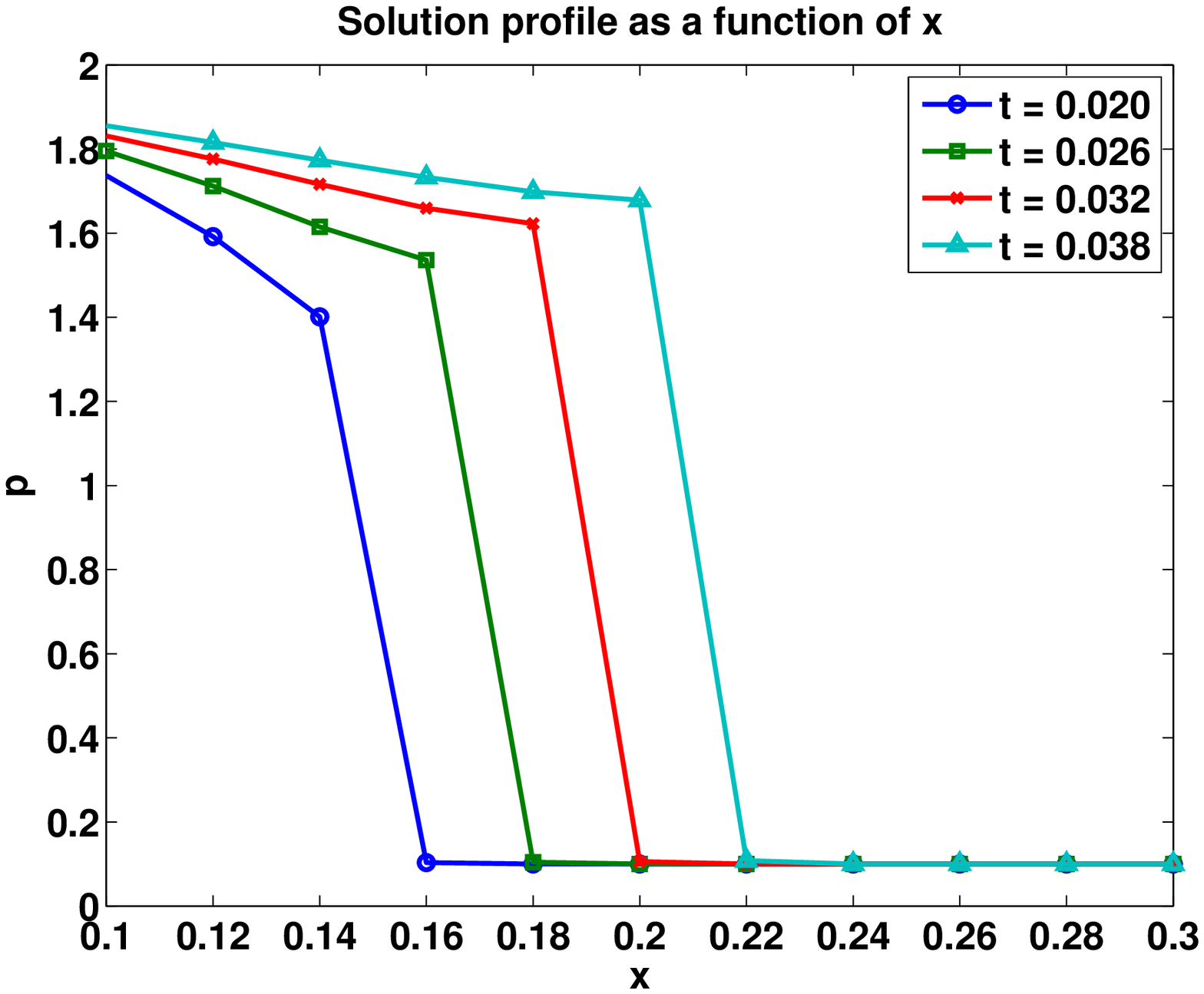}	
			\caption{$p(x,t)$ at each $t$ defining an oscillation minimum}
			\label{fig:time_nx50pos}	
		\end{subfigure}
		\caption{Numerical solution with harmonic averaging and $\dx = 0.02$.}
\end{figure}
\begin{figure}[H]
		\center
		\begin{subfigure}[H]{.42\textwidth}  
			\includegraphics[width =\textwidth]{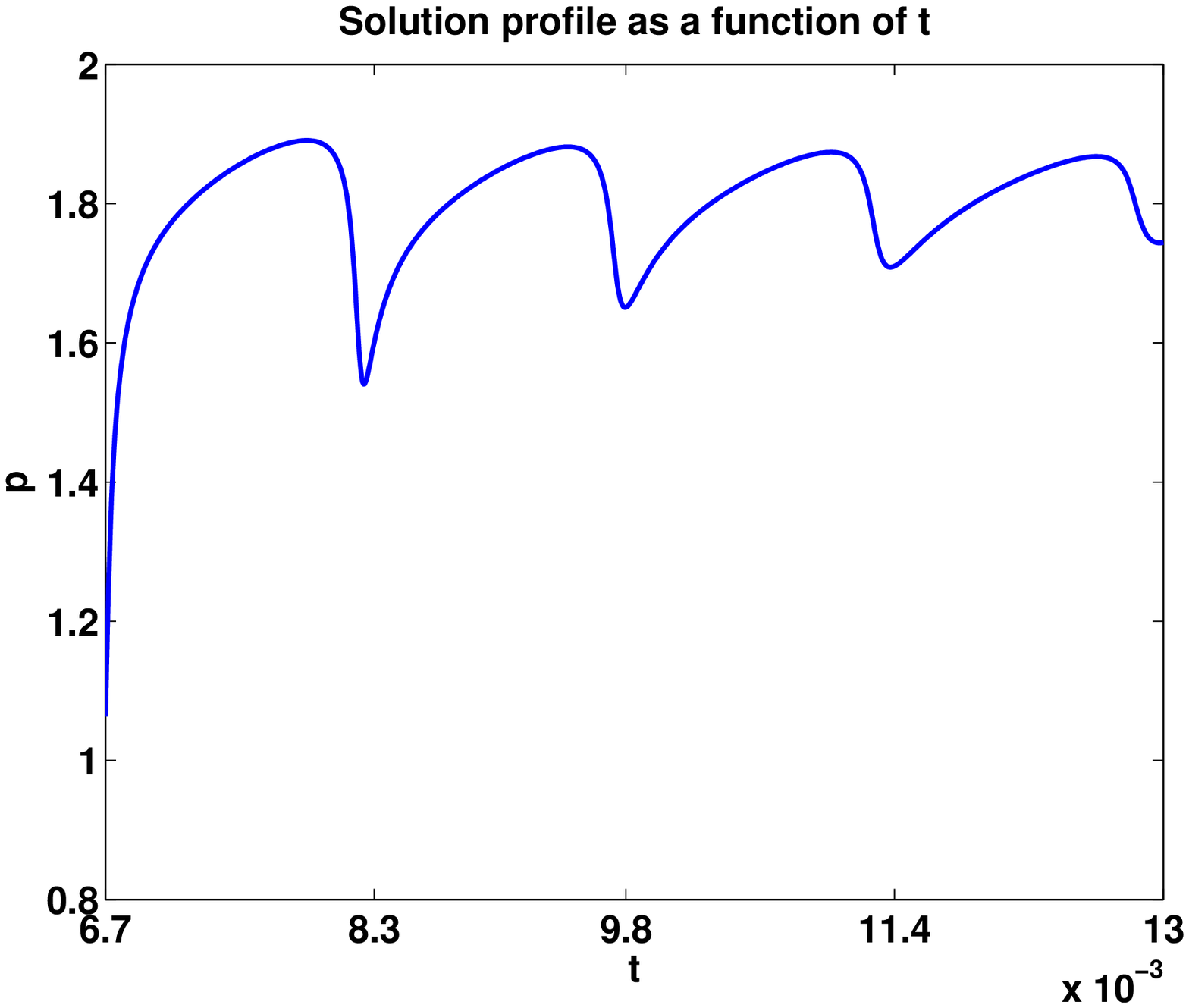}
			\caption{Period of temporal oscillations is 1.6e-3.}
			\label{fig:time_nx100}
		\end{subfigure}	
		\begin{subfigure}[H]{.42\textwidth}  
			\includegraphics[width =\textwidth]{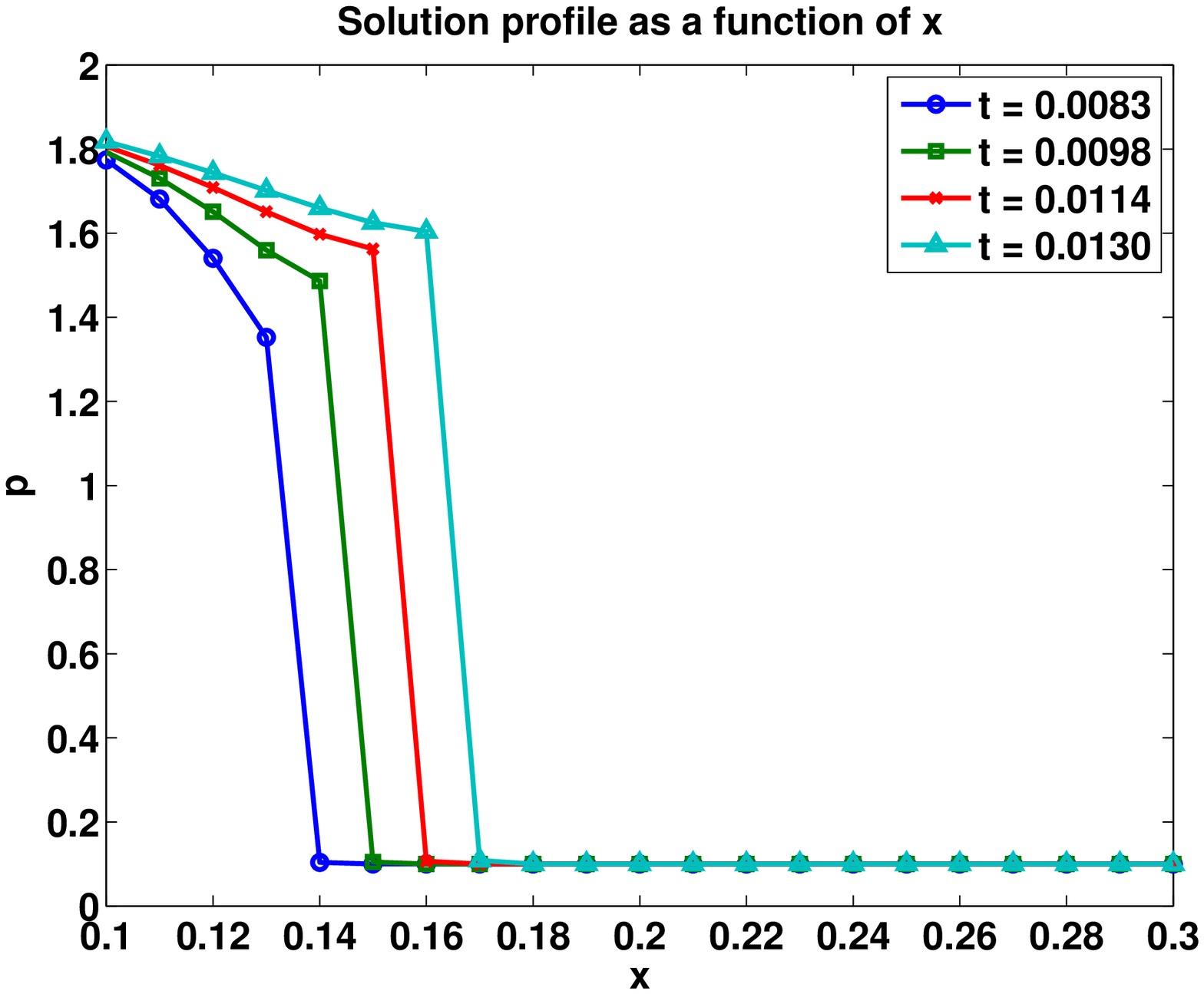}
			\caption{$p(x,t)$ at each $t$ defining an oscillation minimum}
			\label{fig:time_nx100pos}	
		\end{subfigure}
		\caption{Numerical solution with harmonic averaging and $\dx = 0.01$.}
\end{figure}

The temporal oscillations do not behave like those caused by numerical dispersion, as is often discussed in hyperbolic conservation laws.  Figures \ref{fig:time_nx50} and \ref{fig:time_nx100} depict that in the former refinement results in higher frequencies and lower amplitudes, whereas in the latter the amplitude does not decrease.  Figures \ref{fig:time_nx50pos} and \ref{fig:time_nx100pos} show that the formation of a temporal oscillation coincides with the time when the front crosses a grid cell \cite{vandermeer2016}.  Finer grids consist of more grid cells, explaining the higher frequency oscillations present.  Thus, refining the grid does not remove the temporal oscillations.  

Locking, an extreme case of lagging, is defined by the numerical solution not moving in time. 
This locking is illustrated for the PME with $m = 3$ \cite{lipnikov2016} in the locking problem (TLP).  Here, the initial condition $p(x,0) = h(x) = 10^{-3}$, and so $k(p(x,0)) = 10^{-9}$  for all $x \in \Omega$, magnitudes smaller than in previous examples.
With the Dirichlet boundary condition $p(0,t) = (3t)^{1/3}$, an analytic solution  $p_{\text{exact}}(x,t) = k(x,t)^{1/3}$  exists 
with
\begin{equation}
		k(x,t) = 
		\begin{dcases}
  			3(t-x),           & \text{if} \hspace{.1cm}x < t,  \\
    			\epsilon,              & \text{otherwise}.
		\end{dcases}
		\label{eq:mim_exact}
\end{equation}
for $\epsilon = 0$  \cite{lipnikov2016}.  To be consistent with the numerical simulations in \cite{lipnikov2016}, we do not use $\epsilon = 0$ and instead take $\epsilon = 10^{-9}$. 



\begin{figure}[H]
		\center
		\begin{subfigure}[H]{0.32\textwidth}  
			\includegraphics[width =\textwidth]{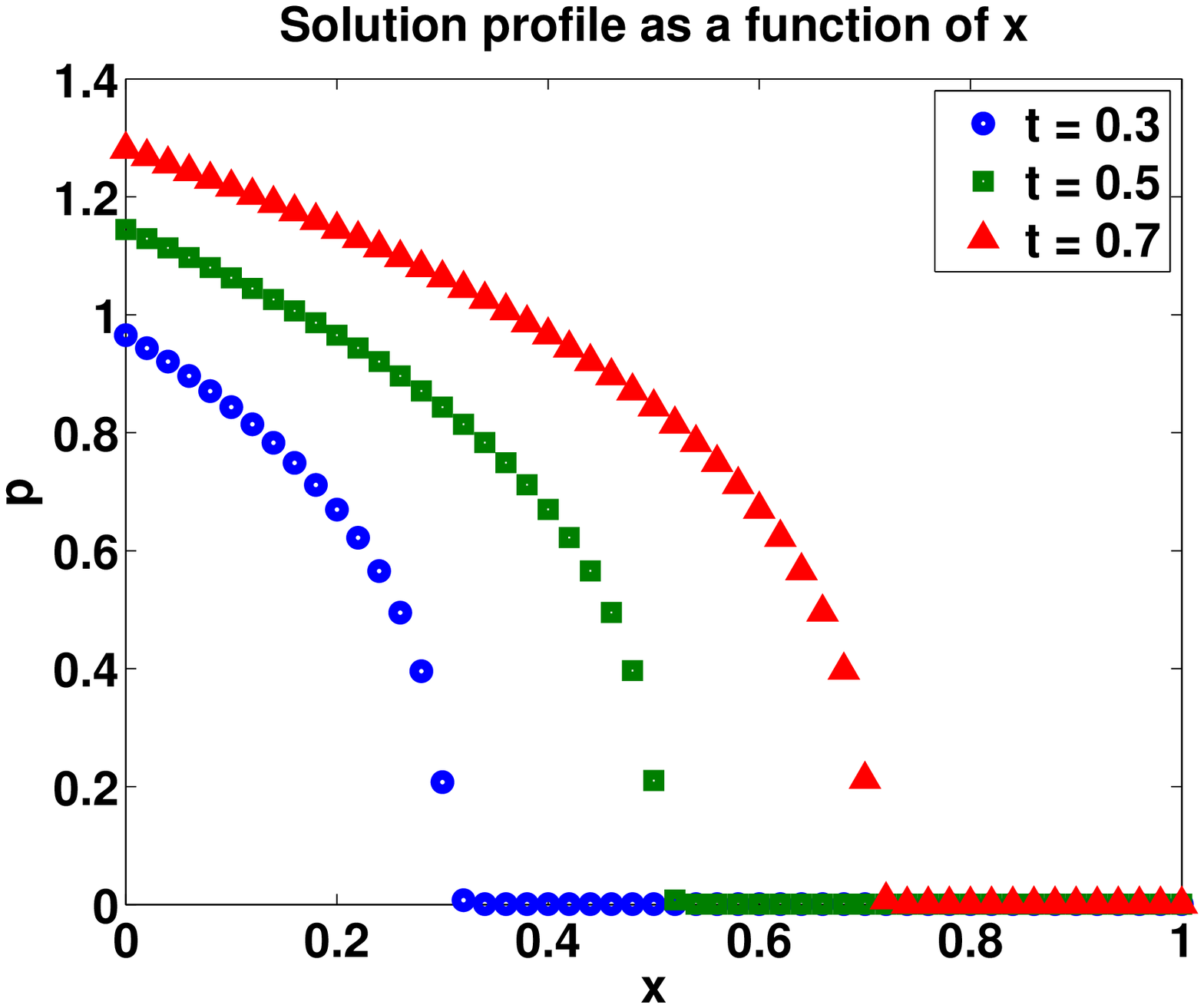}
			\caption {Arithmetic Averaging}
			\label{fig:locking_arith}
		\end{subfigure}
		\begin{subfigure}[H]{.32\textwidth}  
			\includegraphics[width =\textwidth]{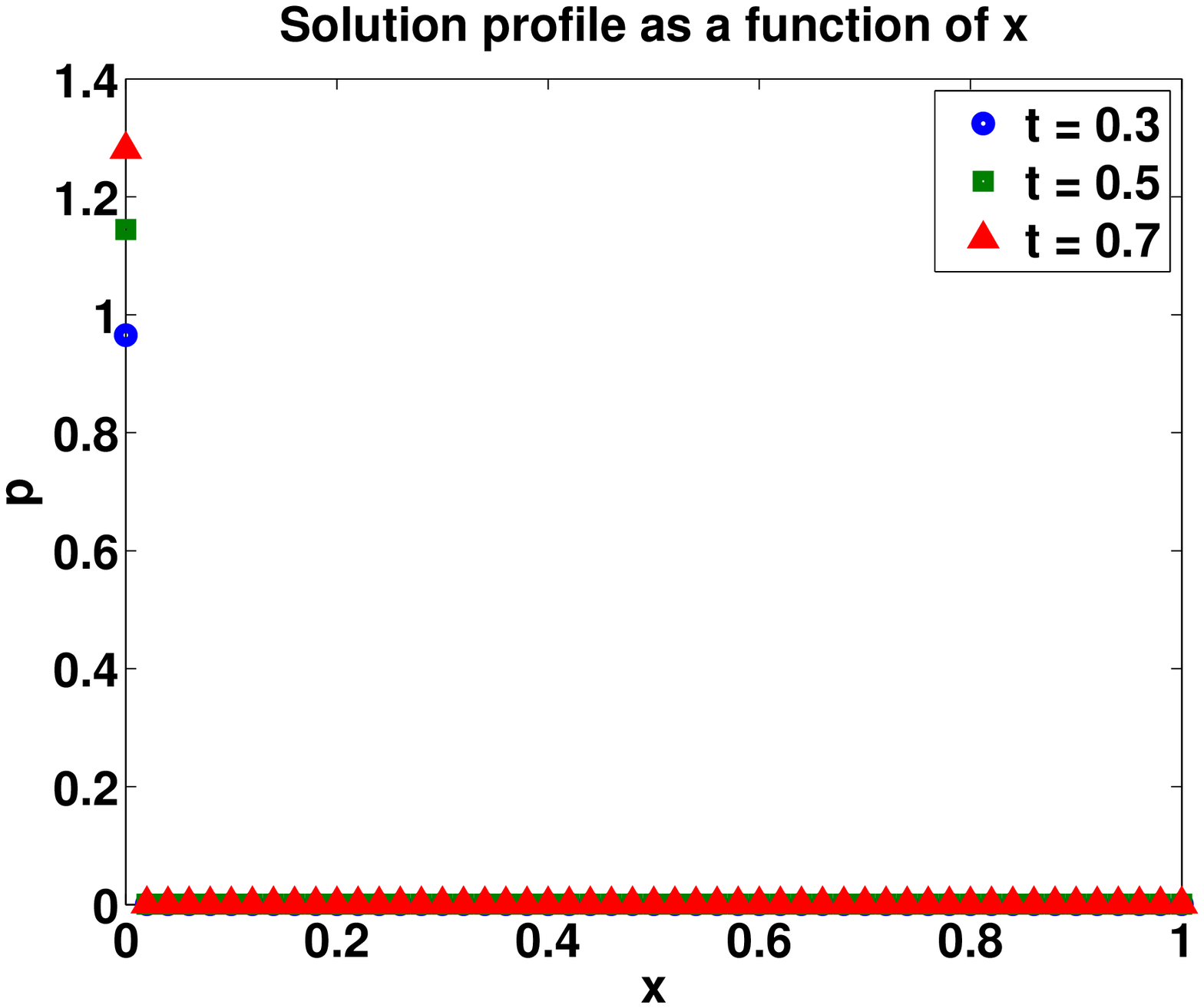}
			\caption {Harmonic Averaging}
			\label{fig:locking_harm}
		\end{subfigure}
		\begin{subfigure}[H]{0.32\textwidth}  
			\includegraphics[width =\textwidth]{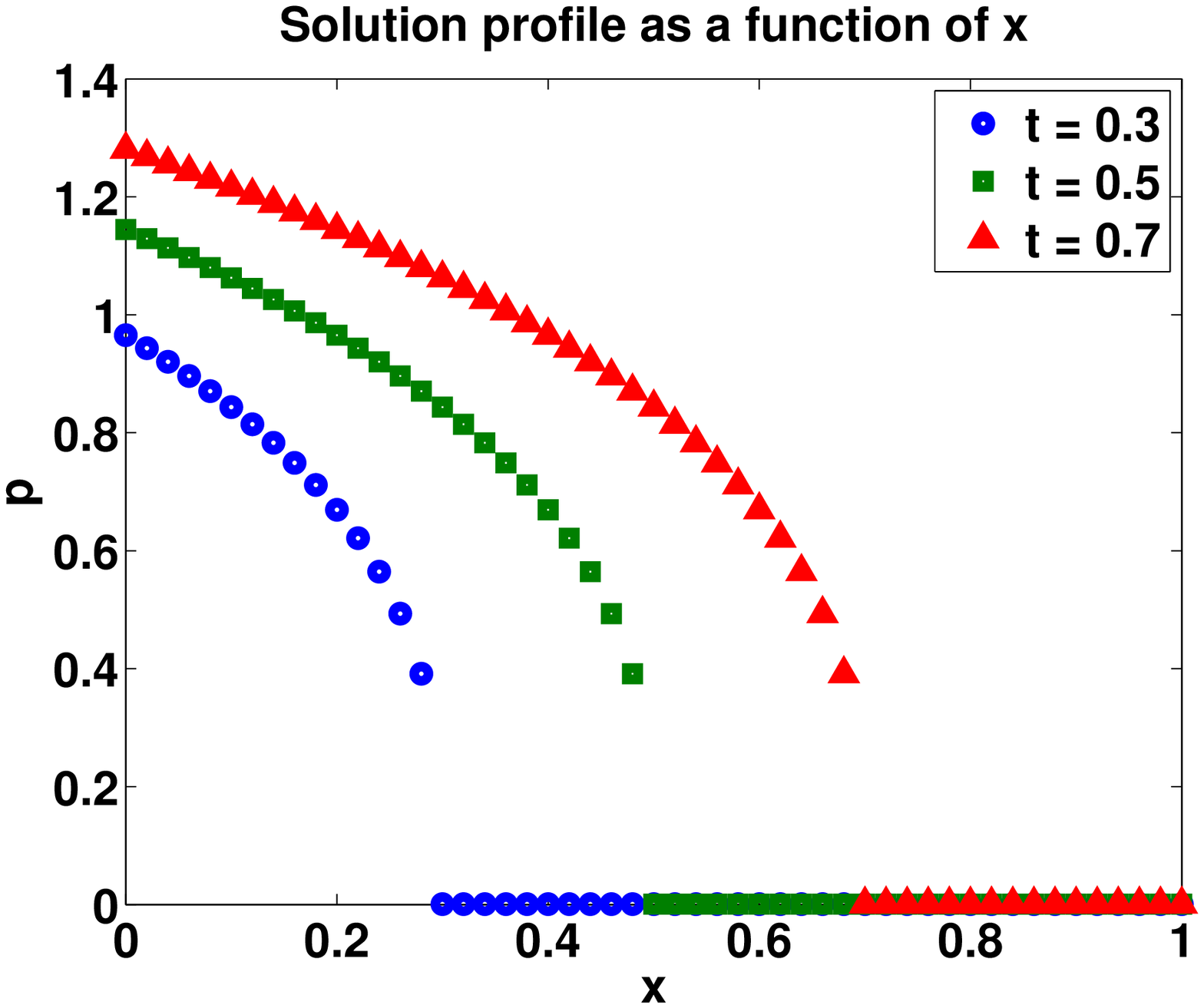}
			\caption {Exact Solution}
			\label{fig:locking_exact}
		\end{subfigure}
		\caption{Comparison of arithmetic and harmonic averaging for $k(p) = p^3$ with $N = 50$ grid points in the more extreme case, where artificial locking is observed with harmonic averaging. The results match Figure 1 in \cite{lipnikov2016}.}
\end{figure}

Figure \ref{fig:locking_harm} illustrates that the solution with harmonic averaging locks.  Here, the temporal profile is constant at the initial value, and so no temporal oscillations are observed.  
Figure \ref{fig:locking_arith} reveals that the numerical solution with arithmetic averaging is more accurate than the solution with harmonic averaging and that it tracks the exact solution in Figure \ref{fig:locking_exact} well.



\subsection{Effects of Temporal Discretizations on the Numerical Artifacts}
\label{TVD}
Since the oscillations are occurring in time, it seems reasonable to assume that the time-stepping method would have an effect on their presence or absence.  In this subsection, we explore the use of implicit Euler, which is stable for this problem, as well as second order total variation diminishing (TVD) discretizations to examine whether changing the stability and accuracy could affect the numerical artifacts.

\begin{figure}[H]
		\center
		\begin{subfigure}[H]{.32\textwidth}  
			\includegraphics[width =\textwidth]{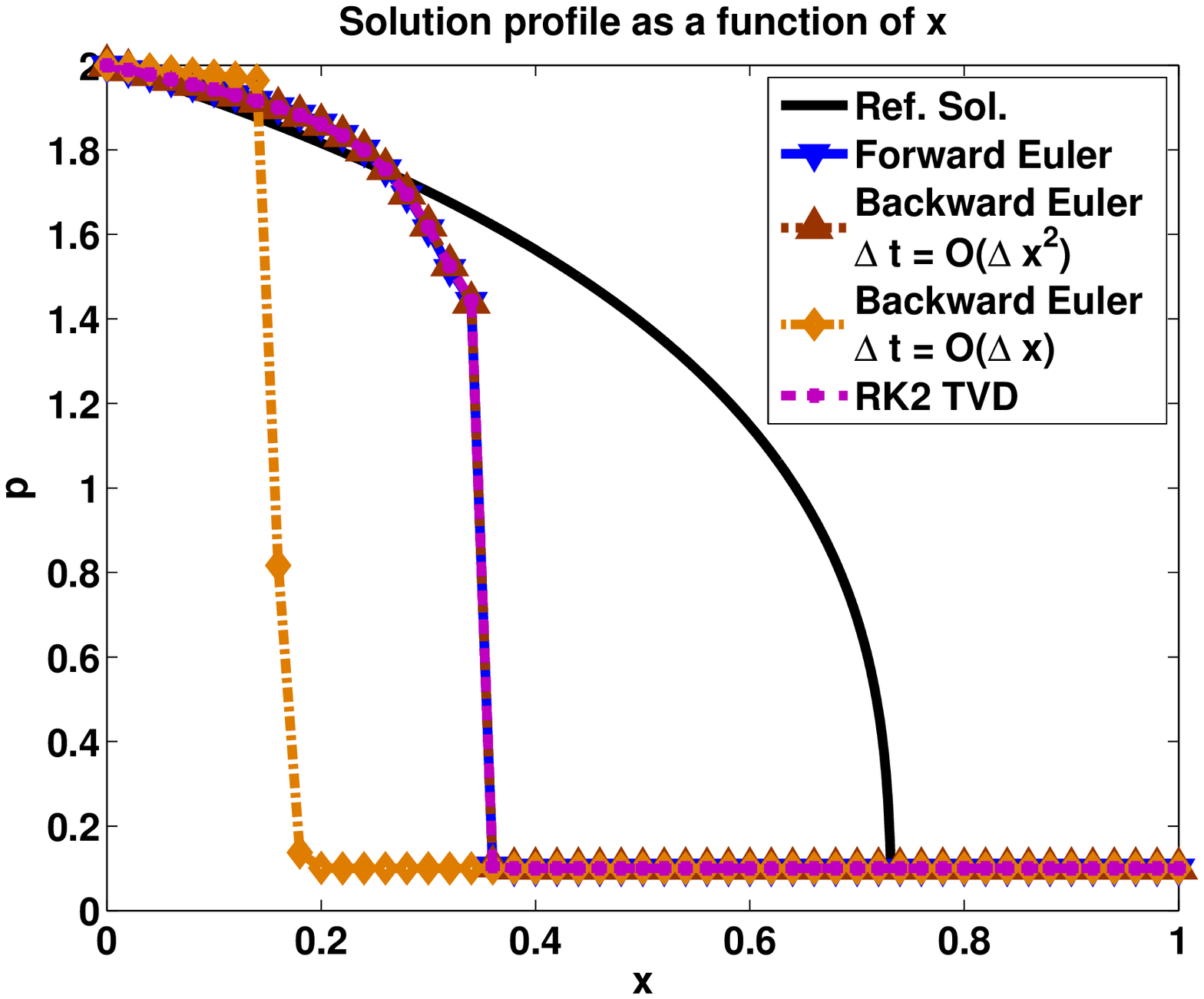}
			\caption {$t = 0.08$}
			\label{fig:timepos}
		\end{subfigure}
		\begin{subfigure}[H]{0.32\textwidth}  
			\includegraphics[width =\textwidth]{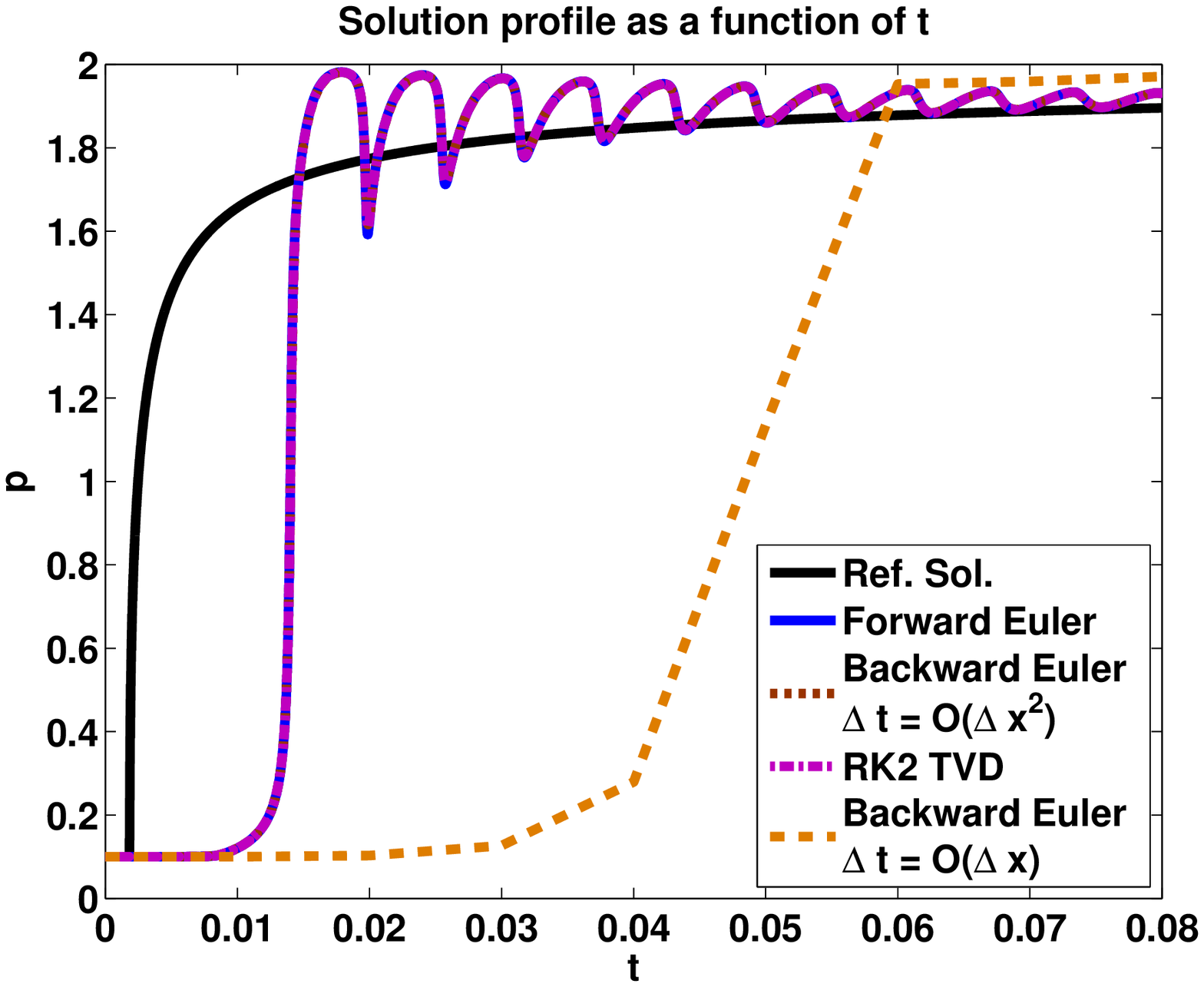}
			\caption {$x = 0.12$}
			\label{fig:timetemp}
		\end{subfigure}
		\begin{subfigure}[H]{0.34\textwidth}  
			\includegraphics[width =\textwidth]{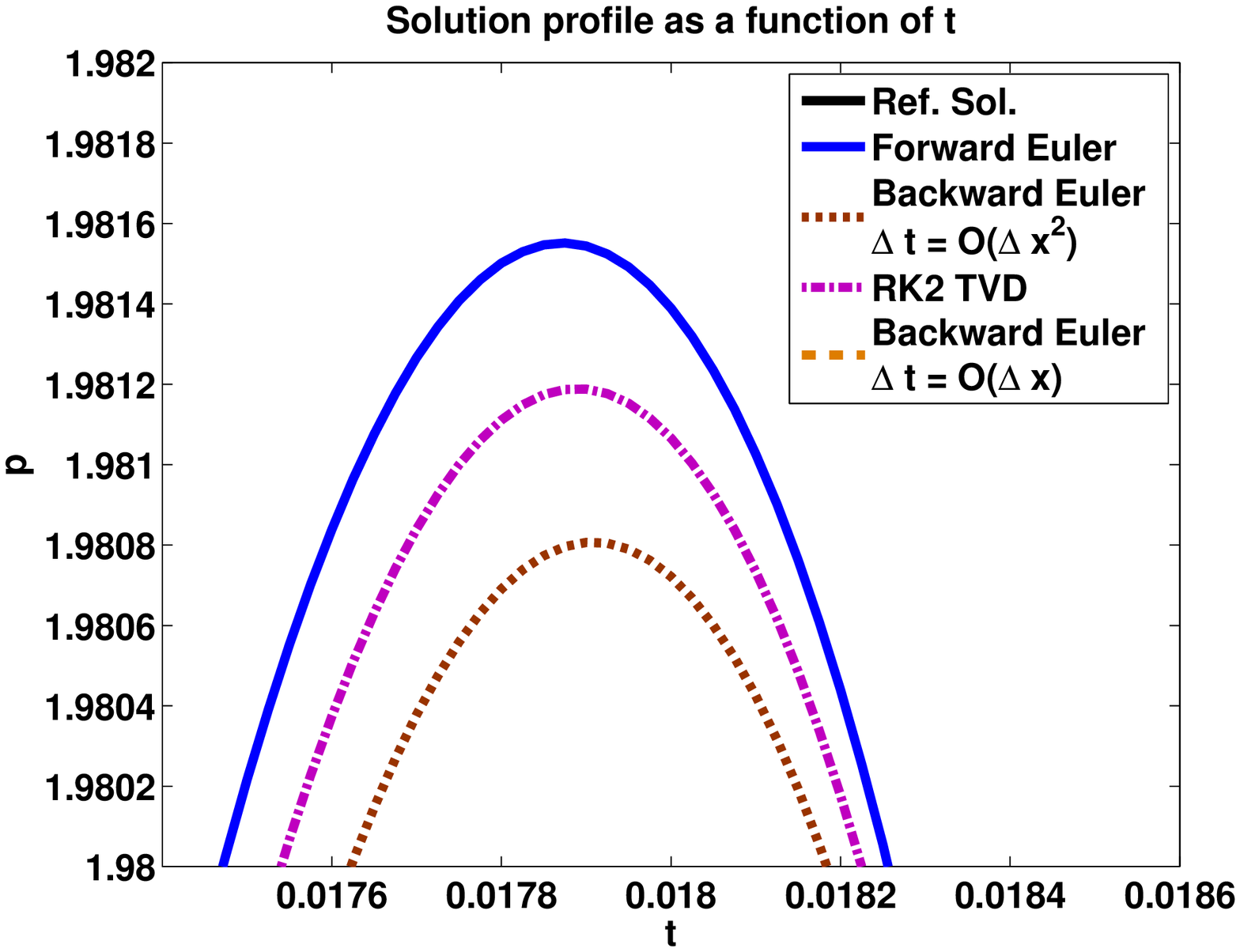}
			\caption {Zoomed-in Time}
			\label{fig:time_zoom}
		\end{subfigure}
		\caption{Numerical solution of the second order finite volume method with harmonic averaging for various temporal schemes.  The results are for $k(p) = p^3$ and $N = 50$ grid points.}
		\label{fig:time}
\end{figure}

For this problem, extensive numerical stability tests give confidence in the stability region of Backward Euler and that larger time steps can be taken than for Forward Euler, where $\dt = \mathcal{O}(\dx^2)$.  As expected, the numerical convergence studies reveal that in general the time step $\dt \propto \mathcal{O}(\dx^2)$ for Forward Euler in the PME is decreased as $m$ is increased due to sharper solutions.  Figure \ref{fig:timetemp} shows that there are no temporal oscillations for the scheme with Backward Euler and $\dt = \mathcal{O}(\dx)$, but for $\dt = \mathcal{O}(\dx^2)$ they are still present.  
However, taking Backward Euler with a large time step results in severe lagging and an incorrect shock position, as seen in 
Figure \ref{fig:timepos}. 
Clearly, there is a trade-off between spatial accuracy and a smooth temporal profile for this implicit temporal method.

We implement the TVD Runge Kutta (RK) 2 scheme introduced in \citet{shu98}, which was designed to reduce spatial oscillations,
but does not guarantee anything about the temporal profile.  
Stability estimates give that $\dt = \mathcal{O}(\dx^2)$ must be used here and the results are comparable with Forward and Backward Euler time-stepping for the same time step.  Figure \ref{fig:time_zoom} displays that the amplitude of the temporal oscillations with TVD RK2 is less than that with Forward Euler and larger than that with diffusive Backward Euler.

Since the artifacts still persist with the above temporal discretizations, in the next section we focus on Forward Euler time-stepping and 
use the Modified Equation Analysis to explain the difference in behaviors between the FTCS schemes with arithmetic and harmonic averaging.


\section{Modified Equation Analysis of the GPME for differentiable $k(p)$}
\label{MEA}

In this section, we derive the one-dimensional modified equations of the FTCS discretizations of Eqn. \eqref{eq:GPME} with arithmetic and harmonic averaging, respectively.  In the analysis that uses Taylor series expansions, we make the common assumption that $p(x,t)$ and $k(p(x,t))$ are infinitely differentiable with bounded derivatives.  %
We show that the leading error terms differ in two important ways.
We use the numerical discretization given by
 \begin{equation}
 	\begin{aligned}
	\frac{p_i^{n+1} - p^n_i}{\dt} &= \frac{1}{\dx}\bigg[k_{i+1/2}\Big(\frac{p^n_{i+1}-p^n_i}{\dx}\Big) -k_{i-1/2}\Big(\frac{p^n_{i}-p^n_{i-1}}{\dx}\Big)\bigg] = \frac{1}{\dx}\bigg[F_{i-1/2} -F_{i+1/2}\bigg],
	\label{eq:scheme}
	\end{aligned}
\end{equation}
where $k_{i+1/2}$ is a representation of $k(p)$ on the interval $[x_i,x_{i+1}]$ at time $t^n$, 
and the flux $F_{i+1/2}$ is defined in Eqn. \eqref{eq:fluxesdef}.  The superscript $n$ refers to the time step.  The $i\!-\!1/2$ terms have analogous definitions. 

\subsection{Modified Equation for the Scheme with Arithmetic Averaging}
For the semi-discrete portion, we are ultimately looking for an expression for the flux difference at cell center $i$.  This expression is derived in two steps, namely the fluxes are replaced with Taylor expansions about cell faces $i\!+\!1/2$ and $i\!-\!1/2$, and then the flux difference is expressed as a Taylor expansion about cell center $i$.

The truncation error for the arithmetic average $k_{i+1/2}^A$ of $k_i^n \equiv k(x_i,t^n)$ and $k_{i+1}^n \equiv k(x_{i+1}, t^n)$ is computed in Eqn. \eqref{eq:arith_trun}.  
Analogous to the definition of $k^+$ at time $t^n$ in Section \ref{arithharm_trunc}, we define 
\[
	 p^+ \equiv p(x_i + \dx/2,t^n), \hspace{.2cm} k^- \equiv k(x_i - \dx/2,t^n) \hspace{.1cm} \text{ and } \hspace{.1cm} p^- \equiv p(x_i - \dx/2,t^n),
\] to simplify the notation.  Let 
\[
	v = p^+  +  \frac{\dx^2}{8}p^+_{xx}  +  \mathcal{O}(\dx^4),
\] be the even terms and 
\[
	w = \frac{\dx}{2}p^+_x +  \frac{\dx^3}{48}p^+_{xxx} + \mathcal{O}(\dx^5),
\] be the odd terms in the expansion of $p^n_{i+1}$ about $x_i + \dx/2$.  Substituting $v$ and $w$ into the definition of $u_{i+1/2}$ at time $t^n$ in Eqn. \eqref{eq:vel}, we have 
\begin{equation}
	u^n_{i+1/2} = -\frac{(v+w) - (v-w)}{\dx} = -\frac{2w}{\dx} = -\Big(p^+_x +  \frac{\dx^2}{24}p^+_{xxx}\Big) + \mathcal{O}(\dx^4).
	\label{eq:velocity}
\end{equation}
Multiplying $u^n_{i+1/2}$ and $k^+\!+\!\tau^A$ from Eqn. \eqref{eq:arith_trun} leads to 
\begin{equation}
	 F^A_{i+1/2} = -\Big(k^+p_x^+ + \frac{\dx^2}{24}k^+p^+_{xxx} +  \frac{\dx^2}{8}k^+_{xx}p^+_x\Big) + \mathcal{O}(\dx^4).
	 \label{eq:term1}
\end{equation}
In the same way, we can derive
\begin{equation}
	F^A_{i-1/2} =  -\Big(k^-p_x^- + \frac{\dx^2}{24}k^-p^-_{xxx} +  \frac{\dx^2}{8}k^-_{xx}p^-_x\Big) + \mathcal{O}(\dx^4),
	\label{eq:term2}
\end{equation}
with the indices shifted. 
For the second semi-discrete step, we Taylor expand the fluxes $F_{i+1/2}$ and $F_{i-1/2}$ about $x_i$ to cell center $i$.
 Let
 \[
 	v = F(x_i,t^n) + \frac{\dx^2}{8}F_{xx}(x_i,t^n) +  \mathcal{O}(\dx^4),
\]   be the even terms and
\[
	w = \frac{\dx}{2}F_x(x_i,t^n) + \frac{\dx^3}{48}F_{xxx}(x_i,t^n) + \mathcal{O}(\dx^5),
\] be the odd terms in the expansion of $F_{i+1/2}$.  Then, $F_{i+1/2} = v+w$ and $F_{i-1/2} = v-w$.  Dividing the flux difference by $\dx$ gives the right hand side of Eqn. \eqref{eq:scheme} as
 \begin{equation}
 	-\frac{2w}{\dx} = -F_x(x_i,t^n) - \frac{\dx^2}{24}F_{xxx}(x_i,t^n) + \mathcal{O}(\dx^4).
	\label{eq:fluxes}
\end{equation}
 The arithmetic flux $F^A(x_i,t^n)$ is defined by the expressions in Eqn.\eqref{eq:term1} and Eqn. \eqref{eq:term2}  evaluated at $x_i$ and $t^n$.
 Substituting the expression for $F^A(x_i,t^n)$ into Eqn. \eqref{eq:fluxes} gives 
\[
	(kp_x)_x + \frac{\dx^2}{24}(kp_{xxx})_x +  \frac{\dx^2}{8}(k_{xx}p_x)_x + \frac{\dx^2}{24}(kp_{x})_{xxx} + \mathcal{O}(\dx^4),
\]
or
\begin{equation}
	(kp_x)_x + \dx^2\Big(\frac{1}{12} kp_{xxxx} + \frac{1}{6} k_x p_{xxx} +  \frac{1}{4} k_{xx}p_{xx} + \frac{1}{6}k_{xxx}p_x\Big)  + \mathcal{O}(\dx^4),
\label{eq:semidiscrete_pre}
\end{equation}
as the right hand side of Eqn. \eqref{eq:scheme},  where $p \equiv p(x_i,t^n)$ and $k \equiv k(p(x_i,t^n))$.

Now that the spatial semi-discrete part is complete, the next part in the Modified Equation Analysis is to Taylor expand $p_i^{n+1} = p(x_i, t^{n+1})$ in time about $t^n$.  Substituting this Taylor expansion of $p_i^{n+1}$ and Eqn. \eqref{eq:semidiscrete_pre} into Eqn. \eqref{eq:scheme} leads to \begin{equation}
  	\begin{aligned}
		\frac{1}{\dt} \Big(p + \dt p_t + \frac{\dt^2}{2} p_{tt} - p\Big) &= (kp_x)_x + \dx^2\Big(\frac{1}{12} kp_{xxxx}  + \frac{1}{6} k_x p_{xxx} +  \frac{1}{4} k_{xx}p_{xx} + \frac{1}{6}k_{xxx}p_x\Big) +\mathcal{O}(\dt^2 + \dx^4).
		\nonumber
	\end{aligned}
  \end{equation}
With multiple applications of the chain and product rules, the above equation can be written as
  \begin{equation}
  	\begin{aligned}
  	p_t -(kp_x)_x  = &-\frac{\dt}{2} p_{tt}  + \dx^2\Big(\frac{1}{12} kp_{xxxx} + \frac{1}{3} k_pp_x p_{xxx} +  \frac{3}{4} k_{pp}p_x^2p_{xx} + \frac{1}{4}k_pp_{xx}^2 + \frac{1}{6}k_{ppp}p_x^4\Big) +\mathcal{O}(\dt^2 + \dx^4).
	\end{aligned}
	\label{eq:semidiscrete}
  \end{equation}

We use the governing equation \eqref{eq:GPME} to convert $p_{tt}$ in Eqn. \eqref{eq:semidiscrete} within the order of accuracy into spatial derivatives as 
  \begin{equation}
  	\begin{aligned}
  p_{tt} &=  (kp_x)_{xt} = (kp_x)_{tx} = (k_tp_x + kp_{tx})_x = (k_tp_x + k(kp_x)_{xx})_x  \\
  &= k_tp_{xx} + k_{tx}p_x + k_x(kp_x)_{xx} + k(kp_x)_{xxx},
  	\label{eq:ptt}
	\end{aligned}
\end{equation}
where
 \[
  	\begin{aligned}	
				&k_t = k_pp_t = k_p(kp_x)_x =  k_p^2p_x^2+kk_pp_{xx}, \\
		&k_{tx} = 2k_pk_{pp}p_x^3 + (3k_p^2+ kk_{pp})p_xp_{xx}  + kk_pp_{xxx}. \\
	\end{aligned}
  \]
Eqn. \eqref{eq:ptt} combined with the above expressions and the chain and product rules gives
 \begin{equation}
 	p_{tt} = 4kk_pp_{xx}^2 + (3k_pk_{pp} + kk_{ppp})p_x^4 + 7(k_p^2+ kk_{pp})p_x^2p_{xx}  + 6kk_pp_{xxx}p_x + k^2p_{xxxx}.
	\label{eq:ptt_final}
\end{equation}
Because of the type of PDE, conversion from temporal to spatial derivatives is more straightforward and efficient.
 
Lastly, we substitute Eqn. \eqref{eq:ptt_final} into Eqn. \eqref{eq:semidiscrete} to obtain the modified equation of the scheme with arithmetic averaging given by 
\begin{equation}
\begin{aligned}
	p_t -(kp_x)_x  &= 
	{p_x^2p_{xx}}\Big[\!-\!\frac{7\dt}{2}(k_p^2+ kk_{pp}) + \frac{3\dx^2}{4} k_{pp}\Big] + p_{xx}^2\Big[k_p\Big(\!-\!2k\dt +\frac{\dx^2}{4}\Big)\Big]
	+p_{xxx}p_x\Big[k_p\Big(-3k\dt + \frac{\dx^2}{3}\Big)\Big] \\
	&+ {p_x^4}\Big[\!-\!\frac{\dt}{2}( 3k_pk_{pp} + kk_{ppp})  + \frac{\dx^2}{6}k_{ppp}\Big] + p_{xxxx}\Big[k\Big(\!-\!k\frac{\dt}{2} + \frac{\dx^2}{12}\Big)\Big] + \mathcal{O}(\dt^2 + \dx^4).
	\label{eq:arith_modeqtn}
\end{aligned}
\end{equation}

At this point, it is important to understand the order of the scheme.  As seen in Eqn. \eqref{eq:arith_modeqtn}, the leading terms are either of order $\dt$ or $\dx^2$.  In practice, we find after thorough numerical stability experimentation that $\dt$ is proportional to $\dx^2$ and the leading terms are of the same order $\dx^2$.

\subsection{Modified Equation for the Scheme with Harmonic Averaging} 
In deriving the analogous of Eqn. \eqref{eq:arith_modeqtn} for the scheme with harmonic averaging, we find that the modified equations  can be written in the same general form as 
\begin{equation}
	\begin{aligned}
	p_t -(kp_x)_x  &= 
	{p_x^2p_{xx}}\big[\mathcal{A} + \mathcal{B}\big] + p_{xx}^2\mathcal{C}
	+p_{xxx}p_x\mathcal{D}
	+ {p_x^4}\big[\mathcal{E}  + \mathcal{F}\big] + p_{xxxx}\mathcal{G} + \mathcal{O}(\dt^2 + \dx^4).
	\label{eq:gen_modeqtn}
\end{aligned}
\end{equation}
\vspace{-.5cm}
\begin{table}[H]
\begin{center}
 \label{table}
\center
\renewcommand{\arraystretch}{1.5}
\begin{tabular}{|c|c|c|}  
\hline 
& Arithmetic \\ \hline
$\mathcal{A}$	& $-\frac{7\dt}{2}(k_p^2+ kk_{pp})$ 		
\\  \hline
$\mathcal{B}$     &  $\frac{3\dx^2}{4} k_{pp}$     
 \\ \hline
$\mathcal{C}$	&$k_p\Big(\!-\!2k\dt +\frac{\dx^2}{4}\Big)$	
 \\ \hline
$\mathcal{D}$ & $k_p\Big(\!-\!3k\dt + \frac{\dx^2}{3}\Big)$ 
\\ \hline 
$\mathcal{E}$ & $-\frac{\dt}{2}\Big( 3k_pk_{pp} + kk_{ppp}\Big)$  
 \\ \hline
$\mathcal{F}$ & $\frac{\dx^2}{6}k_{ppp}$ 
 \\ \hline
$\mathcal{G}$ & $k\Big(\!-\!k\frac{\dt}{2} + \frac{\dx^2}{12}\Big)$ 
\\ \hline
\end{tabular}
\label{eq:harm_modeqtn}
 \end{center}
 \caption{Definitions of the coefficients in Eqn. \eqref{eq:gen_modeqtn} for the modified equation of the FTCS discretization with arithmetic averaging.}
 \label{table}
\end{table} 
For the harmonic average, the only coefficients that differ from those in Table \ref{table} are $\mathcal{B}$ and $\mathcal{F}$.  We define $\mathcal{B}^H \equiv \mathcal{B} + \delta^H\mathcal{B}$, where 
\[
	\delta^H\mathcal{B} = -{\frac{3\dx^2}{4} k_p^2k^{-1}},
\]
and  $\mathcal{F}^H \equiv \mathcal{F} + \delta^H\mathcal{F}$, where 
\[
	\delta^H\mathcal{F} = -\frac{\dx^2}{4}\big(2k_pk_{pp}k^{-1} - k_p^3k^{-2}\big).
\]
These differences in the $\dx^2$ terms originate from 
$ -(\dx^2/4)[(k_x^+)^2/k^+]$
in Eqn. \eqref{eq:harm_trun}, and can be derived by first multiplying $u^n_{i+1/2}$ from Eqn. \eqref{eq:velocity} and $k^+\!+\!\tau^H$.  The expression for the harmonic flux $F^H_{i+1/2}$ is then given as
\begin{equation}
\begin{aligned}
	F^H_{i+1/2} &= u_{i+1/2}\Big[(k^+ + \tau^A) -\frac{\dx^2}{4} \frac{(k_x^+)^2}{k^+}\Big] \\
				&= F^A_{i+1/2} +\frac{\dx^2}{4}\frac{(k_x^+)^2}{k^+}p^+_x + \mathcal{O}(\dx^4).
	\end{aligned}
	 \label{eq:term12_harm}
\end{equation}
The analogous expression holds for $F^H_{i-1/2}$.
We follow the same procedure as done in the prior subsection and substitute Eqn. \eqref{eq:term12_harm} evaluated at $x_i$ into Eqn. \eqref{eq:fluxes}.  The arithmetic flux term leads to Eqn. \eqref{eq:semidiscrete_pre}, while the second term leads to 
\begin{equation}
	\begin{aligned}
		 -\frac{\dx^2}{4}\Big[k_p^2k^{-1}p_x^3\Big]_x &= 
		\underbracket{\Big[\!-\!\frac{3\dx^2}{4}k_p^2k^{-1}\Big]}_{\delta^H\mathcal{B}}p_x^2p_{xx} + \underbracket{\Big[\!-\!\frac{\dx^2}{4} (2k_pk_{pp}k^{-1} - k_p^3k^{-2})\Big]}_{\delta^H\mathcal{F}}p_x^4 \\
	\end{aligned}
	\label{eq:harm_terms}
\end{equation}
Eqn. \eqref{eq:harm_terms} identifies the additional terms in $\mathcal{B^H}$ and $\mathcal{F^H}$ for the scheme with harmonic averaging.

\subsubsection{Modified Equation Analysis for Backward Euler}
\label{BE}
The Modified Equation Analysis for the implicit Backward Euler method is similar to the Modified Equation Analysis shown for the explicit Forward Euler method.  For implicit Euler, the only difference is that Eqn. \eqref{eq:semidiscrete_pre} is defined at the time $t^{n+1}$.  Expanding about $t^n$ gives additional $\dt$ terms as leading error terms in the expansion.  
For $\dt = \mathcal{O}(\dx)$, the second-order $\delta\mathcal{B}^Hp_x^2p_{xx}$ term added by the harmonic average is no longer on the order of the leading error terms in the truncation error. We will see in the next section that $\delta\mathcal{B}^Hp_x^2p_{xx}$ has an anti-diffusive nature. 
It is then expected that temporal oscillations do not occur with $\dt = \mathcal{O}(\dx)$, whereas they do occur for $\dt = \mathcal{O}(\dx^2)$, as seen in Figure \ref{fig:timetemp}.



\section{Modified Harmonic Method (MHM)} 
\label{MH}
We conjecture and later show that $\mathcal{B}^Hp_x^2p_{xx}$ and $\mathcal{F}^Hp_x^4$ are the critical terms that contribute to the numerical artifacts.  For that reason, we propose the Modified Harmonic Method (MHM) that counteracts $\mathcal{B}^Hp_x^2p_{xx}$ and $\mathcal{F}^Hp_x^4$.  The MHM differs from the scheme with both harmonic and arithmetic averaging.  In the MHM, we take Eqn. \eqref{eq:scheme} with harmonic averaging as the base scheme and modify the right hand side with 
\begin{equation}
	\begin{aligned} 
		\underbrace{-\mathcal{B}^H\big[D^-(p_i^n)\big]^2\big[D^+(D^-(p_i^n))\big]}_{\RN{1}}\underbrace{ -\mathcal{F}^H\big[D^-(p_i^n)\big]^4}_{\RN{2}}, 
		\label{eq:MHM}
	\end{aligned}
\end{equation}
where
\[
	D^-(p_i^n) = \frac{p^n_i - p^n_{i-1}}{\dx} \hspace{.2cm} \text{ and } \hspace{.2cm} D^+(p_i^n) = \frac{p^n_{i+1} - p^n_{i}}{\dx},
\] 
 refer to the backward and forward difference operators, respectively.  Consistent with the discretization of the other terms in the FTCS scheme, the gradient is discretized with upwinding $D^-(p_i^n)$ and the second derivative is discretized with the second order central difference $D^+(D^-(p_i^n))$.  
 

The MHM only requires that the first three derivatives of $k(p)$ with respect to $p$ are defined.  These derivatives are substituted into $\mathcal{B}^H$ and $\mathcal{F}^H$, and then $\mathcal{B}^H$ and $\mathcal{F}^H$ are evaluated at $p_i^n$.
Since the term in Eqn. \eqref{eq:MHM} is proportional to $\dx^2$, any numerical errors made in the finite difference discretization are higher order error terms.  The MHM is also consistent since the term in Eqn. \eqref{eq:MHM} tends to zero as $\dx$ tends to zero.

 The following subsection provides three key findings, namely why and where the numerical artifacts are occurring, and also how to modify the discretization there.  The MHM can of course be applied just locally with a switch that turns it on or off.  An alternative would be to not switch from harmonic to MHM, but from harmonic to arithmetic when needed. 

\subsection{Effect and behavior of $\mathcal{B}^Hp_x^2p_{xx}$ and $\mathcal{F}^Hp_x^4$}

\begin{figure}[H]
		\center
 \begin{subfigure}{.42\textwidth}  
			\includegraphics[width =\textwidth]{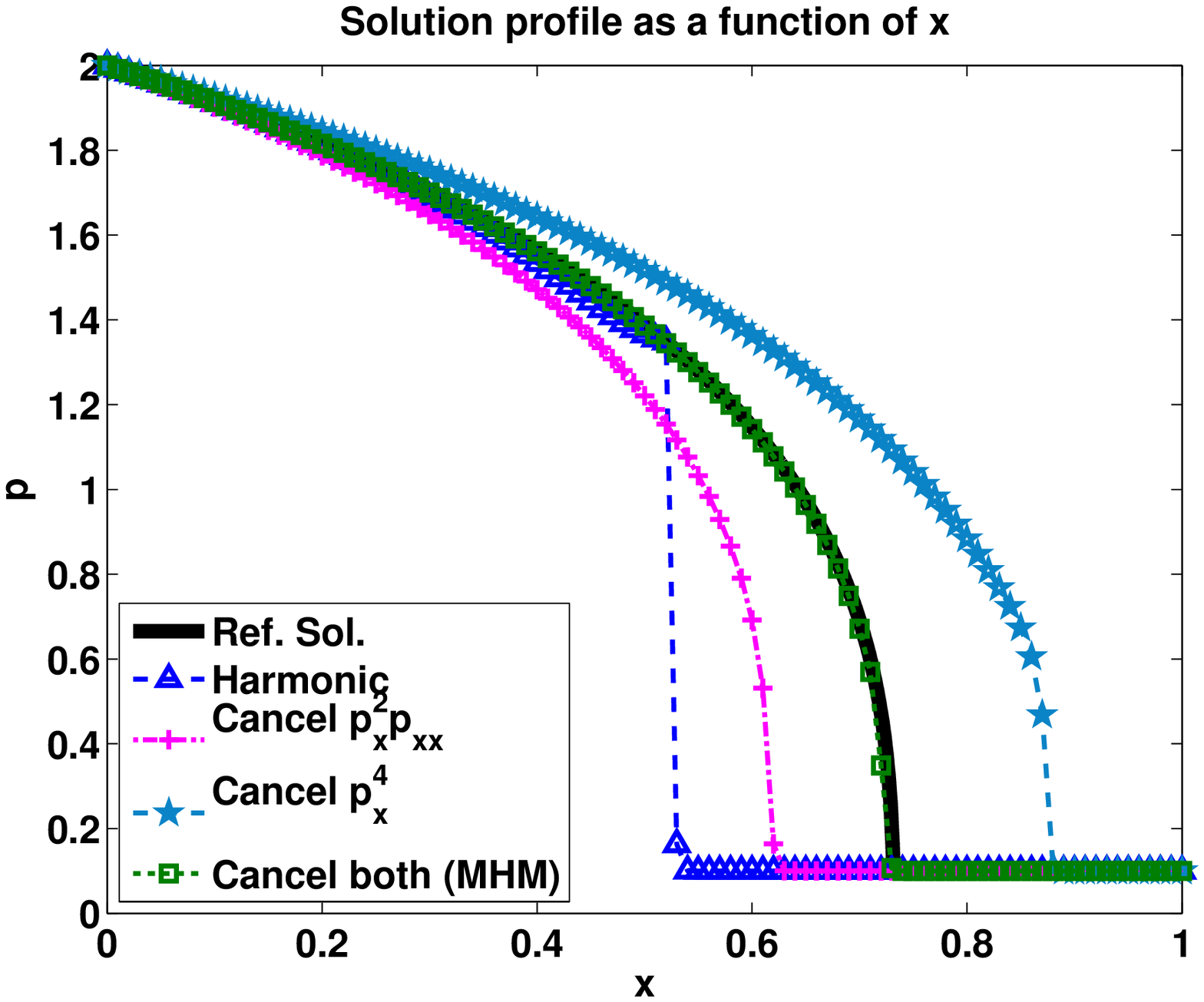}
			\caption {$t = 0.08$}
			\label{fig:px_pos}
		\end{subfigure}
		\begin{subfigure}{0.42\textwidth}  
			\includegraphics[width =\textwidth]{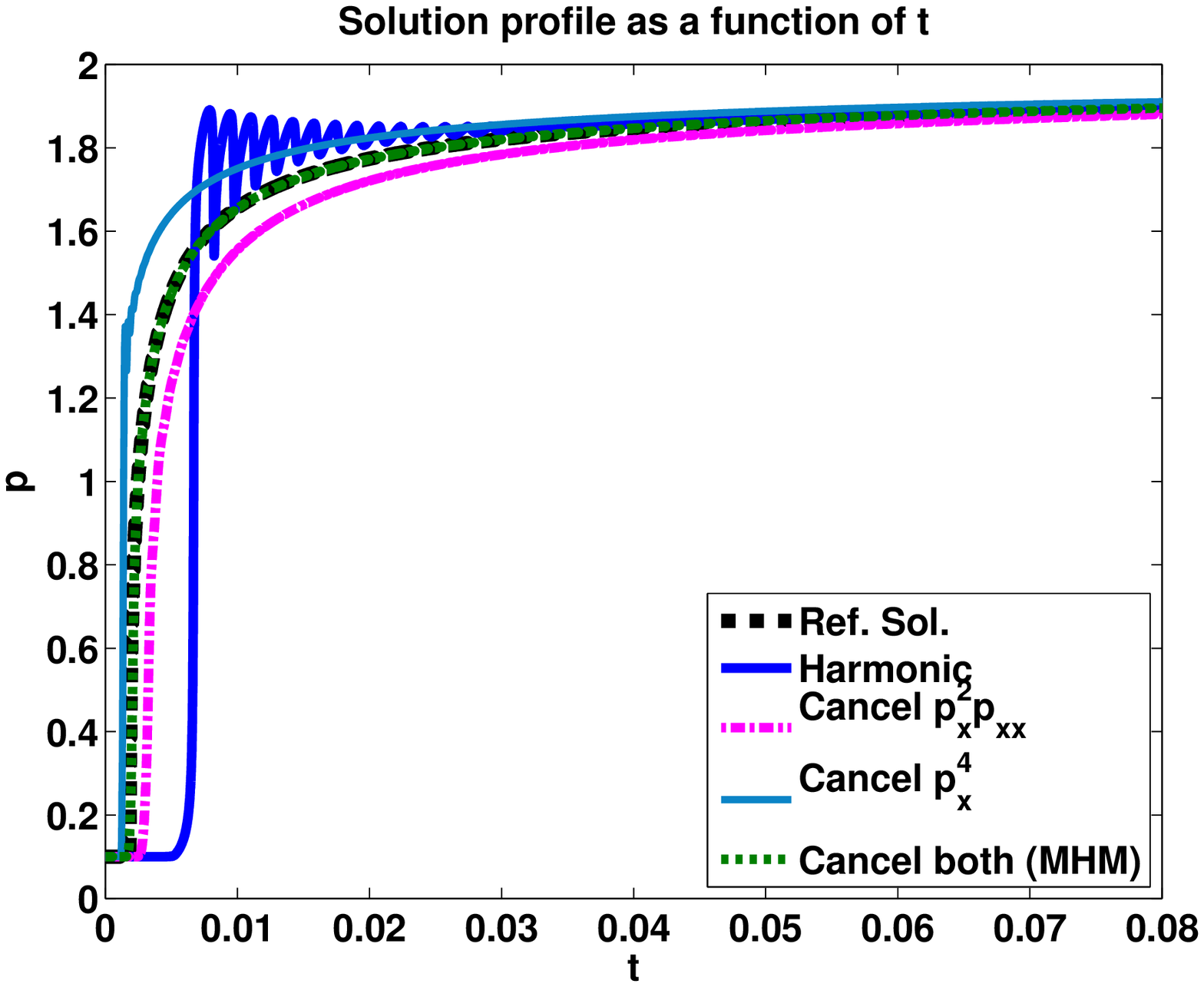}
			\caption {$x = 0.12$}
			\label{fig:px_time}
		\end{subfigure}
		\caption{Harmonic averaging results of $k(p) = p^3$ with $N = 100$ grid points and $\dt  = \dx^2/16$.}
		\label{fig:MHM}
\end{figure}

Figure \ref{fig:px_pos} shows that term $\RN{1}$, which counteracts $\mathcal{B}^Hp_x^2p_{xx}$, improves the speed of the moving interface.  
The lagging has been reduced, but is still present.  Thus, counteracting $\mathcal{B}^Hp_x^2p_{xx}$ alone is not enough for accurate front location detection. From Figure \ref{fig:px_time}, it is also clear that term $\RN{1}$ removes the temporal oscillations.

Figure \ref{fig:px_pos} also illustrates that term $\RN{2}$, which counteracts $\mathcal{F}^Hp_x^4$, has a greater effect on the interface speed than term $\RN{1}$.  In fact, the solution with only term $\RN{2}$ has a speed that is too high.  In this case, the lagging has been overcompensated.  Temporal oscillations are still present with term $\RN{2}$. 

From the above results, we can understand the behavior of $\mathcal{B}^Hp_x^2p_{xx}$ and $\mathcal{F}^Hp_x^4$.  Looking at these results, it seems that $\mathcal{F}^Hp_x^4$ is of an advective nature, which can be seen by writing the term as $\big[\mathcal{F}^Hp_x^3\big]p_x$.  The $\mathcal{B}^Hp_x^2p_{xx}$ term acts anti-diffusive and advective.  The two effects are of no surprise since this is a mixed term and can be seen by writing the term as $(\mathcal{B}^Hp_x^2)p_{xx}$  and $(\mathcal{B}^Hp_xp_{xx})p_x$, respectively.  

To take a closer look at the cause of the temporal oscillations, we further analyze the $\mathcal{B}^Hp_x^2p_{xx}$ term.  Recall that $\mathcal{B}^H$ is comprised of two pieces, namely $\mathcal{B}$ and $\delta^H\mathcal{B}$.  The term $\mathcal{B}p_x^2p_{xx} = \big[(3\dx^2/4)k_{pp}p_x^2\big]p_{xx}$ that is common in both modified equations for harmonic and arithmetic averaging has a positive diffusive effect.  For the PME and superslow diffusion, it can be verified that $k_{pp} > 0$ at the front.  Since the MHM counteracts $\mathcal{B}p_x^2p_{xx}$, it results in a scheme that is a bit less diffusive than the scheme with arithmetic averaging.
The problems occur with the additional term $\delta^H\mathcal{B}p_x^2p_{xx} = -\big[(3\dx^2/4)k_p^2k^{-1}p_x^2\big]p_{xx}$.  It has an anti-diffusive effect for all $k(p) > 0$, since the coefficient of $p_{xx}$ is negative.  This term is inversely proportional to $k(p)$ and so for the small $k(p)$ values tested, the magnitude of $\delta^H\mathcal{B}p_x^2p_{xx}$ is large.
Through the above experimentation, we have seen that counteracting $\delta\mathcal{B}^Hp_x^2p_{xx}$ locally adds diffusion in the desired shock region where the gradient has the largest magnitude.  
The MHM gives an indication of the precise amount of diffusion to add without any parameter tuning to result in a monotone temporal profile. 

We have seen that $\mathcal{B}^H$ is comprised of the above two parts: one that has a negative coefficient and the other with a positive coefficient of $p_{xx}$.  The temporal oscillations appear when the gradients are sharp enough so that the entire coefficient of $p_{xx}$ is negative at any grid point.  The one-dimensional governing equation \eqref{eq:GPME} can be written in non-conservative form as $p_t = kp_{xx} + k_xp_x$ and so the order one coefficient of $p_{xx}$ is $k(p)$.  Up to second order terms we have $[k(p) + (\mathcal{A}+\mathcal{B}^H)p_x^2]p_{xx}$.  It seems from our studies that the temporal oscillations occur when 
\begin{equation}
	k(p) + (\mathcal{A} + \mathcal{B}^H)p_x^2 < 0
	\label{eq:antidiff}
\end{equation} 
at at least one grid point near the front.  

\section{Numerical Results} 
\label{MH_Results}
In this section, we provide the numerical results comparisons of the Modified Harmonic Method (MHM) to the FTCS discretization with arithmetic and harmonic averages for the PME in Eqn. \eqref{eq:PME} and superslow diffusion in Eqn. \eqref{eq:superslow}.  We test the numerically challenging cases required for degeneracy in Section \ref{degen}, where $p<\!<1$. For very small $p$, $k(p)$ goes to zero in the test equations and the parabolic part of Eqn. \eqref{eq:k_eqtn} is masked.  Unless an exact solution is available, we use a reference solution that is computed on a mesh with $N = 3200$ grid points.  
The initial profile is chosen as the curve for $t = 0$ in Figure \ref{fig:movinginterface}, unless otherwise stated. 
The boundary conditions are the Dirichlet conditions
\[
p(0, t) = 2.0, \text{ }  p(1, t) = 0.1, \text{ } \hspace{.1cm} \forall t.\] 
We show results for $N = 50, 100, 200$ and $400$ grid points.  The convergence results are provided in \ref{conv}.
 \subsection{Porous Medium Equation (PME)}
The results for the PME with $m = 1, 2,3$ are shown in this subsection.  For Eqn. \eqref{eq:PME}, the first three derivatives of the monomial $k(p)$ with respect to $p$ that are used in the MHM are $k^{(j)}(p) = \prod_{l=0}^{j-1}(m-l)p^{m-j}|_{j=1,2,3}$. 

\subsubsection{$k(p) = p$, \hspace{.1cm}$ m = 1$}
\begin{figure}[H]
		\center
 \begin{subfigure}{.42\textwidth}  
			\includegraphics[width =\textwidth]{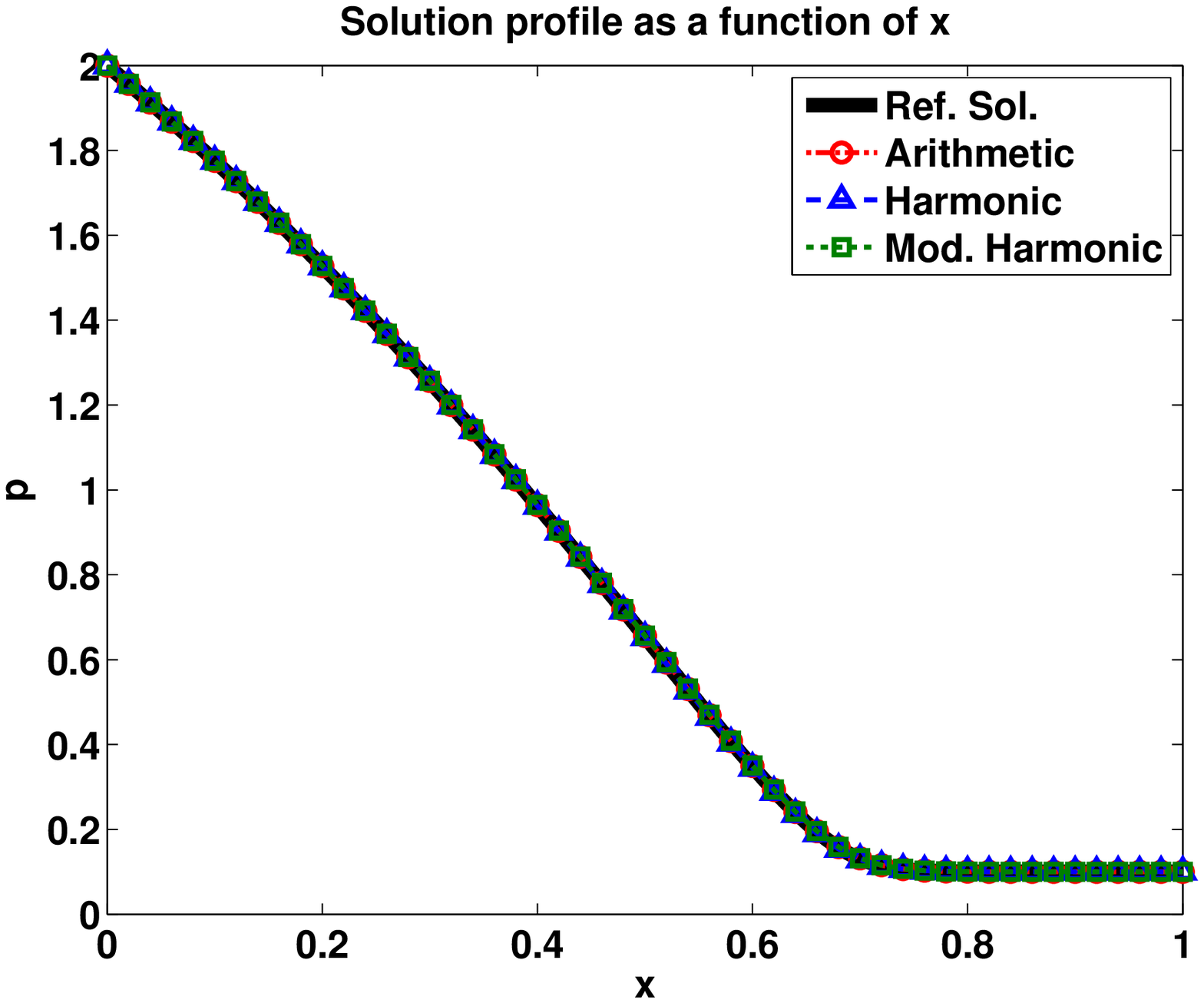}
			\caption {$t = 0.08$}
			\label{fig:m=1_x}
		\end{subfigure}
		\begin{subfigure}{0.42\textwidth}  
			\includegraphics[width =\textwidth]{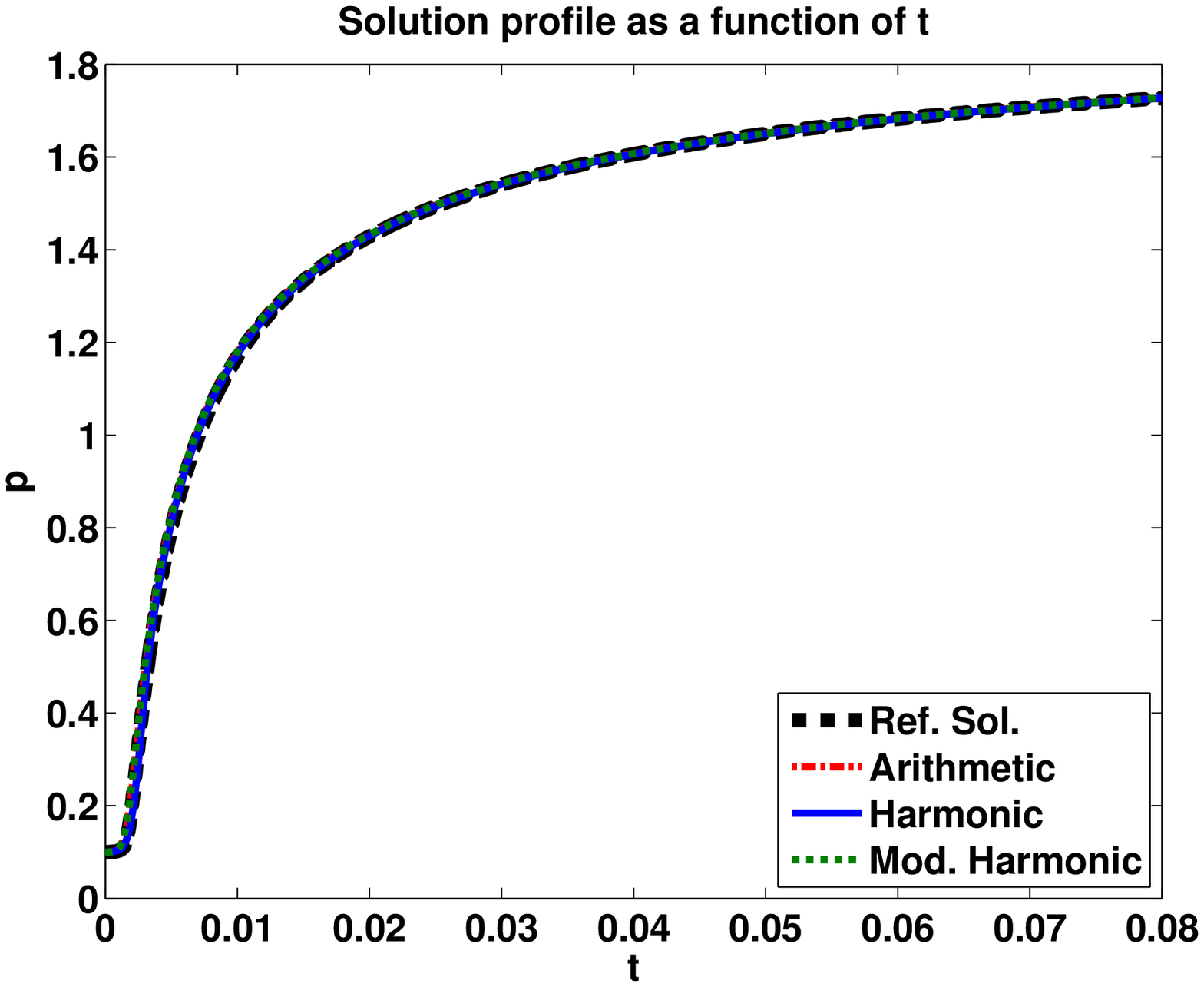}
			\caption {$x = 0.12$}
			\label{fig:m=1_t}
		\end{subfigure}
		\caption{Results for $k(p) = p$ with $N = 50$ grid points and $\dt  = \dx^2/4$.}
		\label{fig:m1}
\end{figure}

The results in Figures \ref{fig:m=1_x} and \ref{fig:m=1_t} agree with the theoretical property that the gradient is finite for $m = 1$, as discussed in Section \ref{degen}. 
For  $m = 1$, the solution is smooth and the magnitude of the anti-diffusive part $\delta\mathcal{B}^H$ of $\mathcal{B}^H$ is small.  The overall coefficient of $p_{xx}$ in Eqn. \eqref{eq:antidiff} remains positive at every grid point and so no temporal oscillations occur. 
For this smooth problem, each scheme agrees with the reference solution, also on coarser grids.  

\subsubsection{$k(p) = p^2$, \hspace{.1cm}$m = 2$}
\begin{figure}[H]
		\center
		\begin{subfigure}{.42\textwidth}  
			\includegraphics[width =\textwidth]{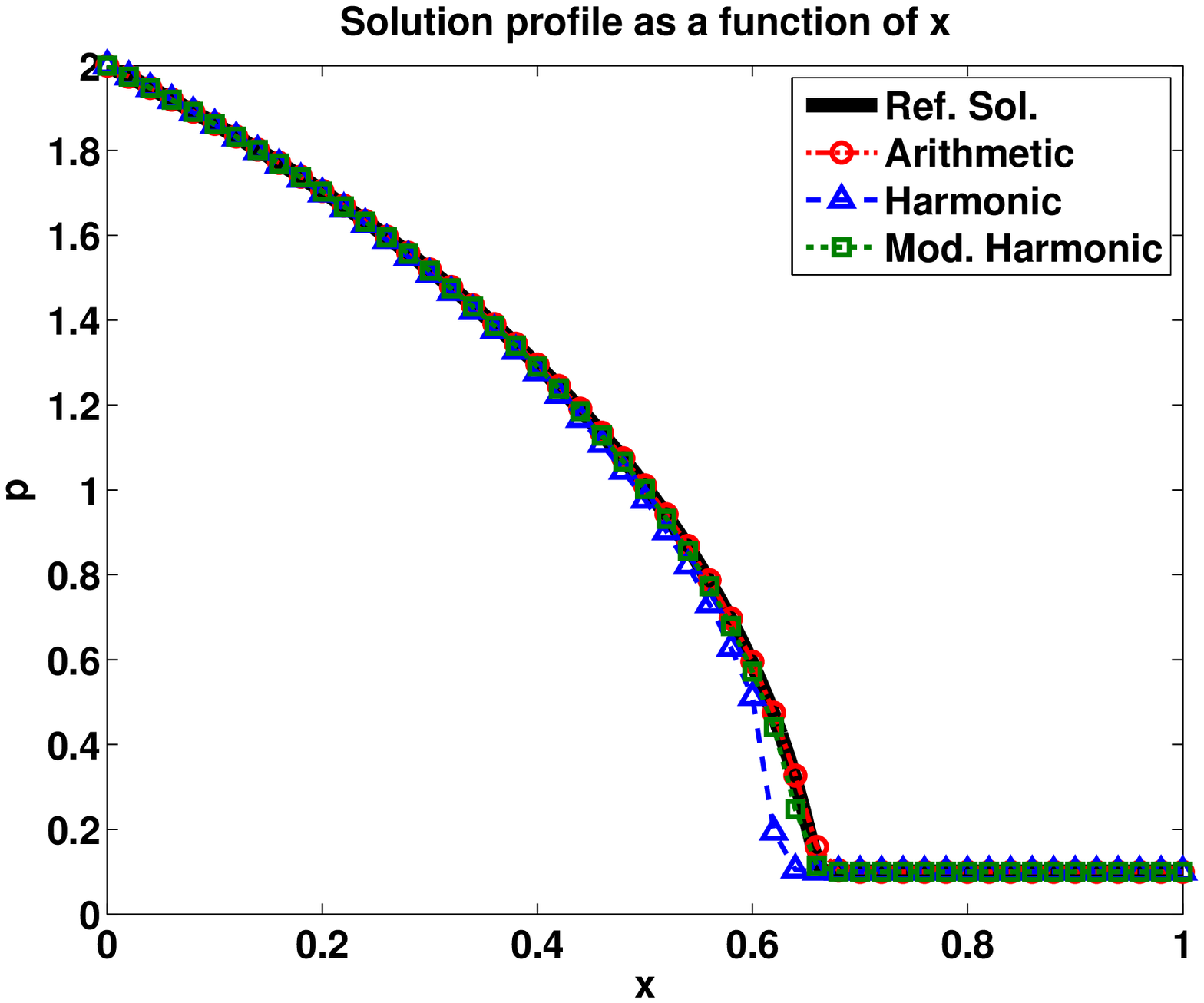}
			\caption {$t = 0.08$}
			\label{fig:m=2_x}
		\end{subfigure}
		\begin{subfigure}{0.42\textwidth}  
			\includegraphics[width =\textwidth]{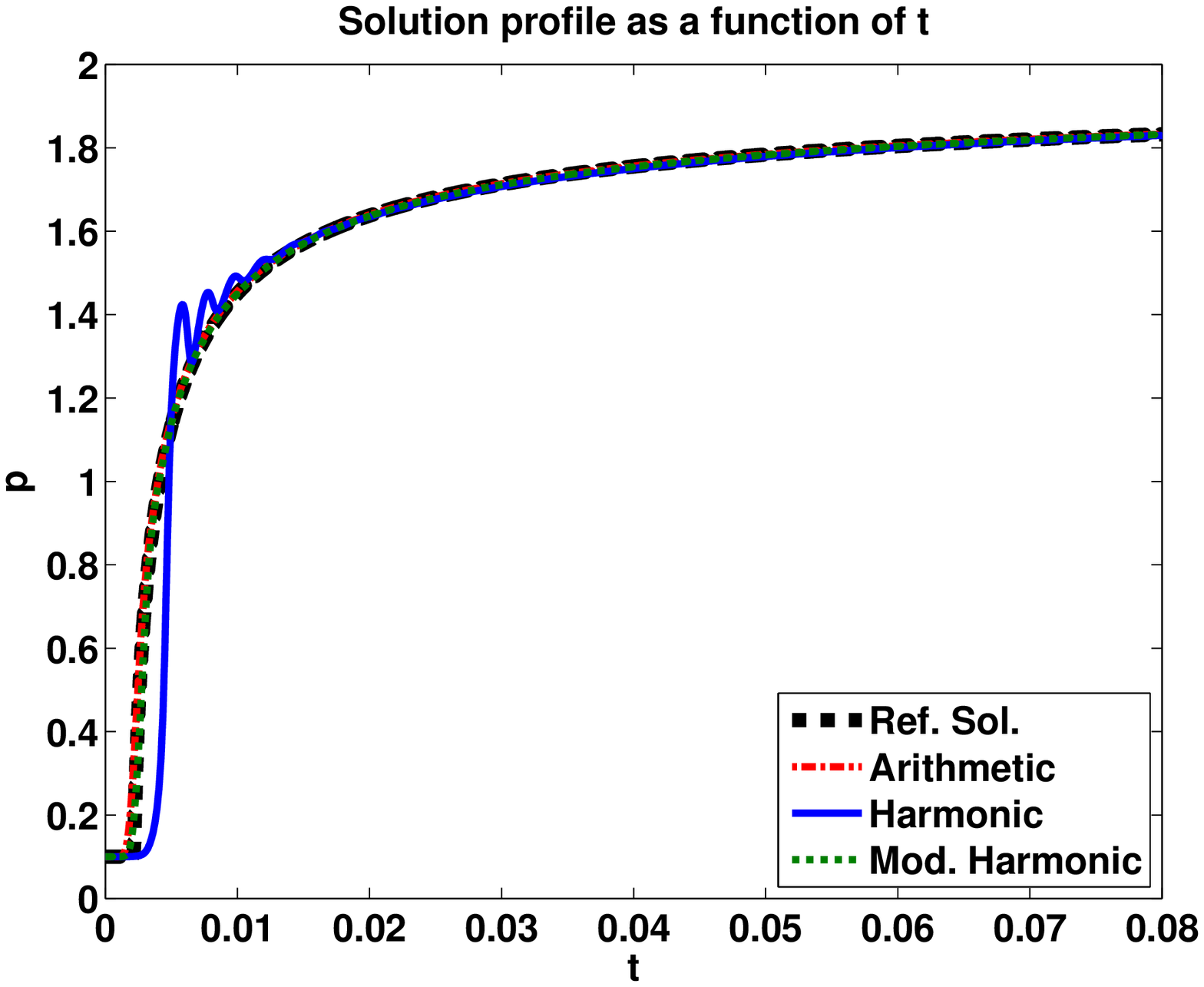}
			\caption {$x = 0.12$}
			\label{fig:m=2_t50}
		\end{subfigure}
		\caption{Results for $k(p) = p^2$ with $N = 50$ grid points and $\dt  = \dx^2/8$.}
\end{figure}

\begin{figure}[H]
		\center
		\begin{subfigure}{.42\textwidth}  
			\includegraphics[width =\textwidth]{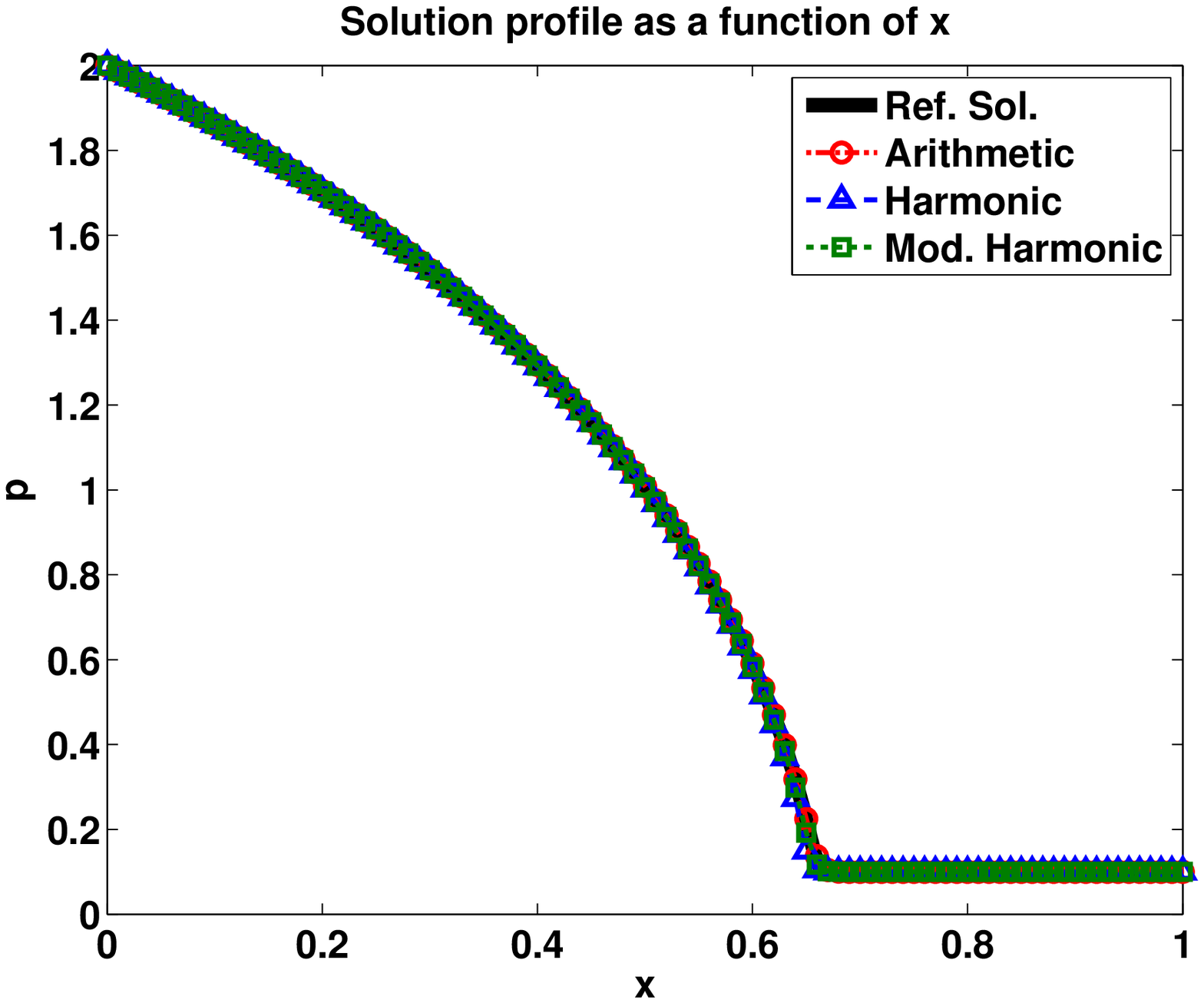}
			\caption {$t = 0.08$}
			\label{fig:m=2100}
		\end{subfigure}
		\begin{subfigure}{0.42\textwidth}  
			\includegraphics[width =\textwidth]{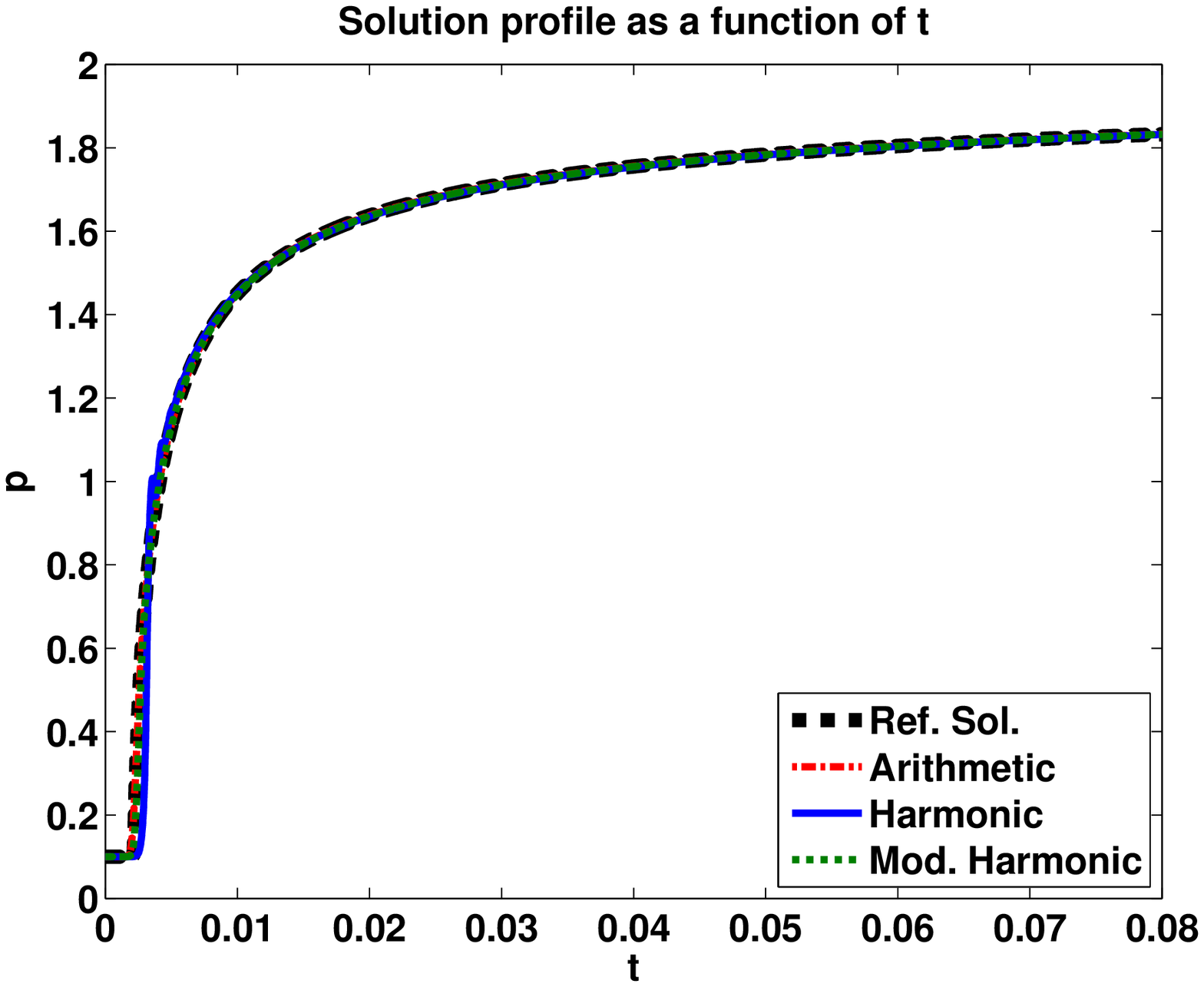}
			\caption {$x = 0.12$}
			\label{fig:m=2_t100}
		\end{subfigure}
		\caption{Results for $k(p) = p^2$ with $N = 100$ grid points and $\dt  = \dx^2/8$.}
		\label{fig:m=2}
\end{figure}

For $m > 1$, sharp gradients and decreased $k(p)$ values with increased $m$ contribute to the numerical artifacts present in the solution with harmonic averaging.  On the coarser grid in Figure \ref{fig:m=2_x}, the numerical solution with harmonic averaging is slightly lagging behind the reference solution.  Figures \ref{fig:m=2_t50} shows temporal oscillations for this method.  Refining the mesh from $N = 50$ to $N = 100$ grid points reduces the lagging of the shock position in Figure  \ref{fig:m=2100}, but lower amplitude temporal oscillations are still present in Figure  \ref{fig:m=2_t100}.

For this problem, $\mathcal{B^H} =  -3\dx^2/2$ and $\mathcal{F^H} = 0$, and so Eqn. \eqref{eq:MHM} simplifies to $(3\dx^2/2)p_x^2p_{xx}$.  The negative value of $\mathcal{B^H}$ clearly shows the aforementioned anti-diffusion.  Figures \ref{fig:m=2100} and \ref{fig:m=2_t100} illustrate that adding the above partially advective and diffusive term in the MHM improves the shock position and results in a smooth, monotone temporal profile for both grids. 
\justify
\vspace{-.4cm}

\subsubsection{$k(p) = p^3$, \hspace{.05cm} $m = 3$}
\begin{figure}[H]
		\center
		\begin{subfigure}[H]{.42\textwidth}  
			\includegraphics[width =\textwidth]{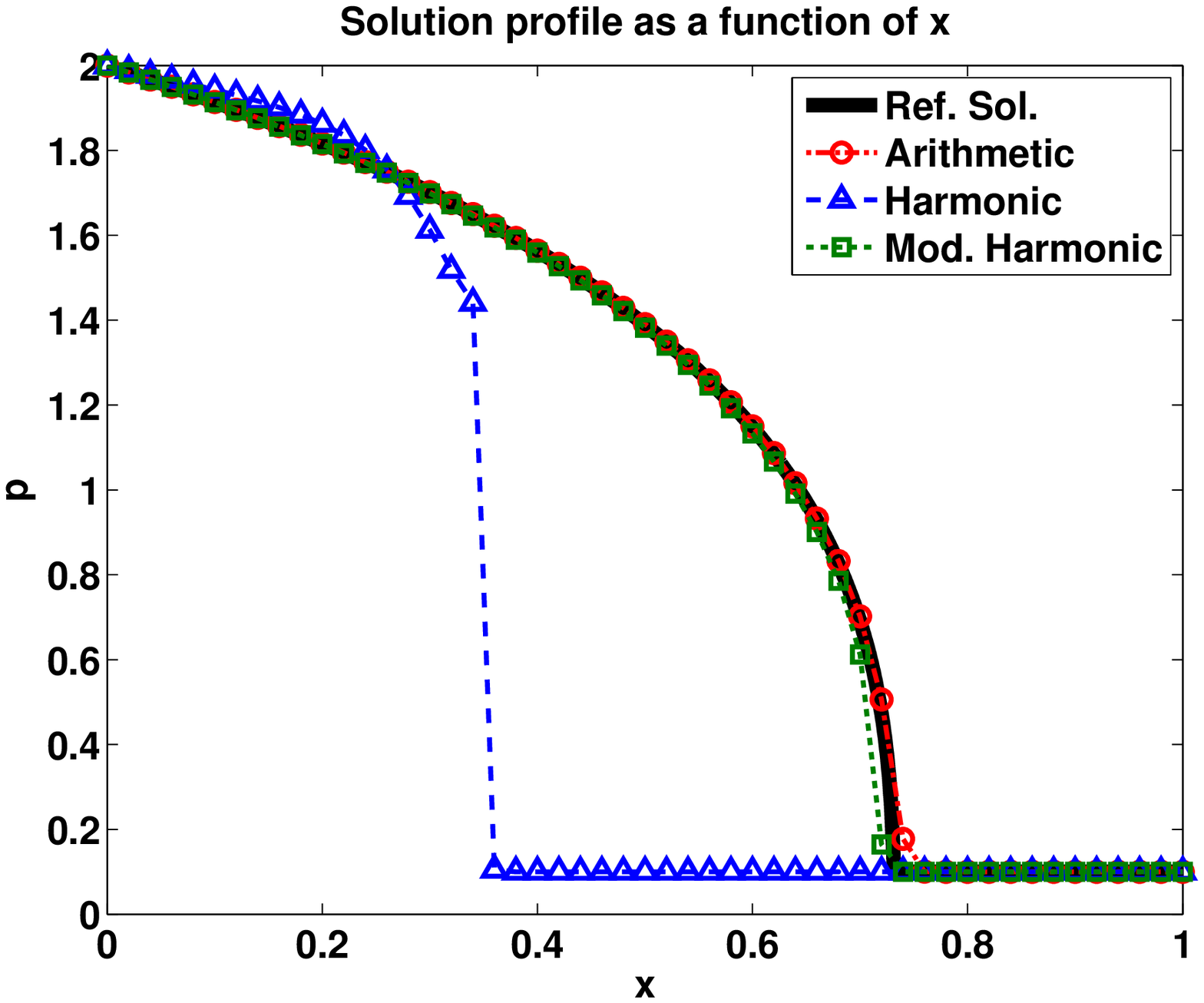}
			\caption {$t = 0.08$}
			\label{fig:m=3_x50}
		\end{subfigure}
		\begin{subfigure}[H]{0.42\textwidth}  
			\includegraphics[width =\textwidth]{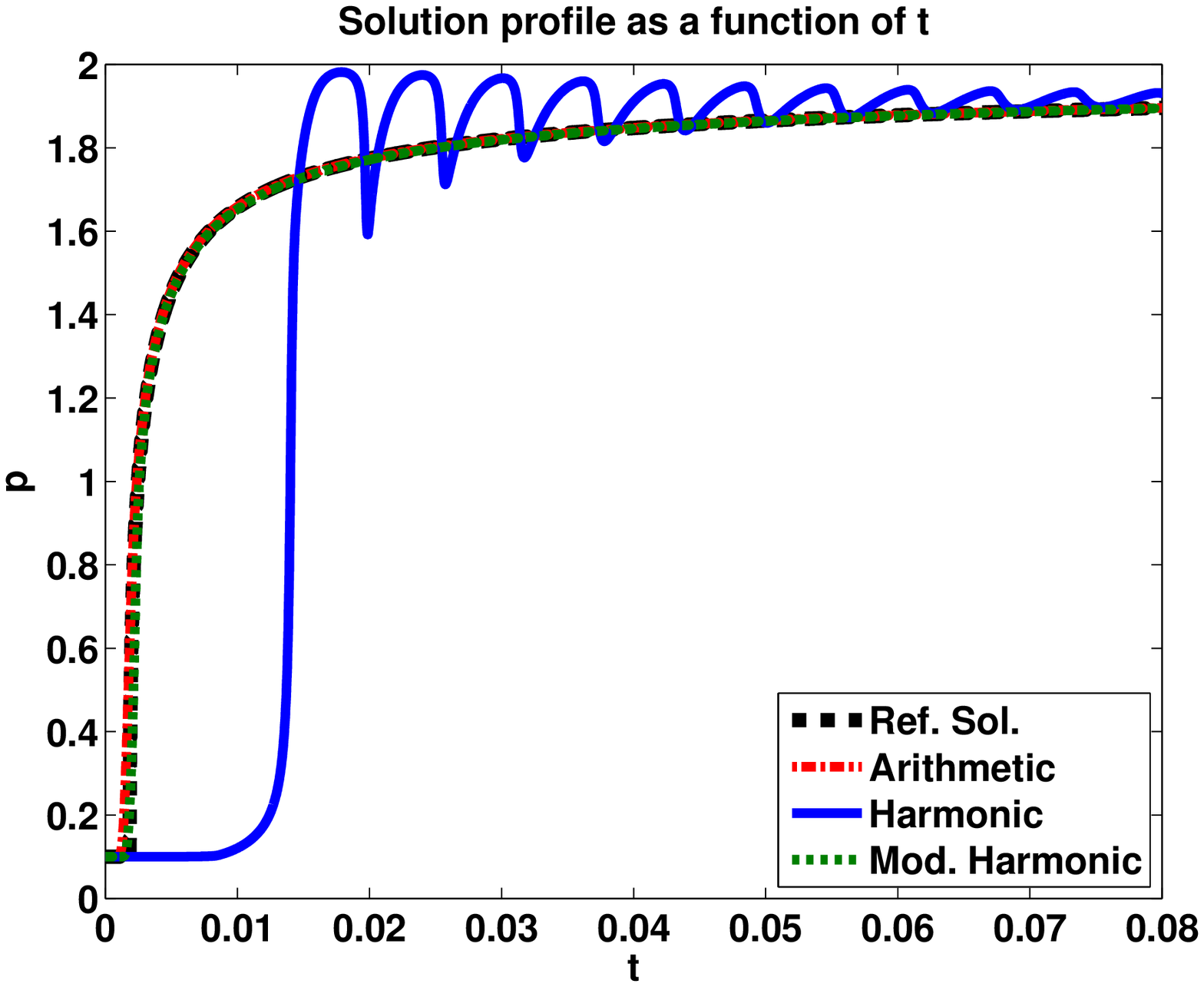}
			\caption {$x = 0.12$}
			\label{fig:m=3_t50}
		\end{subfigure}
		\caption{Results for $k(p) = p^3$ with $N = 50$ grid points and $\dt  = \dx^2/16$.}
		\label{fig:m=350}
\end{figure}

\begin{figure}[H]
		\center
		\begin{subfigure}[h]{.42\textwidth}  
			\includegraphics[width =\textwidth]{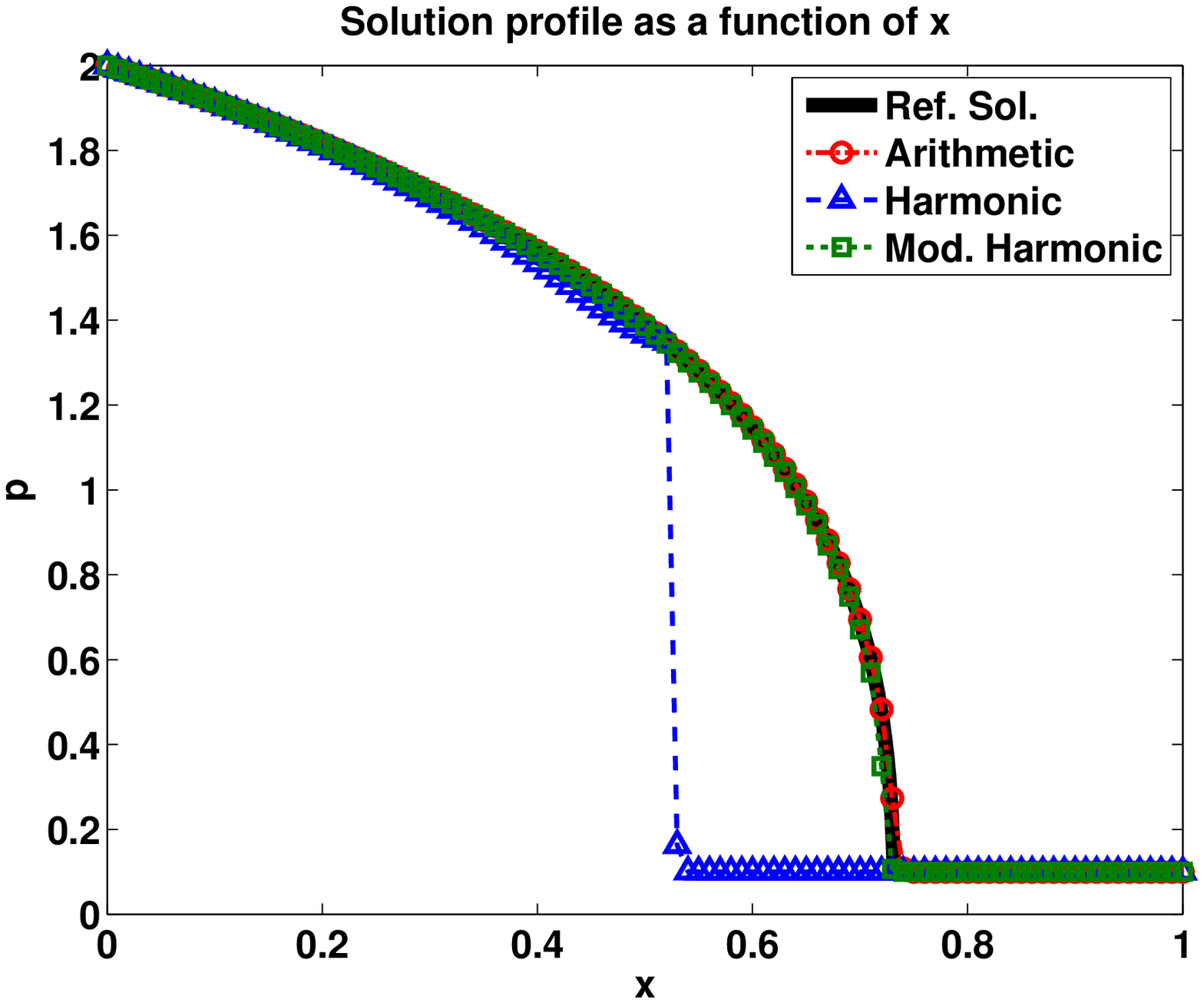}
			\caption {$t = 0.08$}
			\label{fig:m=3_x100}
		\end{subfigure}
		\begin{subfigure}[h]{0.42\textwidth}  
			\includegraphics[width =\textwidth]{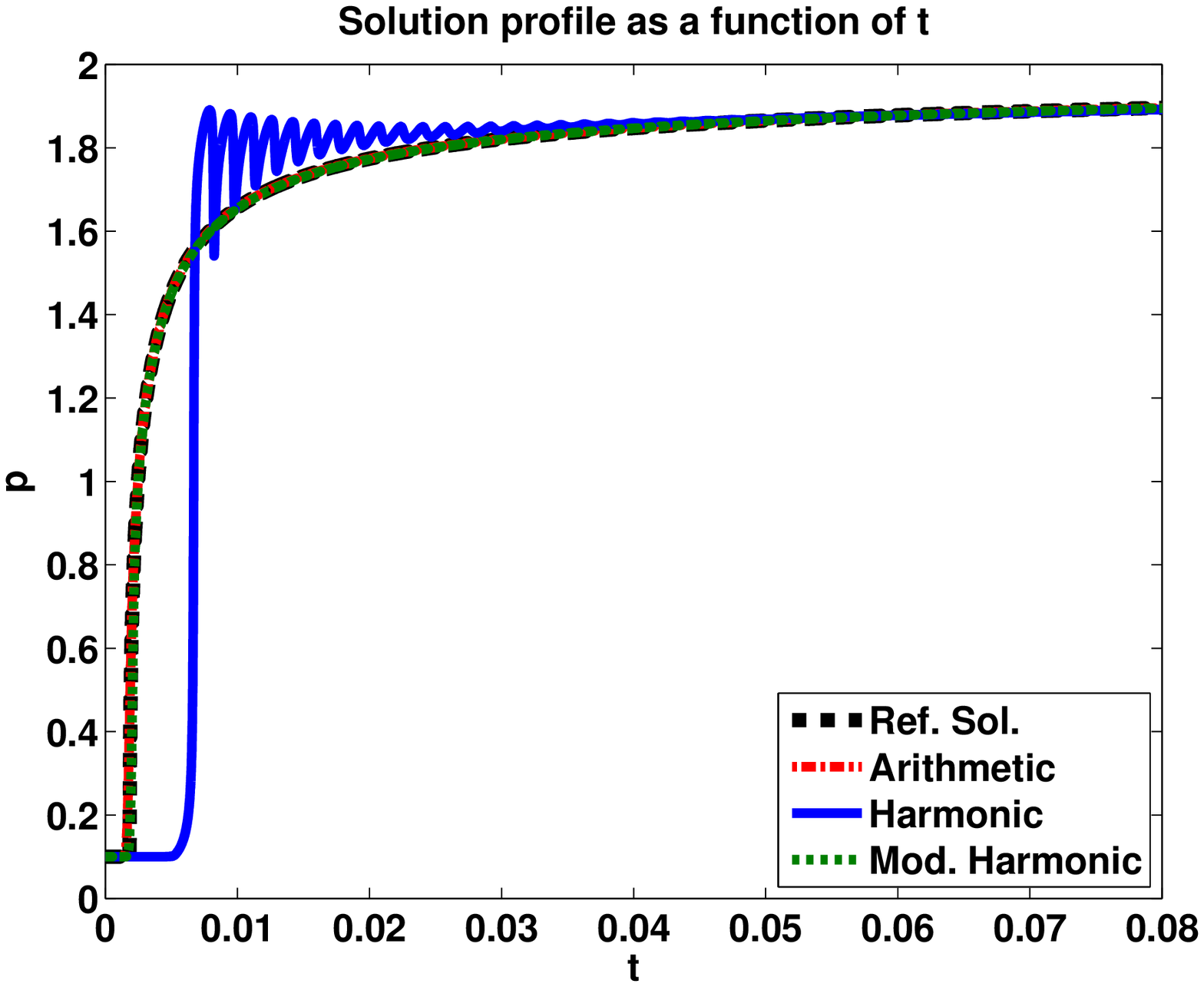}
			\caption {$x = 0.12$}
			\label{fig:m=3_t100}
		\end{subfigure}
		\caption{Results for $k(p) = p^3$ with $N = 100$ grid points with $\dt  = \dx^2/16$. }
		\label{fig:m=3100}
\end{figure}

\begin{figure}[H]
		\center
		\begin{subfigure}[H]{.42\textwidth}  
			\includegraphics[width =\textwidth]{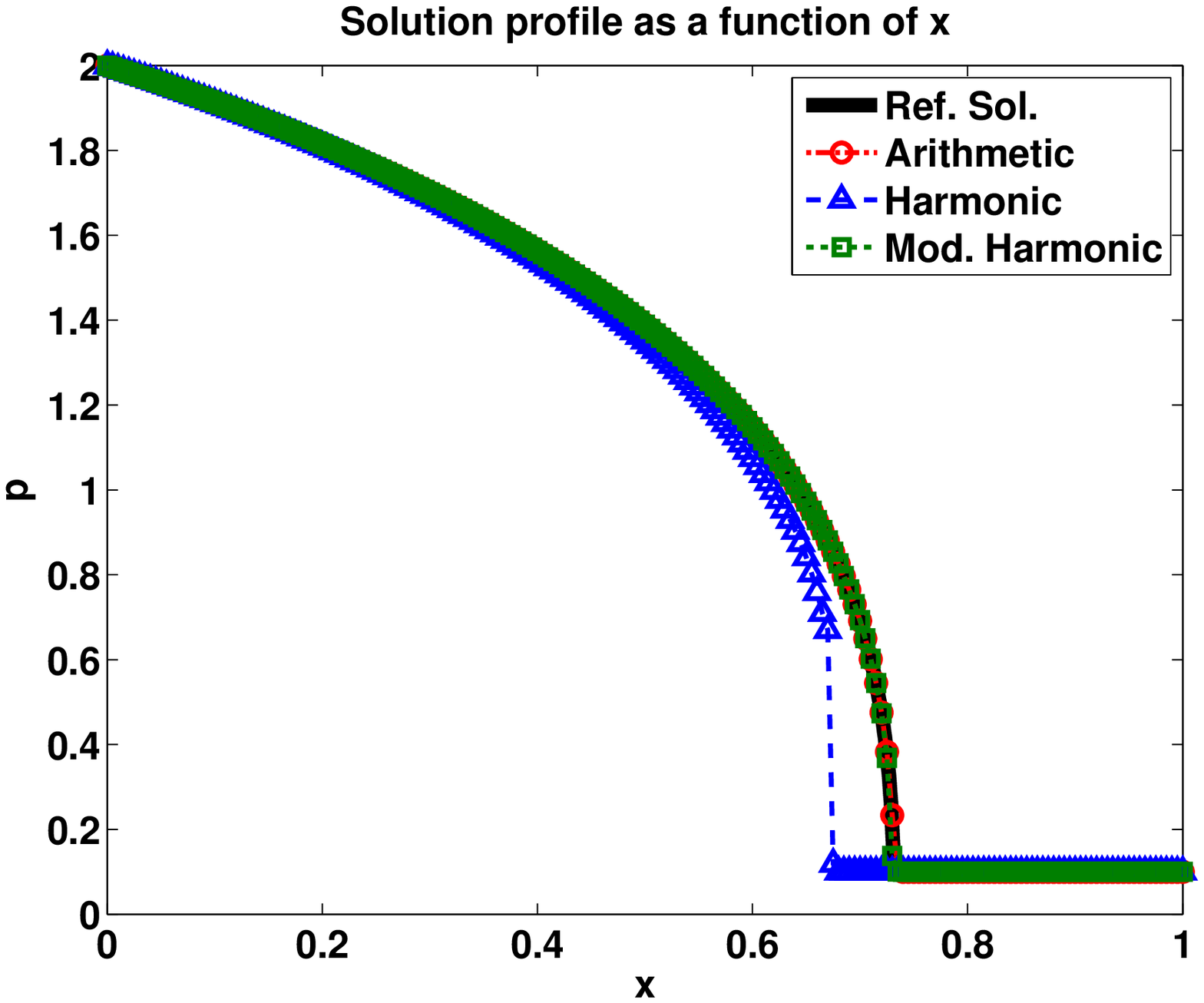}
			\caption {$t = 0.08$}
			\label{fig:m=3_x200}
		\end{subfigure}
		\begin{subfigure}[H]{0.42\textwidth}  
			\includegraphics[width =\textwidth]{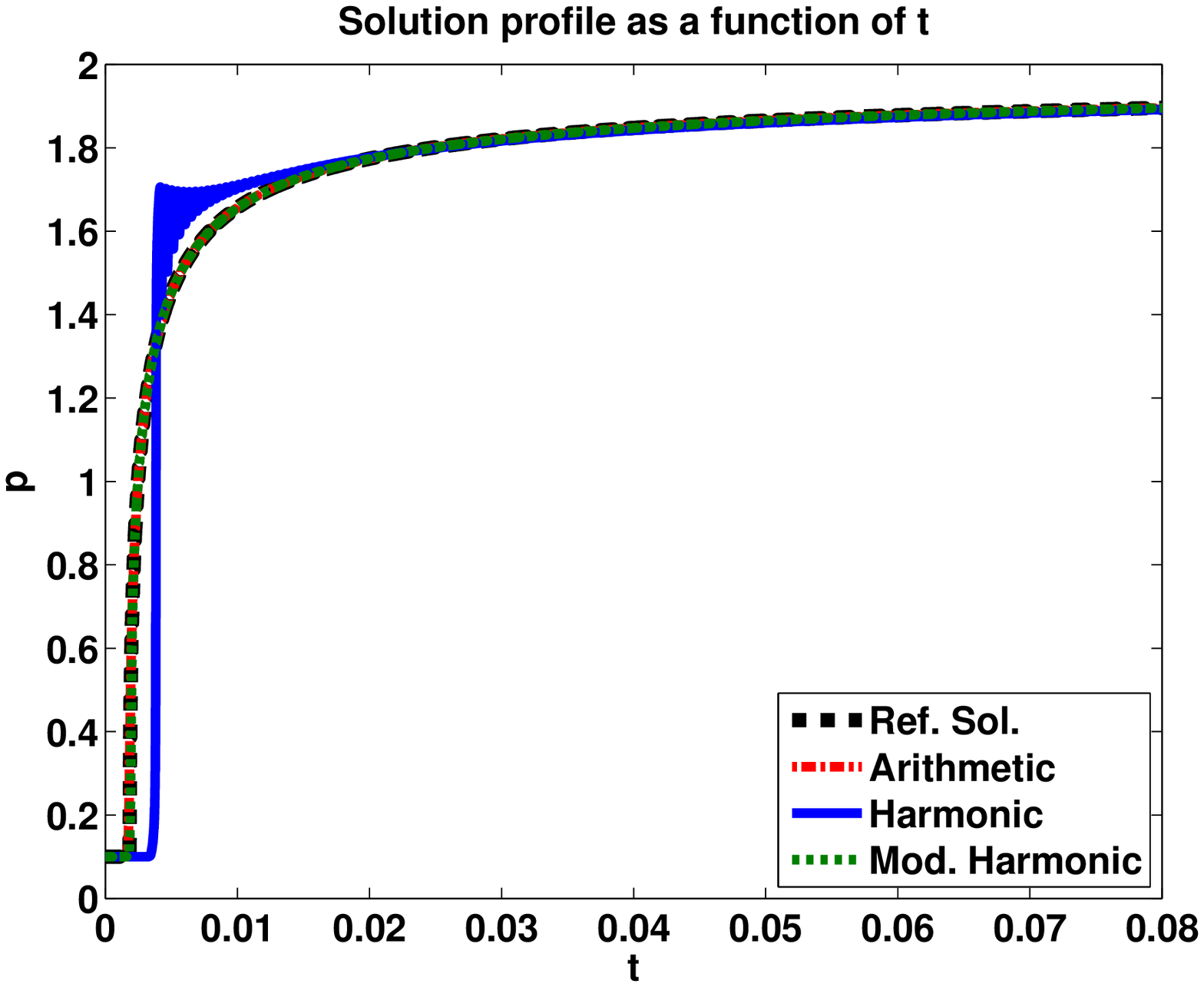}
			\caption {$x = 0.12$}
			\label{fig:m=3_t200}
		\end{subfigure}
		\caption{Results for $k(p) = p^3$ with $N = 200$ grid points with $\dt  = \dx^2/16$.}
		\label{fig:m=3200}
\end{figure}

 As discussed in Section \ref{degen}, we expect more problems for $m = 3$ than $m =1$ or 2 since for the same small $p$ value, $k(p)$ is decreased. 
 Figures \ref{fig:m=3_x50} and \ref{fig:m=3_x100} show the lagging in the numerical solution with harmonic averaging.  This still holds even on the finer grid of $N = 200$ grid points in Figure \ref{fig:m=3_x200}.
For $m = 3$, sharp gradients occur at the front, causing the overall coefficient of $p_{xx}$ in Eqn. \eqref{eq:antidiff} to become negative at one or two grid points.  This negative diffusive coefficient results in temporal oscillations.  Figures \ref{fig:m=3_t50}, \ref{fig:m=3_t100} and \ref{fig:m=3_t200} illustrate that the period and the amplitude decrease as the grid is refined.  

In this case, $\mathcal{B}^H = -(9\dx^2/4)p$ and $\mathcal{F}^H = -5\dx^2/4$, and so Eqn. \eqref{eq:MHM} simplifies to $(\dx^2/4)[9pp_x^2p_{xx}+5p_x^4]$, revealing its locally diffusive and advective parts. The anti-diffusive effect of $\mathcal{B}^H$ in its simplified form is again clear. The improvement with the MHM is apparent in the subsequent figures with the correct shock position and smooth temporal results. 



\vspace{-0.575cm}
\vspace{-0.575cm}
\begin{figure}[H]
		\center
		\begin{subfigure}[H]{.42\textwidth}  
			\includegraphics[width =\textwidth]{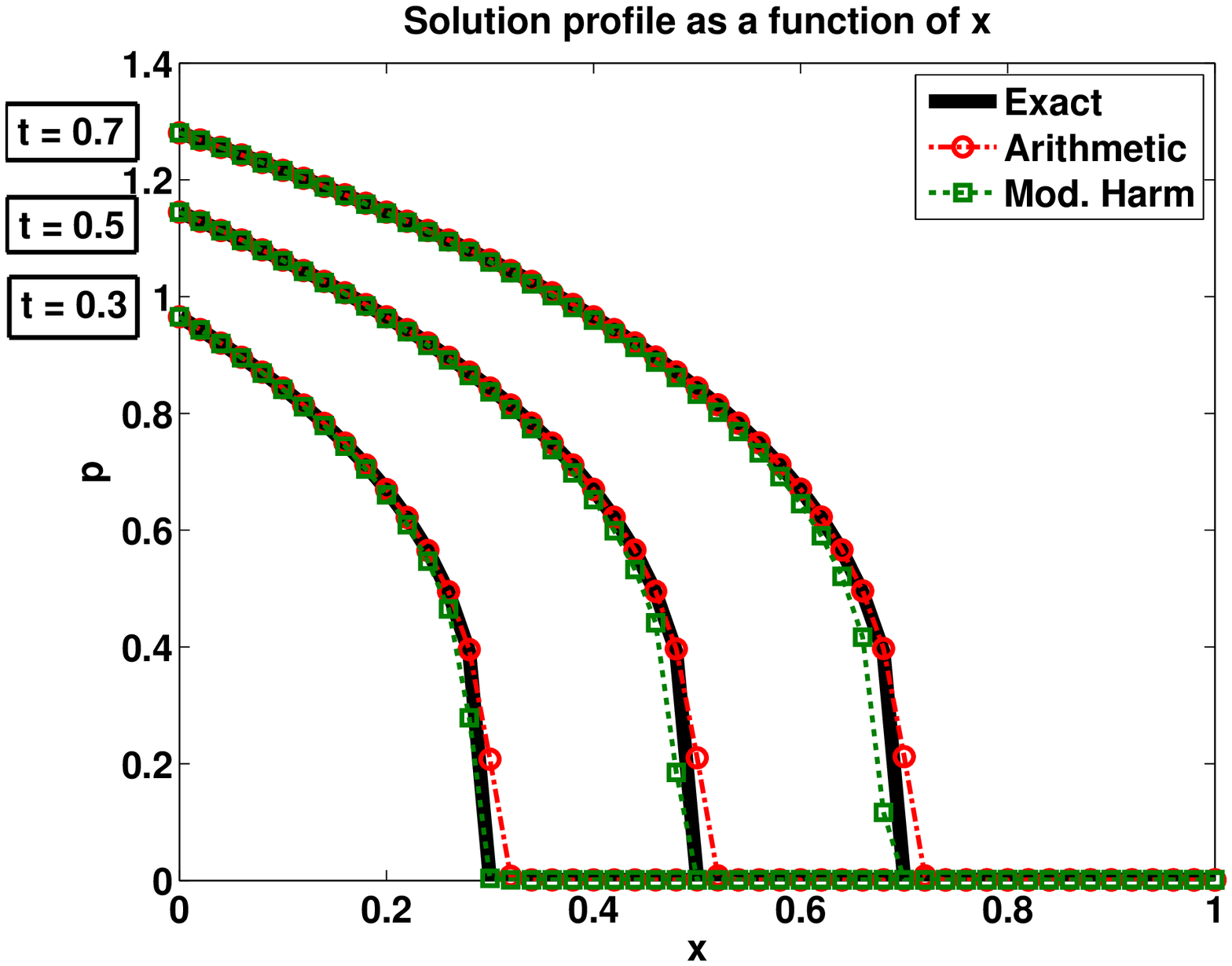}
			\caption {$t = 0.3, 0.5, 0.7$}
			\label{fig:lock_x}
		\end{subfigure}
		\begin{subfigure}[H]{0.41\textwidth}  
			\includegraphics[width =\textwidth]{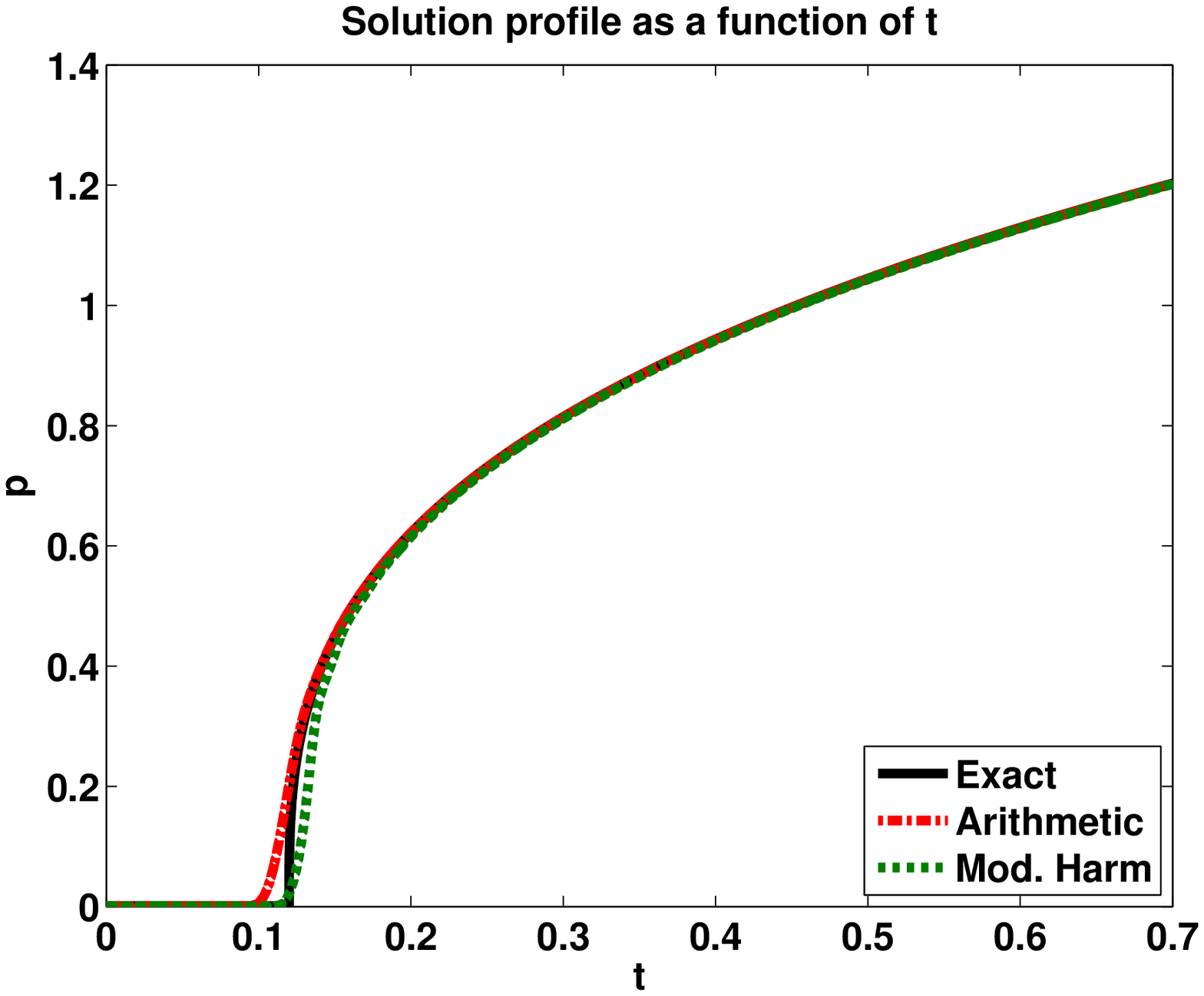}
			\caption {$x = 0.12$}
			\label{fig:m=lock_t}
		\end{subfigure}
		\label{fig}
		\caption{Results for $k(p) = p^3$ with $N = 50$ grid points and $\dt  = \dx^2/16$. The results with harmonic averaging are not shown, since this numerical solution locks as shown in the TLP in Section \ref{TLP}.
}
		\label{fig:mim_lock}
\end{figure}

In the aforementioned examples, we have seen that the MHM resolves the lagging and temporal oscillations associated with harmonic averaging for small $p$ values.  To test the effect of the MHM on the extreme case of lagging, we now consider the locking problem (TLP) from Section \ref{TLP} with an exact solution given in Eqn. \eqref{eq:mim_exact}.  Figure \ref{fig:lock_x} shows that with the MHM, the solution no longer locks, and that its shock position is correct.  Figure \ref{fig:lock_x} also reveals that the numerical solution with arithmetic averaging is slightly ahead of the true shock position.  In Figure \ref{fig:m=lock_t}, we see that the MHM solution has a smooth temporal profile aligning with the true solution.  Even in extreme cases, we find that correcting both terms in the MHM is sufficient.

\subsection{ Superslow diffusion: $k(p) = \exp(-1/p)$}
\begin{figure}[H]
		\center
		\begin{subfigure}{.42\textwidth}  
			\includegraphics[width =\textwidth]{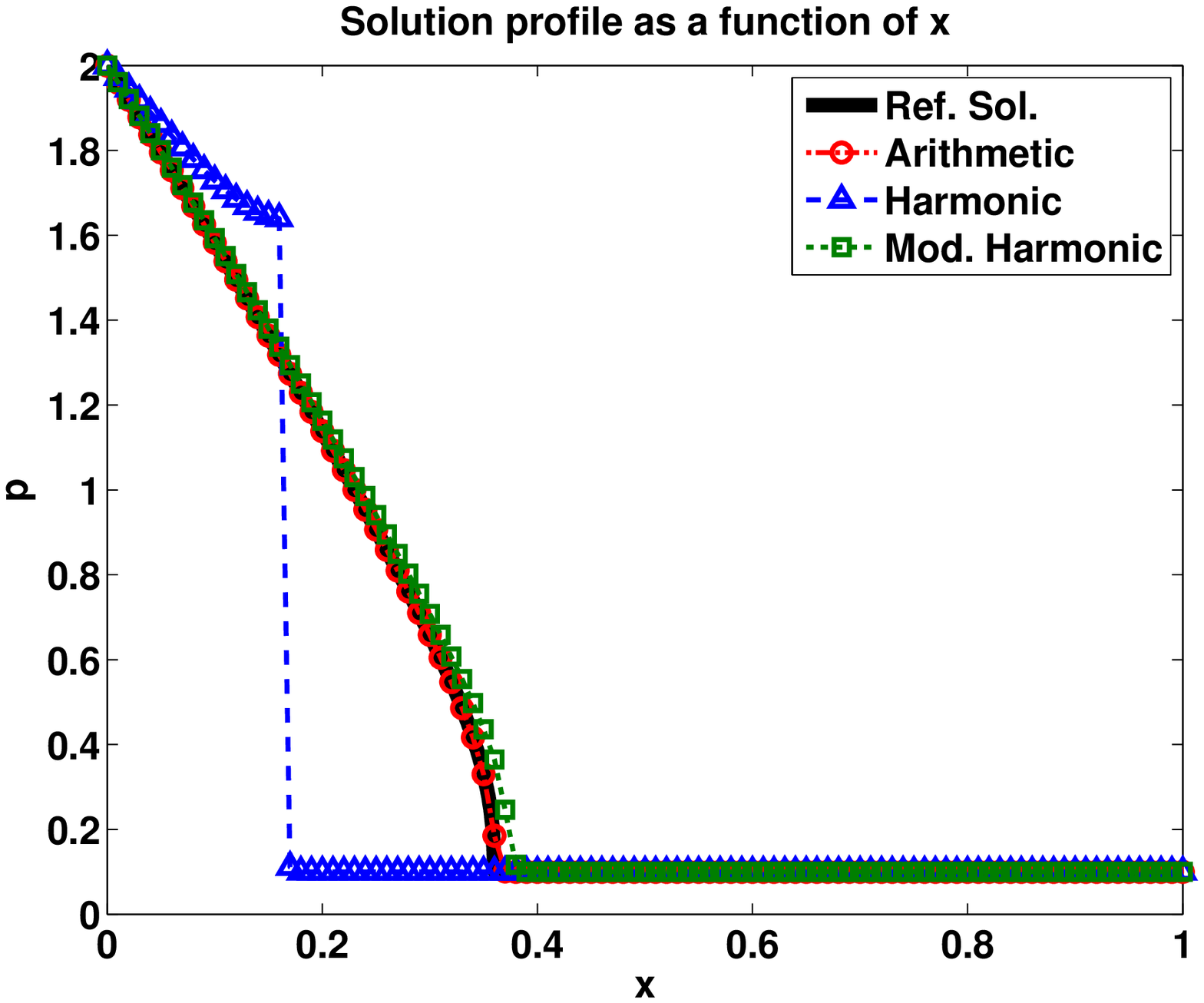}
			\caption {$t = 0.08$}
			\label{sd_100x}
		\end{subfigure}
		\begin{subfigure}{0.42\textwidth}  
			\includegraphics[width =\textwidth]{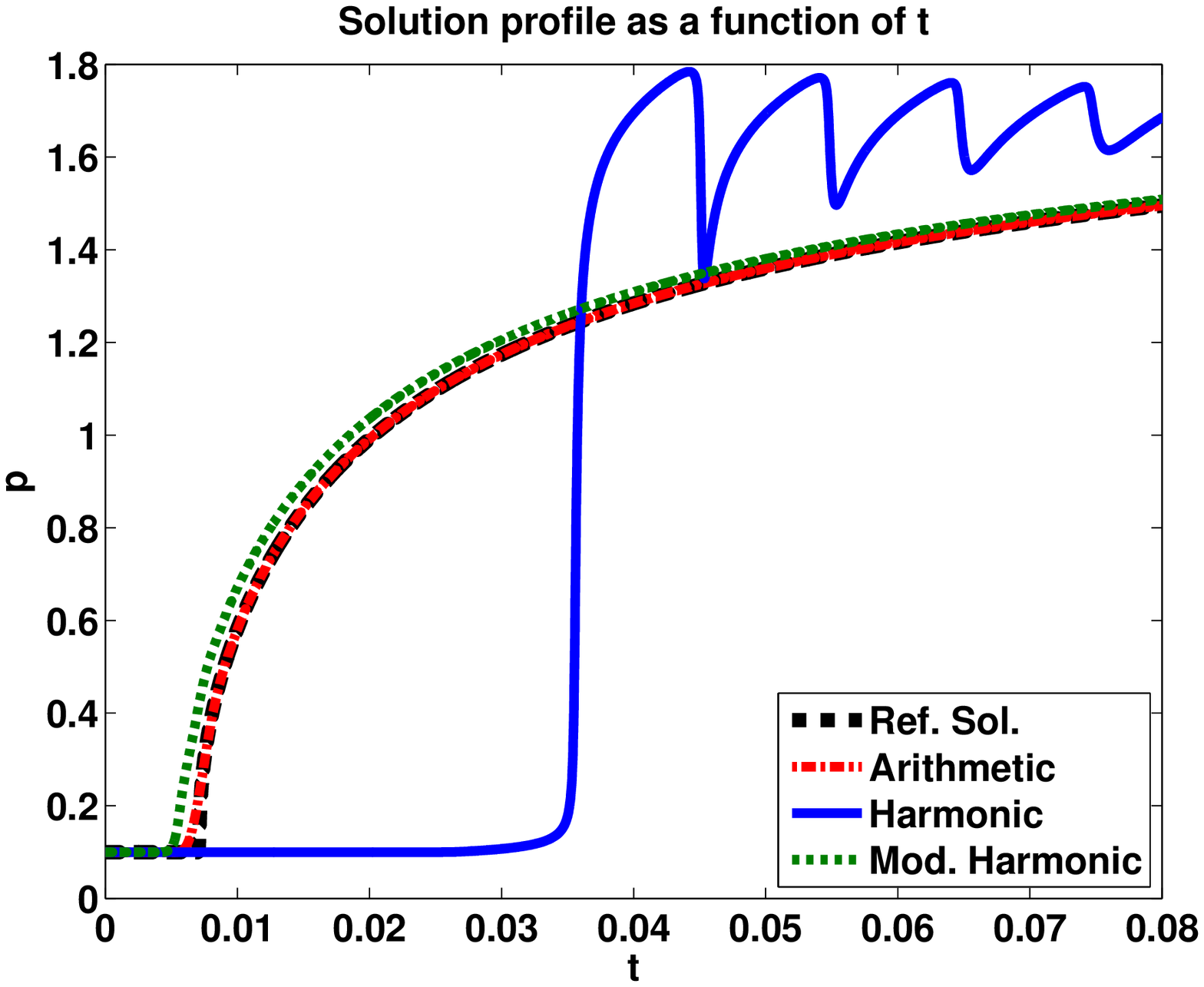}
			\caption {$x = 0.12$}
			\label{sd_100t}
		\end{subfigure}
		\caption{Results for $k(p) = \exp(-1/p)$ with $N = 100$ grid points and $\dt  = \dx^2/2$.}
\end{figure}

\begin{figure}[H]
		\center
		\begin{subfigure}{.42\textwidth}  
			\includegraphics[width =\textwidth]{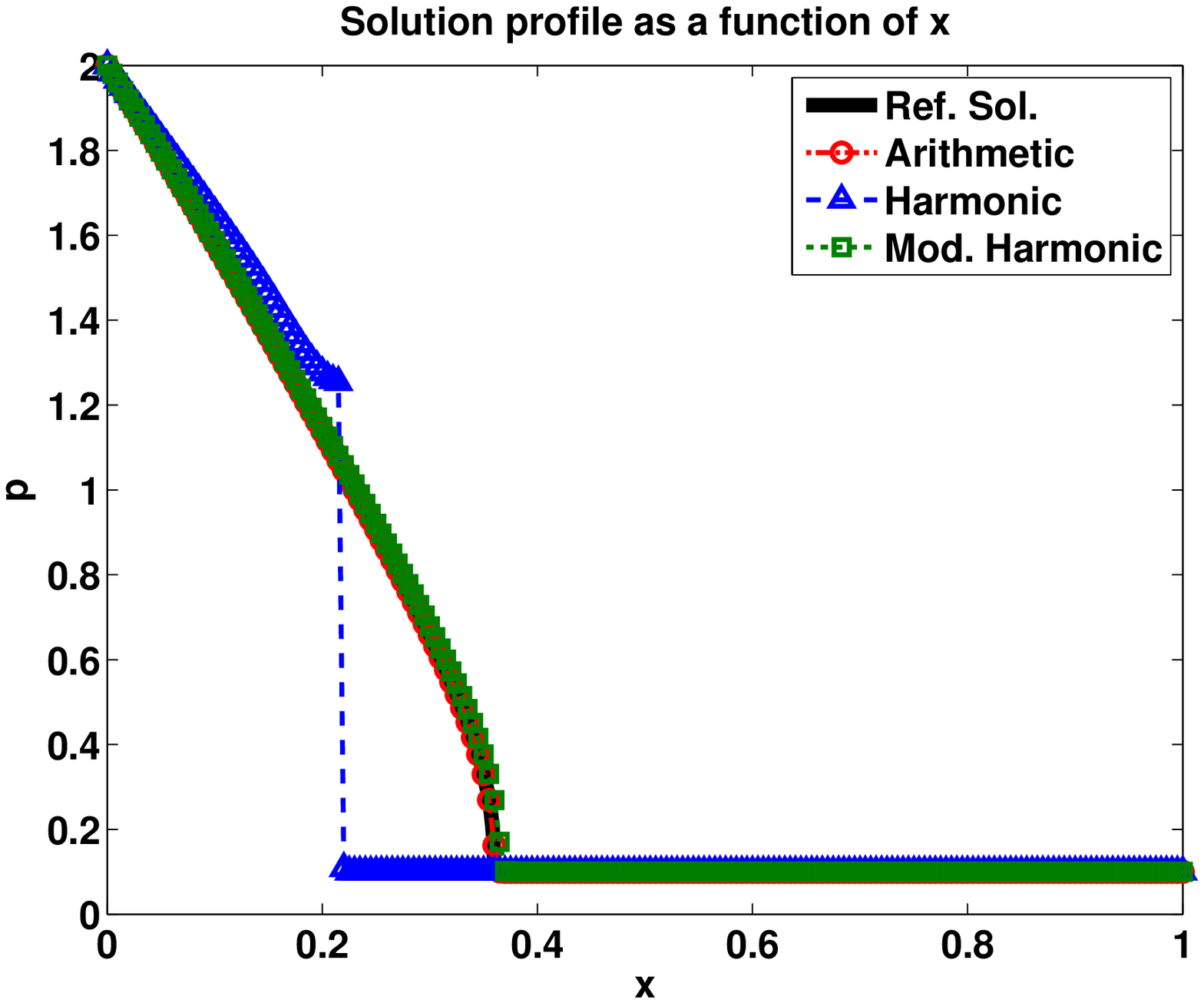}
			\caption {$t = 0.08$}
			\label{sd_200x}
		\end{subfigure}
		\begin{subfigure}{0.42\textwidth}  
			\includegraphics[width =\textwidth]{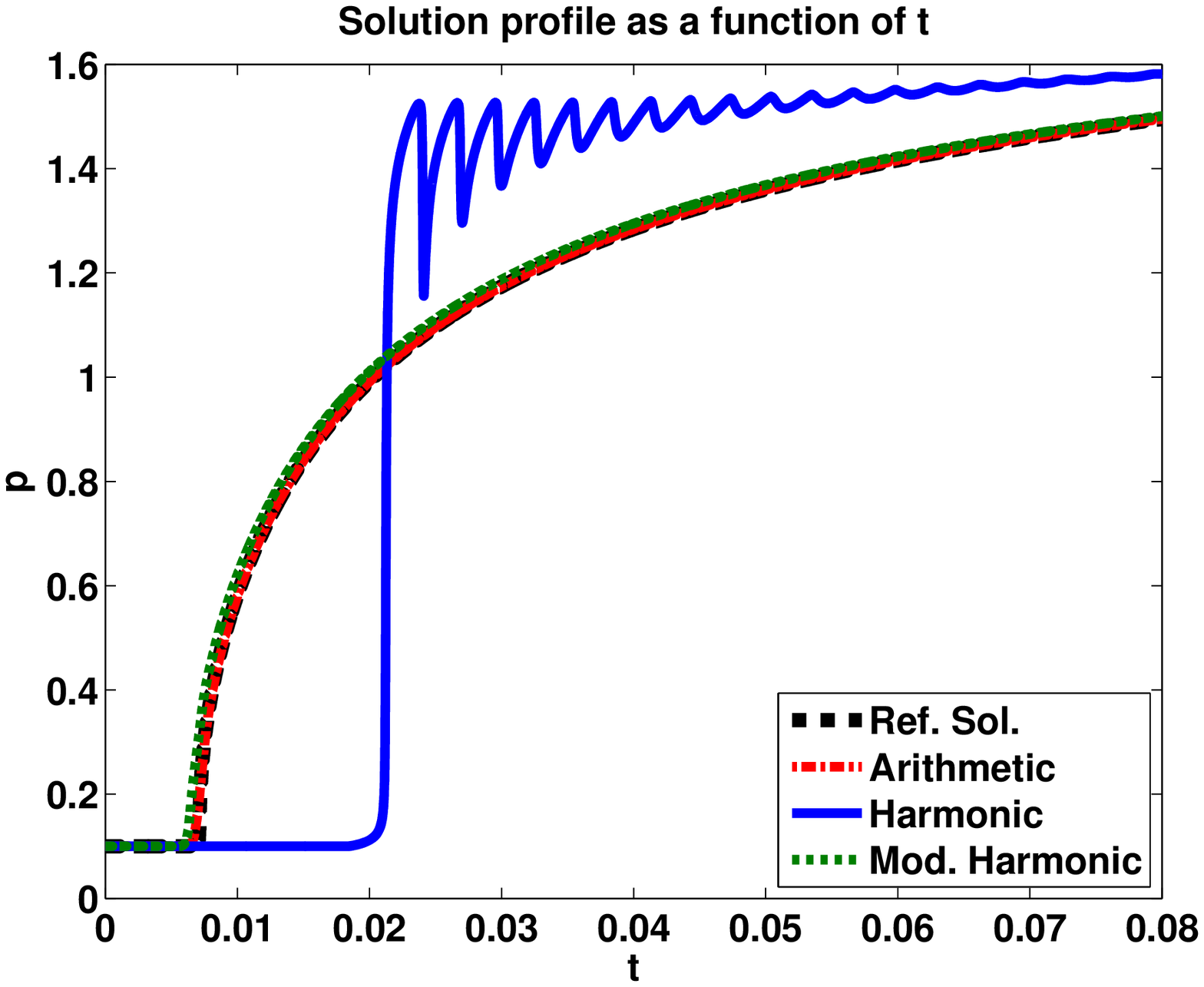}
			\caption {$x = 0.12$}
			\label{sd_200t}
		\end{subfigure}
		\caption{Results for $k(p) = \exp(-1/p)$ with $N = 200$ grid points and $\dt  = \dx^2/2$.}
\end{figure}

\begin{figure}[H]
		\center
		\begin{subfigure}{.42\textwidth}  
			\includegraphics[width =\textwidth]{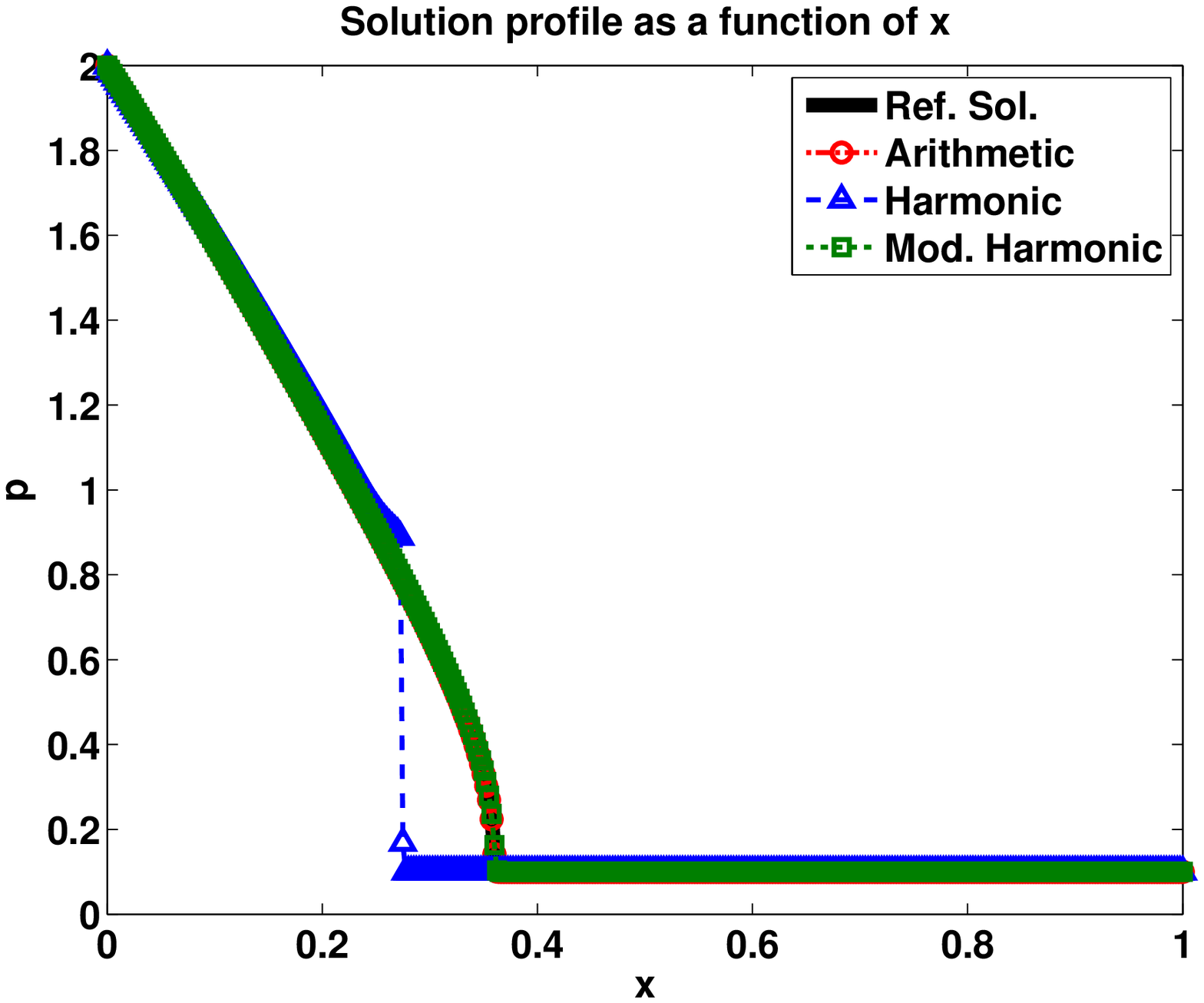}
			\caption {$t = 0.08$}
			\label{sd_400x}
		\end{subfigure}
		\begin{subfigure}{0.42\textwidth}  
			\includegraphics[width =\textwidth]{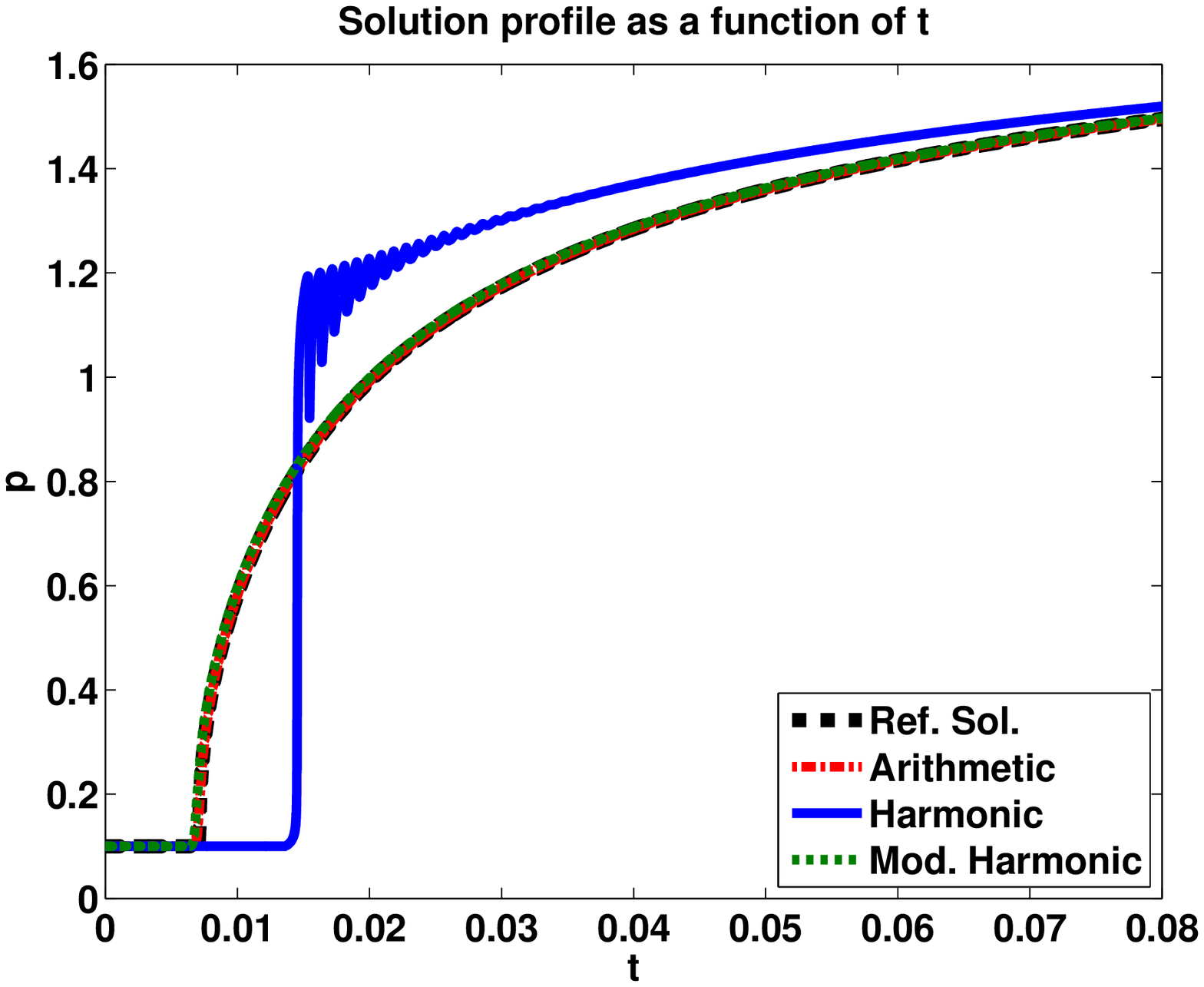}
			\caption {$x = 0.12$}
			\label{sd_400t}
		\end{subfigure}
		\caption{Results for $k(p) = \exp(-1/p)$ with $N = 400$ grid points and $\dt  = \dx^2/2$.}
		\label{fig:superslow}
\end{figure}

Superslow diffusion is the most challenging case studied, since the negative exponential goes to zero faster than any monomial of $p$ for small $p$.  Not surprisingly, Figures \ref{sd_100x} and \ref{sd_200x} show significant lagging, and Figures \ref{sd_100t} and \ref{sd_200t} reveal temporal oscillations in the numerical solution with harmonic averaging. Figures \ref{sd_400x} and \ref{sd_400t} illustrate that even on a fine mesh of $N = 400$ grid points, the numerical solution with harmonic averaging is lagging and high frequency temporal oscillations are present. 

		
To form $\mathcal{B^H}$ and $\mathcal{F}^H$ for this problem, it can be verified that $k_p = kp^{-2}, k_{pp} = k(p^{-4} - 2p^{-3})$ and $k_{ppp} = k(p^{-6} - 6p^{-5} +6p^{-4})$.  
The following figures show that the MHM is close to the true solution, and no oscillations occur.



\section{Conclusions and Further Work} 
We have seen that solving the Generalized Porous Medium Equation (GPME) with near zero coefficients poses interesting numerical challenges.  Lagging, temporal oscillations and locking can occur in the second-order finite volume discretization with harmonic averaging.  The Modified Equation Analysis is utilized to understand the cause of these artifacts in the commonly used FTCS discretization with harmonic averaging and to identify the numerically challenging terms in the truncation error.  These leading error terms are counteracted in a Lax-Wendroff like approach, by updating the discretization with harmonic averaging.  Counteracting either term individually is not enough to correct for the numerical artifacts.  Our Modified Harmonic Method (MHM) demonstrates that coupling the mitigation of the critical terms removes these artifacts.  

The numerical simulations provide confidence in the MHM.  The removal of the numerical artifacts in the results demonstrates that it is sufficient to only consider the different leading order terms in the modified equations of the arithmetic and harmonic average.  The spatial results illustrate how the MHM corrects the lagging or locking in the shock position present in the scheme with harmonic averaging.  The temporal results show that the MHM provides local diffusion at the front, resulting in smooth temporal profiles.   
The MHM gives the precise amount of artificial diffusion to add, without having to tune the parameter or to overly smooth the shock. 
 
 The MHM can easily be applied to the GPME for arbitrary differentiable $k(p)$, by computing the first three derivatives of $k$ with respect to $p$.  We can also consider the more general case, where the coefficient $k(p)$ is no longer differentiable with respect to $p$. For this type of GPME as seen in Eqn. \eqref{eq:foam}, the Modified Equation Analysis does not hold.  This problem is more challenging because both schemes with arithmetic and harmonic averaging lead to temporal oscillations.  An alternative approach is proposed in our upcoming paper \cite{maddix_foam}.  

For future work, this paper provides a general framework to extend the MHM approach to other discretizations and to higher dimensions.  
The Modified Equation Analysis for Backward Euler is discussed, and a MHM can be developed for that temporal scheme.  For TVD RK2 and other more complex temporal discretizations, the Modified Equation Analysis is more challenging and developing an analogous MHM for these discretizations is a future direction. 
In addition to different temporal discretizations, the MHM framework can also be developed for other averages, such as the geometric or integral average, as well as averages over more than two neighbors.
Since the core of the approach is Taylor expansions, this analysis, although more complex, can also be applied in multiple dimensions.


We conclude that the Modified Harmonic Method (MHM) results in accurate numerical solutions to the degenerate GPME that remove the numerical artifacts associated with the FTCS scheme and harmonic averaging.  We demonstrate that harmonic averaging does not have to be cast aside.  It can be used as the base scheme, with a local switch in critical regions to the MHM.  Using the harmonic average as a base scheme that can then be modified locally with a switch can be attractive in applications, where physical blocking modeled correctly by harmonic averaging is essential.


 \appendix
 \section{Derivation of Eqn. \eqref{eq:k_eqtn}}
 \label{degenerate}
 By utilizing the Chain Rule, $k_t = k_pp_t$ and $\nabla k(p) = k_p \nabla p$.  With the use of the governing equation, we have
 \begin{equation}
 	\begin{aligned}
 		k_t &= k_p(\nabla \cdot (k(p) \nabla p)) \\
		&= (k_p|\nabla p|)^2 + kk_p \Delta p \\
		&= |\nabla k|^2 + kk_p \Delta p.
	\end{aligned}
	\label{eq:k}
 \end{equation}
Now, $\Delta k = \nabla \cdot \nabla k =  \nabla \cdot (k_p \nabla p) = k_{pp}|\nabla p|^2 + k_p \Delta p$. Hence, $k_p \Delta p = \Delta k - k_{pp}|\nabla p|^2.$ Substituting this expression into Eqn. \eqref{eq:k} gives
 \begin{equation}
			\begin{aligned}
				k_t = |\nabla k|^2  -  kk_{pp}|\nabla p|^2  + k\Delta k.
				\label{eq:k_mid}
			\end{aligned}
		\end{equation}
We have the second term of Eqn. \eqref{eq:k_eqtn} and so we must simplify the first two terms of Eqn. \eqref{eq:k_mid}.
To relate $kk_{pp}|\nabla p|^2$ to $C(p)|\nabla k|^2$ for some $C(p)$, we write $\nabla k$ in terms of its derivative with respect to $p$. Thus, $kk_{pp}|\nabla p|^2 = C(p)|\nabla k|^2 \iff kk_{pp}|\nabla p|^2 = C(p)k_p^2|\nabla p|^2 \iff kk_{pp} = C(p)k_p^2$.  Substituting this expression into Eqn. \eqref{eq:k_mid} gives Eqn. \eqref{eq:k_eqtn}.


 \section{Convergence Studies}
\label{conv}
This appendix contains the convergence results for the various problems tested with the same time step sizes as used in Section \ref{MH_Results}, where $\dt \propto \mathcal{O}(\dx^2)$.  The error norms are calculated with respect to the reference solution on the $N = 3200$ grid.  These results highlight that the MHM has the same convergence order as the methods with arithmetic and harmonic averaging.  The convergence is at best quadratic when the asymptotic region is reached.  In the more challenging cases, it is linear or superlinear.  For these problems involving shocks, very fine grids are required to reach the asymptotic region.  This is one reason why we do not consider higher-order methods for this problem class.  The higher order convergence rate would not be observed, unless extremely fine meshes are used.  Depending on the application, this may not be computationally efficient.  The global convergence order is provided.  This includes the points surrounding the shock.  The results are given in various norms, including $l_1, l_2$ and $l_\infty.$  If comparing to the convergence studies in \cite{ngo2016, zhang09}, note that high-order convergence rates are measured in the parabolic region away from the shock, rather than in the entire domain.

\begin{figure}[H]
		\center
		\begin{subfigure}[H]{.32\textwidth}  
			\includegraphics[width =\textwidth]{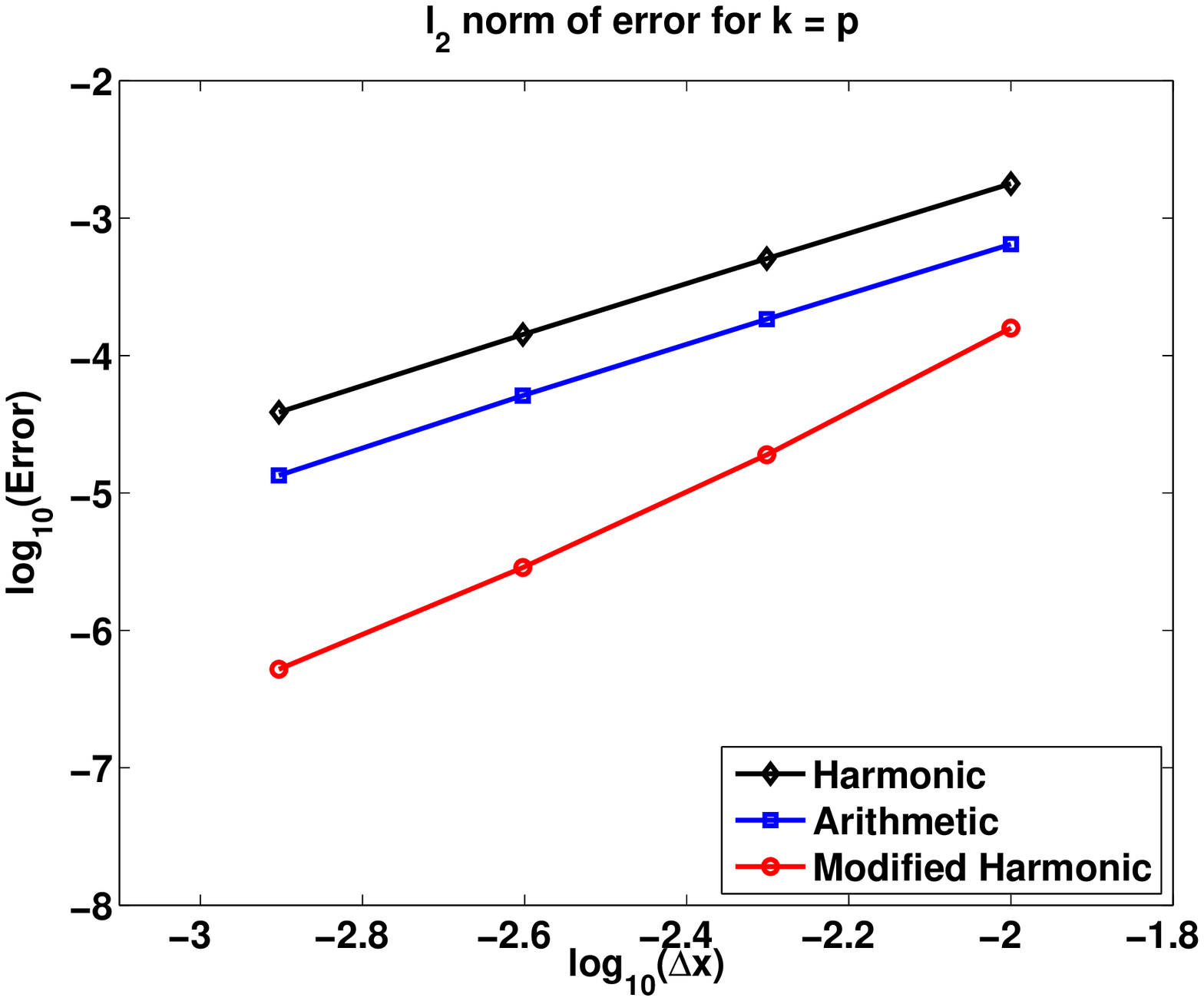}
			\caption {$k(p) = p$}
		\end{subfigure}
		\begin{subfigure}[H]{0.32\textwidth}  
			\includegraphics[width =\textwidth]{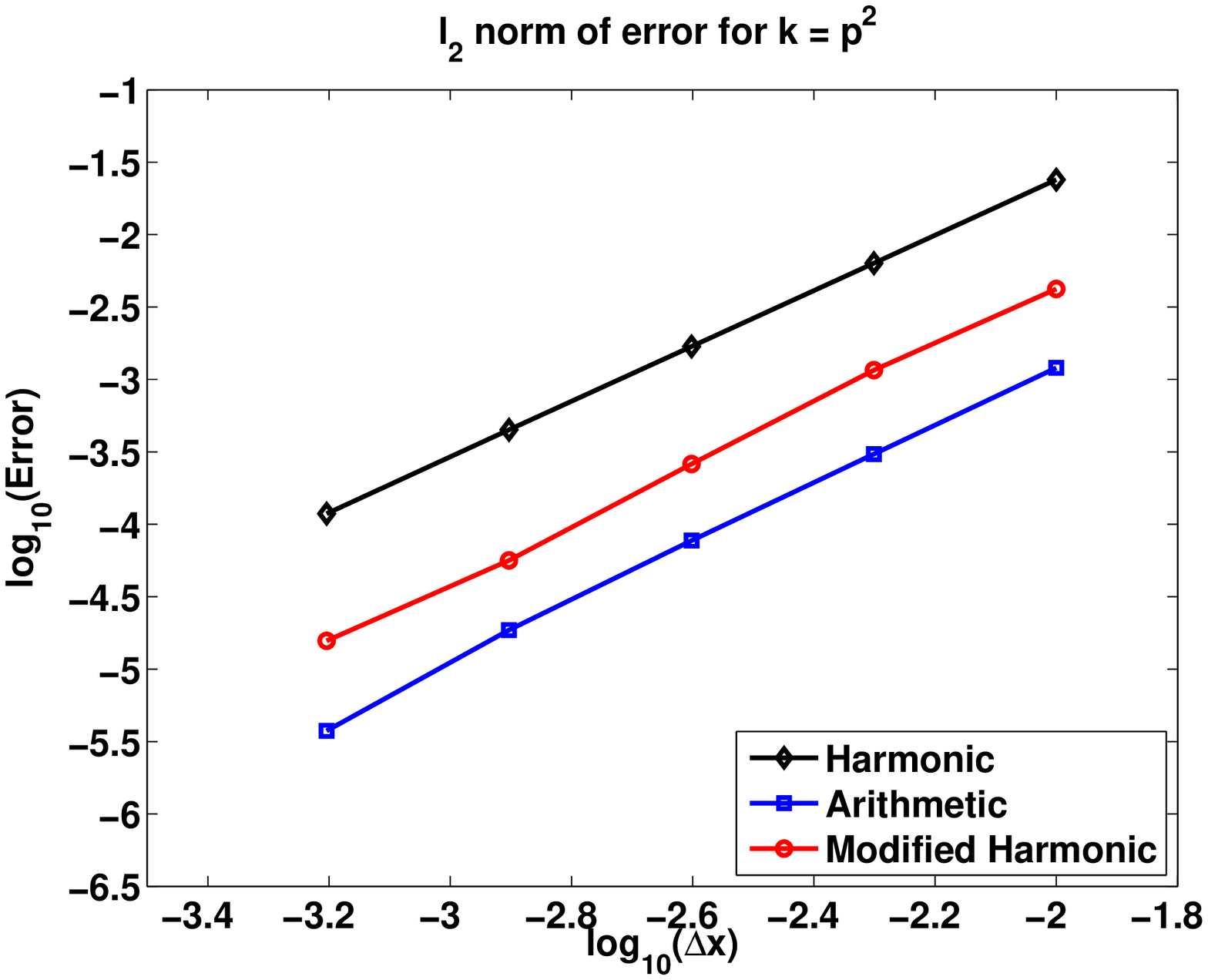}
			\caption {$k(p) = p^2$}
		\end{subfigure}
		\begin{subfigure}[H]{0.32\textwidth}  
			\includegraphics[width =\textwidth]{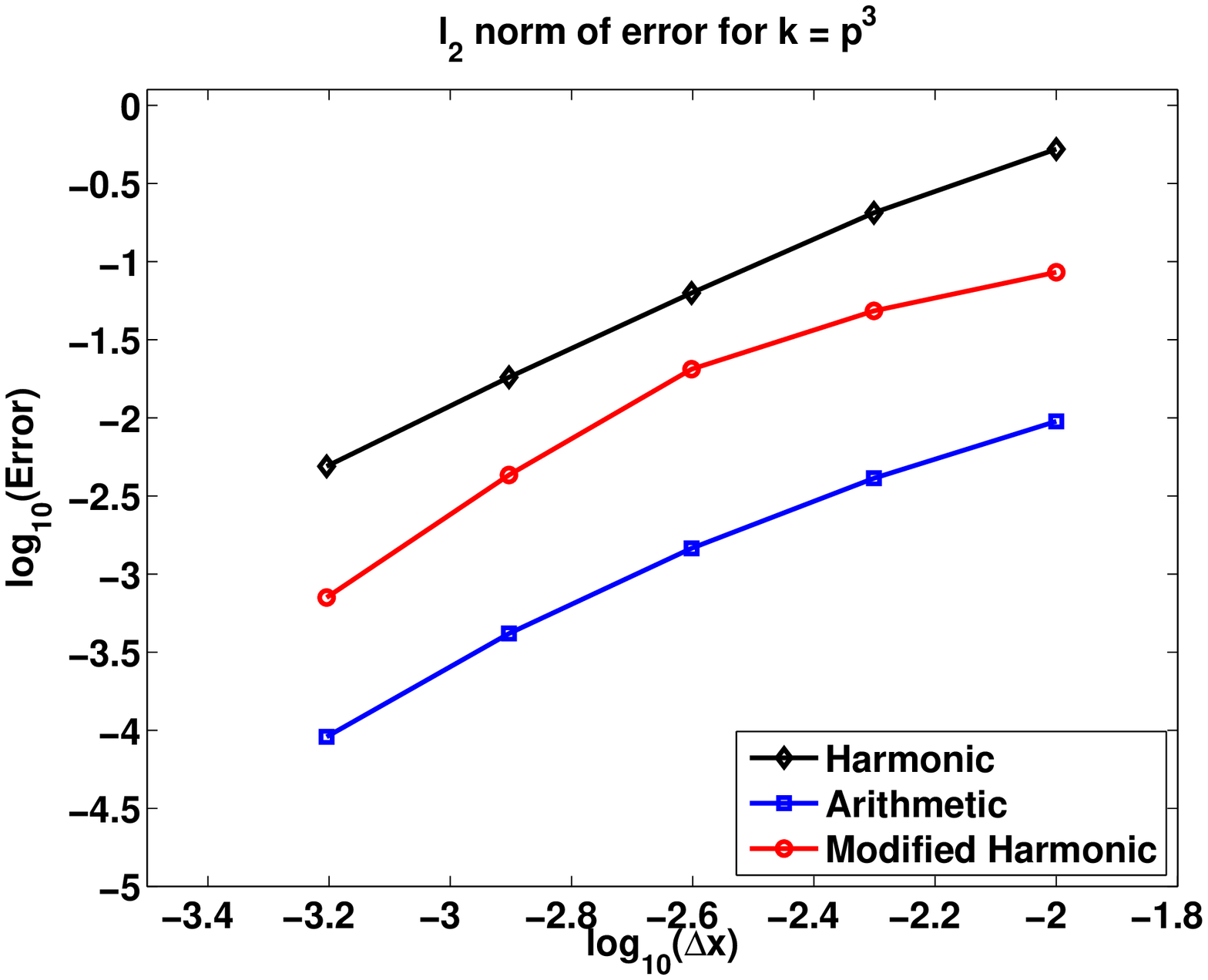}
			\caption {$k(p) = p^3$}
		\end{subfigure}
		\caption{Loglog convergence plots for the PME with various values of $m$, as measured in the $l_2$ norm.}
\end{figure}
\begin{figure}[H]
		\center
		\begin{subfigure}[H]{.32\textwidth}  
			\includegraphics[width =\textwidth]{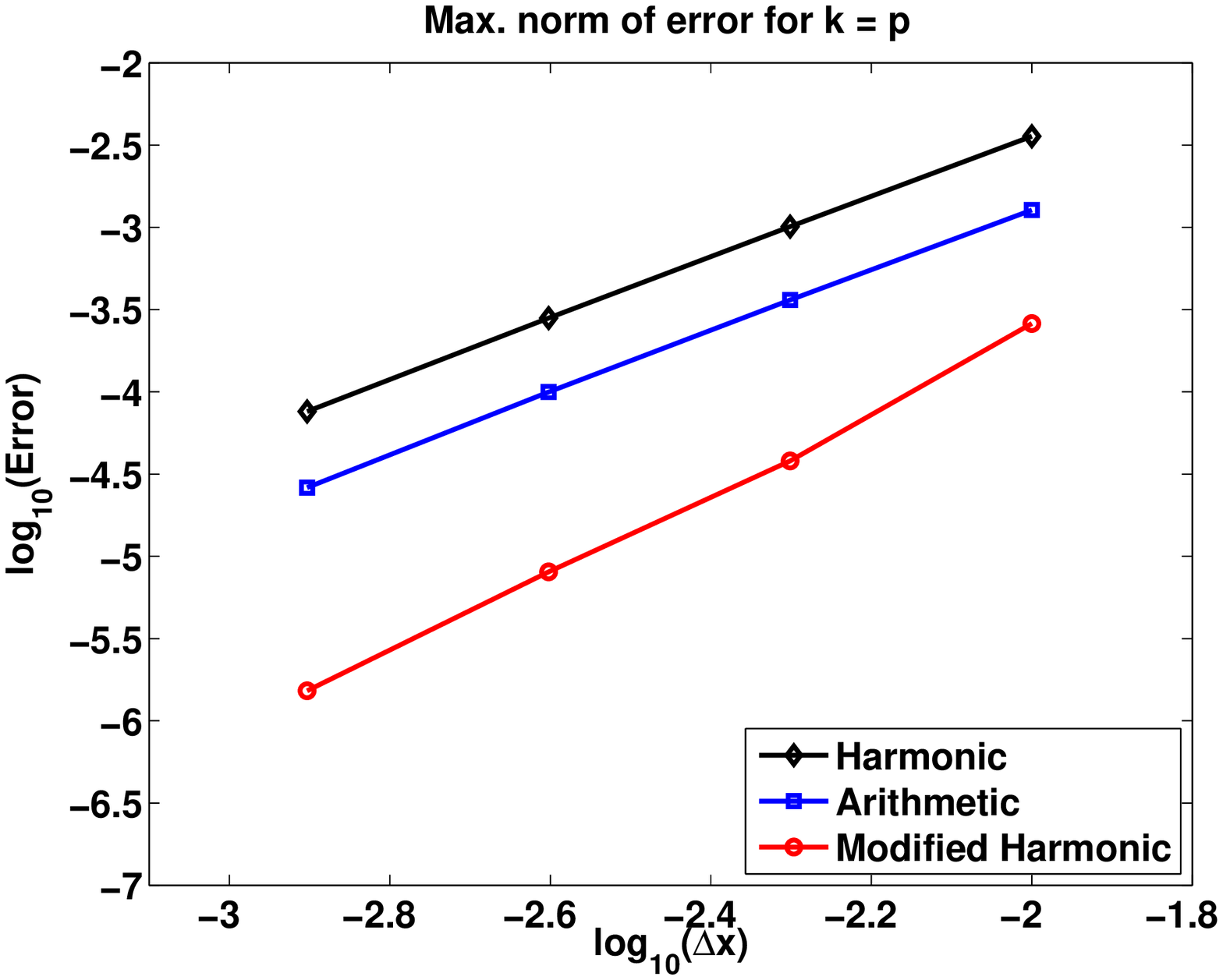}
			\caption {$k(p) = p$}
		\end{subfigure}
		\begin{subfigure}[H]{0.32\textwidth}  
			\includegraphics[width =\textwidth]{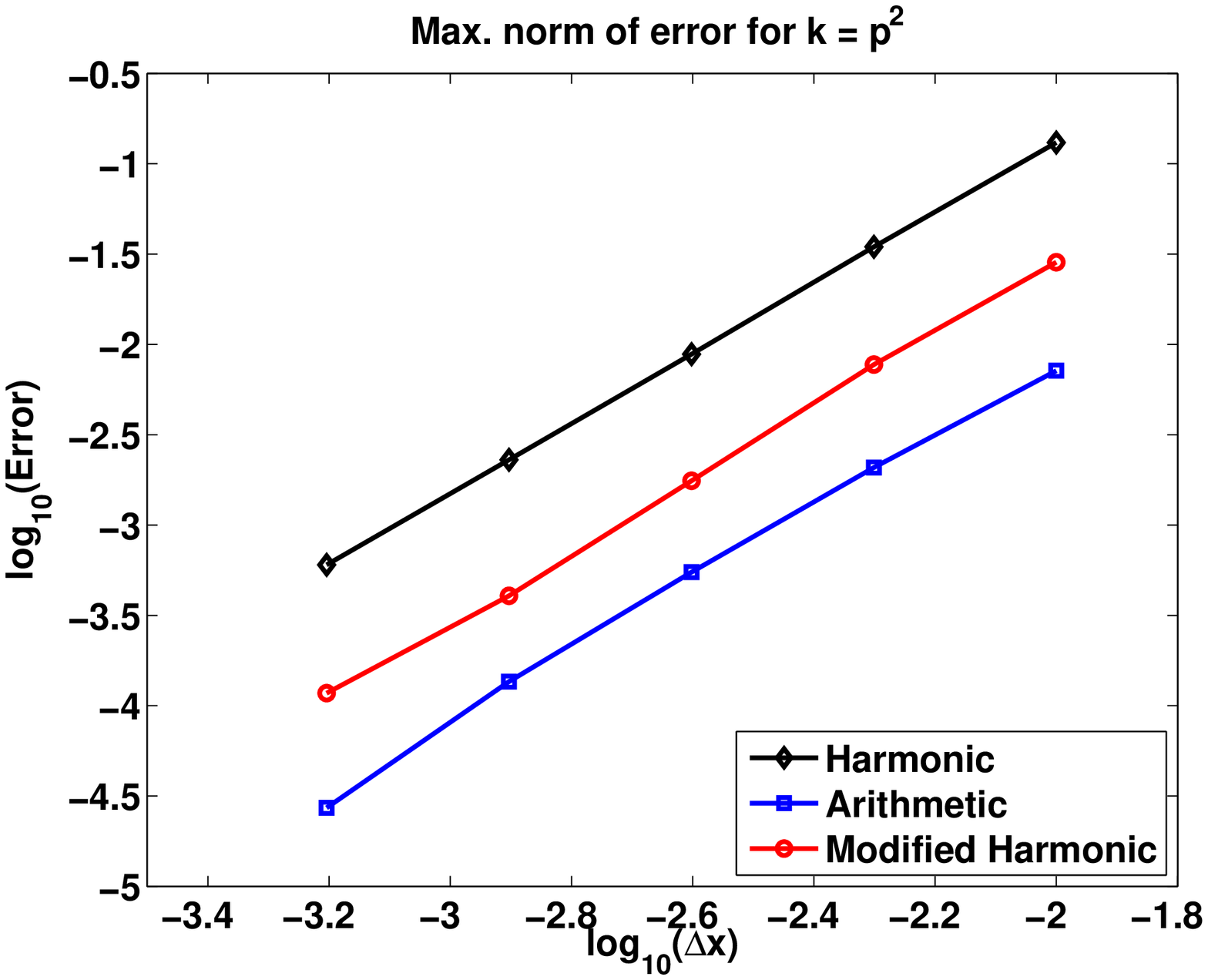}
			\caption {$k(p) = p^2$}
		\end{subfigure}
		\begin{subfigure}[H]{0.32\textwidth}  
			\includegraphics[width =\textwidth]{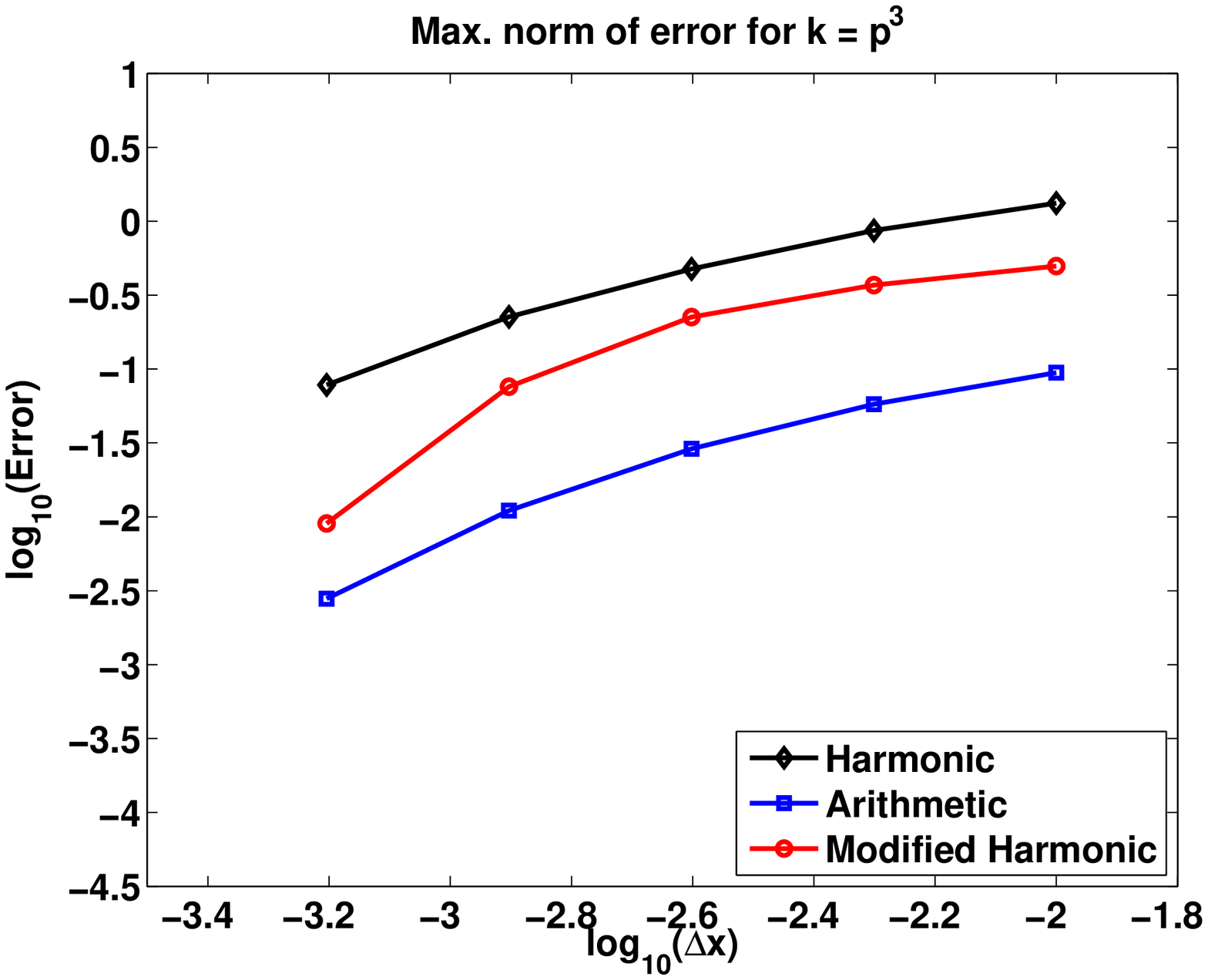}
			\caption {$k(p) = p^3$}
		\end{subfigure}
		\caption{Loglog convergence plots for the PME with various values of $m$, as measured in the $l_\infty$ norm.}
\end{figure}

\begin{table}[H]
\begin{center}
\caption{$l_2$ and $l_{\infty}$ norm errors for $k(p) = p$}
\small
\begin{tabular}{|c|c|c|c|c|c|c|}  
\hline 
& $N = 100$ & $N = 200$ & $N = 400$ &$N = 800$  &Order \\  \hline
Harmonic: $l_2$			&1.7810e-03          &5.0767e-04 &1.4199e-04 &3.8658e-05 & 1.8415 \\ \hline
Arithmetic: $l_2$       &  6.4523e-04        &1.8399e-04 &5.1068e-05 &1.3422e-05  &  1.8610     \\  \hline
Modified Harmonic: $l_2$  		& 1.5818e-04 &1.8926e-05 &2.8700e-06 & 5.2173e-07			
 &2.7454    \\ \hline
 Harmonic: $l_{\infty}$			&3.5787e-03          &1.0099e-03 &2.8041e-04 &7.5947e-05 & 1.8524 \\ \hline
Arithmetic: $l_{\infty}$	      &  1.2739e-03        &3.6128e-04 &9.9831e-05 &2.6149e-05  &  1.8674     \\  \hline
Modified Harmonic: $l_{\infty}$	  		& 2.5986e-04 &3.8020e-05 &8.0411e-06 & 1.5222e-06			
 &2.4488    \\ \hline
\end{tabular}
\label{table:timing}
\end{center}
\end{table} 

\begin{table}[H]
\begin{center}
\caption{$l_2$ and $l_{\infty}$ norm errors for $k(p) = p^2$}
\small
\begin{tabular}{|c|c|c|c|c|c|c|}  
\hline 
& $N = 100$ & $N = 200$ & $N = 400$ &$N = 800$ & $N = 1600$ &Order \\  \hline
Harmonic: $l_2$			&2.3922e-02        &6.3523e-03&1.6892e-03 &4.4868e-04 &1.1845e-04 & 1.9139 \\ \hline
Arithmetic: $l_2$	       &  1.2043e-03         &3.0568e-04 &7.7231e-05 &1.8594e-05 &  3.7444e-06 &  2.0698     \\  \hline 
Modified Harmonic: $l_2$	  		& 4.2076e-03 &1.1588e-03 &2.6080e-04 &5.6385e-05			&1.5678e-05 
 &2.0497    \\ \hline
 Harmonic: $l_{\infty}$			&1.3110e-01         &3.4657e-02 &8.8187e-03 &2.3002e-03 &6.0205e-04 &  1.9446 \\ \hline
Arithmetic: $l_{\infty}$       &  7.1516e-03     &2.0817e-03 &5.4226e-04 &1.3595e-04 &  2.7241e-05 &  2.0009    \\  \hline
Modified Harmonic: $l_{\infty}$  	&2.8494e-02	&7.7394e-03 &1.7575e-03 &4.0618e-04 & 1.1752e-04
 &2.0095    \\ \hline
\end{tabular}
\label{table:timing}
\end{center}
\end{table} 

\begin{table}[H]
\begin{center}
\caption{$l_2$ and $l_{\infty}$ norm errors for $k(p) = p^3$}
\small
\begin{tabular}{|c|c|c|c|c|c|c|}  
\hline 
& $N = 100$ & $N = 200$ & $N = 400$ &$N = 800$ & $N = 1600$ &Order \\  \hline
Harmonic: $l_2$			&5.2439e-01          &2.0508e-01 &6.3005e-02 &1.8166e-02 &4.8935e-03 & 1.6984 \\ \hline
Arithmetic: $l_2$       &  9.5041e-03      &4.0947e-03 &1.4595e-03 &4.1586e-04 &  9.0647e-05 &  1.6724     \\  \hline
Modified Harmonic: $l_2$  		&8.5436e-02 &4.8365e-02 &2.0489e-02 & 4.3061e-03			&7.0659e-04
 &1.7325    \\ \hline
 Harmonic: $l_{\infty}$			&1.3248          &8.6558e-01 &4.7567e-01 &2.2533e-01 &7.8073e-02 &  1.0111 \\ \hline
Arithmetic: $l_{\infty}$       &  9.4499e-02     &5.7749e-02 &2.881e-02 &1.1027e-02 &  2.7981e-03 &  1.2544     \\  \hline
Modified Harmonic: $l_{\infty}$  	&4.9737e-01	&3.6981e-01 &2.2481e-01 &7.6049e-02 & 9.0303e-03
 &1.3849    \\ \hline
\end{tabular}
\label{table:timing}
\end{center}
\end{table} 

For the superslow diffusion case, it requires extremely fine grids for the harmonic average to enter the asymptotic region.  So, we provide a direct comparison of the scheme with arithmetic averaging and the Modified Harmonic method.
\begin{figure}[H]
		\center
		\begin{subfigure}[H]{.32\textwidth}  
			\includegraphics[width =\textwidth]{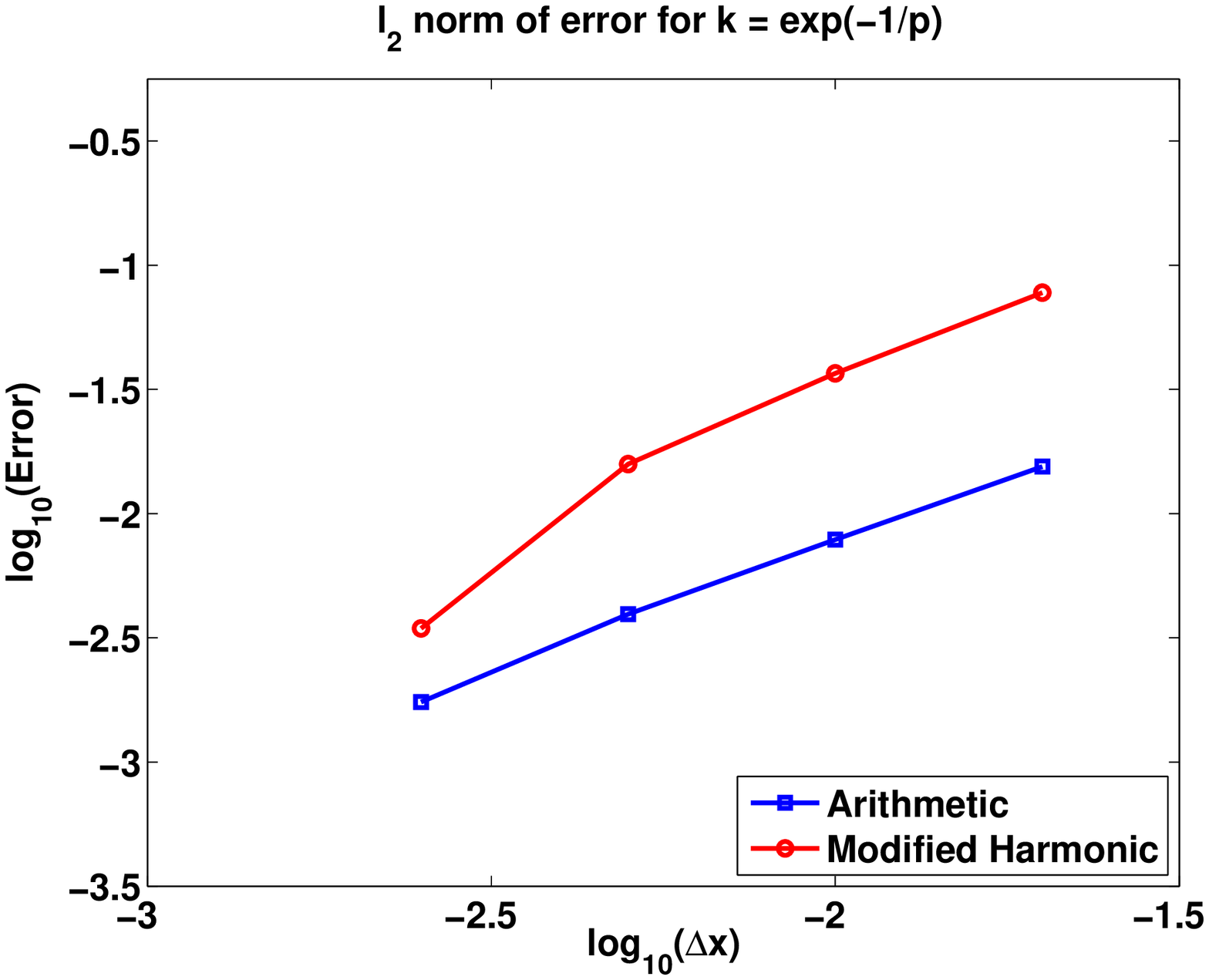}
			\caption {$l_2$ norm}
			\end{subfigure}
			\begin{subfigure}[H]{.32\textwidth}  
			\includegraphics[width =\textwidth]{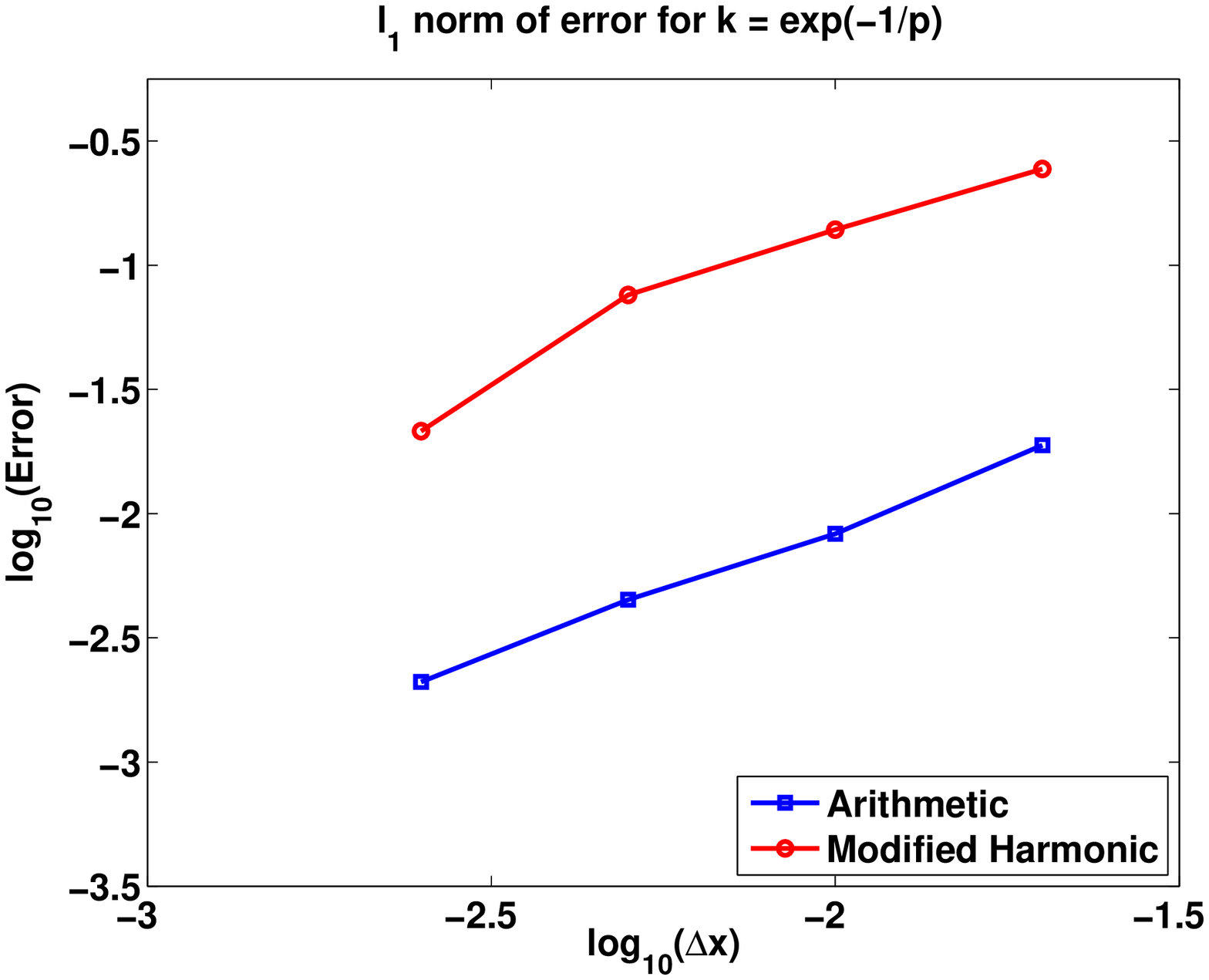}
			\caption {$l_1$ norm}
			\end{subfigure}
		\caption{Loglog convergence plots for superslow diffusion, as measured in the $l_2$ and $l_1$ norms.}
\end{figure}
\begin{table}[H]
\begin{center}
\caption{$l_2$ and $l_1$ norm errors for $k(p) = \exp(-1/p)$}
\center
\small
\begin{tabular}{|c|c|c|c|c|c|c|}  
\hline 
& $N = 50$ & $N = 100$ & $N = 200$ &$N = 400$  &Order \\  \hline
Arithmetic: $l_2$       &  1.5464e-02     &7.8655e-03 &3.9396e-03 &1.7435e-03  &  1.0444     \\  \hline
 Modified Harmonic: $l_2$  	&7.7613e-02	&3.6673e-02 &1.5824e-02 &3.4519e-03 
 &1.4685    \\ \hline
 Arithmetic: $l_1$      &   1.8870e-02     &8.2992e-03 &4.5044e-03 &2.1020e-03  &  1.0381     \\  \hline
 Modified Harmonic: $l_1$  	&2.4406e-01	&1.3910e-01 &7.5899e-02 &2.1495e-02&1.1389    \\ \hline
\end{tabular}
\label{table:timing}
\end{center}
\end{table} 

\section*{Acknowledgements}
This material is based upon work supported by the National Science Foundation Graduate Research Fellowship under Grant No. DGE - 114747.  We would also like to thank Giuseppe Domenico Cervelli for his guidance.

\bibliographystyle{elsarticle-num-names}
 \bibliography{PMEreferences}






\end{document}